\documentclass[11pt,a4paper]{article}
\pdfoutput=1
\usepackage[a4paper]{geometry}
\usepackage{graphicx}
\usepackage{amsmath}
\usepackage{natbib}

\usepackage{amsfonts}
\usepackage{enumitem} 
\usepackage{tikz}
\usetikzlibrary{arrows.meta,arrows}
\usepackage[toc,page,header]{appendix}
\usetikzlibrary{shapes,decorations,arrows,calc,arrows.meta,fit,positioning}
\usetikzlibrary{arrows.meta, automata, positioning, quotes}
\usepackage{graphicx}
\usepackage{xr}
\usepackage[caption=false]{subfig}

\usepackage[noend]{algpseudocode}
\usepackage[ruled,vlined]{algorithm2e}
\usepackage{cases}
\newlist{mylist}{enumerate*}{1}
\setlist[mylist]{label=(\roman*)}

\usepackage{amsthm}

\usepackage{etoolbox}

\makeatletter
\patchcmd{\@makecaption}
{\parbox}
{\advance\@tempdima-\fontdimen2} % decrease the width!
{}{}
\makeatother

\makeatletter
\def\BState{\State\hskip-\ALG@thistlm}
\makeatother

\newtheorem{Satz}{Satz}[section]
\newtheorem{Lemma}[Satz]{Lemma}
\newtheorem{Corollary}[Satz]{Corollary}
\newtheorem{Theorem}[Satz]{Theorem}
\newtheorem{Example}[Satz]{Example}
\newtheorem{Proposition}[Satz]{Propostion}

\theoremstyle{definition}
\newtheorem{Definition}[Satz]{Definition}

\theoremstyle{remark}

\newtheorem{Remark}[Satz]{Remark}

\newenvironment{proofof}[1][]{\begin{trivlist}
\item[\hskip \labelsep {\textit{Proof of #1.}}]}{\hfill{}$\square$\end{trivlist}}

\newcommand\ci{\protect\mathpalette{\protect\independenT}{\perp}}
\def\independenT#1#2{\mathrel{\rlap{$#1#2$}\mkern2mu{#1#2}}}

\DeclareMathOperator{\Forbb}{forb}
\DeclareMathOperator{\pa}{pa}

\DeclareMathOperator{\de}{de}
\DeclareMathOperator{\possDe}{possde}
\DeclareMathOperator{\possde}{possde}
\DeclareMathOperator{\an}{an}

\DeclareMathOperator{\possan}{possan}
\DeclareMathOperator{\CN}{cn}
\DeclareMathOperator{\possCN}{posscn}

\DeclareMathOperator{\Cov}{Cov}

\newcommand{\mpdag}{maximal PDAG}
\newcommand{\Mpdag}{Maximal PDAG}

\newcommand{\pstar}[1][p]{{#1}^{*}}
\newcommand{\g}[1][G]{\mathcal{#1}}
\newcommand{\f}[2][X,Y]{\Forbb(#1,#2)}
\newcommand{\fb}[2][X,Y]{\Forbb(\mathbf{#1},#2)}
\newcommand{\cn}[2][X,Y]{\CN(#1,#2)}
\newcommand{\cnb}[2][X,Y]{\CN(\mathbf{#1},#2)}

\newcommand{\posscn}[2][X,Y]{\possCN(#1,#2)}
\newcommand{\posscnb}[2][X,Y]{\possCN(\mathbf{#1},#2)}
\newcommand{\opt}[2][X,Y]{\mathbf{O}(#1,#2)}
\newcommand{\optb}[2][X,Y]{\mathbf{O}(\mathbf{#1},#2)}
\newcommand{\vsp}{\vspace{.2cm}}

\makeatletter
\def\cctext#1{\expandafter\@cctext\csname c@#1\endcsname}
\def\@cctext#1{\ifcase#1\or \textbf{Amenability}\or \textbf{Forbidden set}\or \textbf{Blocking}\fi}
\makeatother
\RequirePackage{hyperref}

\begin{document}
	
		\title{Graphical criteria for efficient total effect estimation via adjustment in causal linear models}
	
	\author{Leonard Henckel, \\ ETH Zurich \\ 
		\and Emilija Perkovi\'c, \\ University of Washington \\
		\and Marloes H. Maathuis, \\ ETH Zurich \\}

	\maketitle	
	
\begin{abstract}
	Covariate adjustment is a commonly used method for total causal effect estimation. In recent years, graphical criteria have been developed to identify all valid adjustment sets, that is, all covariate sets that can be used for this purpose. Different valid adjustment sets typically provide total effect estimates of varying accuracies. Restricting ourselves to causal linear models, we introduce a graphical criterion to compare the asymptotic variances provided by certain valid adjustment sets. We employ this result to develop two further graphical tools. First, we introduce a simple variance reducing pruning procedure for any given valid adjustment set. Second, we give a graphical characterization of a valid adjustment set that provides the optimal asymptotic variance among all valid adjustment sets. Our results depend only on the graphical structure and not on the specific error variances or edge coefficients of the underlying causal linear model. They can be applied to directed acyclic graphs (DAGs), completed partially directed acyclic graphs (CPDAGs) and maximally oriented partially directed acyclic graphs (maximal PDAGs). We present simulations and a real data example to support our results and show their practical applicability.
\end{abstract}

%keywords{causal inference, graphical models, covariate adjustment, efficiency}

\section{Introduction}
%Covariate adjustment is a commonly used tool to estimate total causal effects from observational data. When an appropriate set of covariates is included in a regression analysis the results can be interpreted causally even when only observational data is available. Whether a given covariate set allows for causal inference is a difficult question and
%%When the underlying causal graph is known it possible to decide whether a covariate set is valid for causal inference. 
%as a result there is a large body of work on graphical criteria that identify such valid adjustment sets. The best-known graphical criterion is the back-door criterion from \citet{pearl1993bayesian}, which is defined for directed acyclic graphs (DAGs). It has been generalized to completed partially directed acyclic graphs (CPDAGs), maximal ancestral graphs (MAGs) and partial ancestral graphs (PAGs) by \citet{maathuis2013generalized}. 
%In this paper, we use the generalized adjustment criterion from \citep{perkovic2015complete,perkovic17,perkovic16} as it is sound and complete for DAGs, CPDAGs, PAGs, MAGs and maximally oriented partially directed acyclic graphs (maximal PDAGs). This criterion was originally introduced by \citet{shpitser2010validity} for DAGs end extended to MAGs by \citet{vanconstructing}.

Covariate adjustment is a popular method for estimating total causal effects from observational data. Given a causal graph, with nodes representing covariates and edges direct effects,  graphical criteria have been developed to read of covariate sets that can be used for this purpose. We refer to such sets as valid adjustment sets. 
The best-known such criterion is probably the back-door criterion \citep{pearl1993bayesian}, which is sufficient for adjustment. A necessary and sufficient criterion was developed by \citet{shpitser2010validity} and \citet{perkovic16}. 

%We refer to the latter as the adjustment criterion. 
% vanconstructing??

%In causal linear models total effects can be estimated via ordinary least squares (OLS) regression given an appropriate conditioning set. This result is well known in the Gaussian case with \citet{shpitser2010validity} and \citet{perkovic16} having fully characterized the class of conditioning sets for which this is the case (see Definition \ref{S-adjustment}). We refer to sets fulfilling their graphical criterion as \textit{valid adjustment sets}. 

Given the complete identification of all valid adjustment sets, the following question naturally arises: If more than one valid adjustment set is available, which one should be used for estimation? 
In practice this choice will often be affected by considerations such as ease and cost of data collection. On the other hand, statistical aspects should also be taken into account as different valid adjustment sets provide estimates with varying accuracy. 
Restricting ourselves to causal linear models, we develop graphical tools to leverage the information encoded in the causal graph to identify adjustment sets that are not only valid but also efficient.

As of now, efficiency considerations have not featured prominently in adjustment set selection.
%The problem which valid adjustment set to use for estimating a total effect, has received relatively little attention so far. 
When the treatment is a single variable $X$, the parent set of $X$, i.e. the set of direct causes of $X$, is often used as an adjustment set  \citep[e.g.,][]{williamson2014introduction,gascon2015prenatal,sunyer2015association}.
% in $\g[D]$ \ema[References... Pearl and some applied articles?]. % since the parent set only depends on the local neighborhood of $X$. 
Although easy to compute and guaranteed to satisfy the backdoor-criterion, the parents of $X$ are typically quite inefficient in terms of the asymptotic variance, as they are usually strongly correlated with $X$ (see Example \ref{Apply-3.1}). 
%There has also been research showing that in the presence of unmeasured confounding, adjusting for certain direct causes of $X$ may amplify the bias \citep{wooldridge2016should,ding2017instrumental}.
Another approach to choosing an adjustment set is to adjust for as few covariates as possible \citep[e.g.,][]{de2011covariate,jonker2012iron,schliep2015effect}. For efficiency, this is also sub-optimal in general, as adjusting for certain additional covariates that explain variance in the outcome $Y$, sometimes called precision variables or risk factors, can be beneficial for efficiency (see Example \ref{Apply-3.1}). 

%There has however been some research into the role that 

%The inattention to efficiency considerations in valid adjustment set selection is reflected in the scarcity of theoretical results in this area. 
Literature on variable selection for efficient total causal effect estimation has been growing in recent years, particularly in the area of propensity score methods. 
For example there are simulation studies \citep{brookhart2006variable,lefebvre2008impact}, results regarding minimum asymptotic variance bounds \citep{rotnitzky1995semiparametric,hahn2004functional,rotnitzky2010note} and theoretical results for certain estimators \citep{robinson1991some,lunceford2004stratification,schnitzer2016variable,wooldridge2016should}. These results indicate that the following two notions hold: First, adding instrumental variables to a given valid adjustment set harms the efficiency and second, adding precision variables improves the efficiency. Model selection procedures taking these notions into account have also been developed \citep{vanderweele2011new,shortreed2017outcome}. 

%While useful and a good starting point, these two principles do suffer from some drawbacks. Firstly, they have, to our knowledge, not been comprehensively shown to hold true for any estimator. The existing theoretical results require stronger assumptions on the given covariates. Secondly, even if we suppose there correctness they are not easy to apply. Since removing and adding covariates to a valid adjustment set can render said set invalid, it is unclear whether and in what form they can be applied greedily. Further, as they rely on conditional associations, it is important to note that whether adding a given covariate to an adjustment set is actually harmful or beneficial can vary depending on the starting adjustment set. It is also not straightforward to use the two principles to compare the performance of two given valid adjustment sets, unless one is a subset of the other. 

While the above notions provide useful heuristics, there are pitfalls to the approach of labeling individual covariates as either good or bad for efficiency. Whether adding a given covariate to an adjustment set is harmful or beneficial can vary depending on the starting adjustment set, i.e. is generally speaking a conditional property. Furthermore, adding or removing a covariate might render a valid adjustment set invalid. As a result some care must be taken in sequentially applying these heuristics. 

%In particular it is not clear which valid adjustment provides a most efficient estimate or whether such a set even exists.

%Leo: decided against example as I going on a discussion about sequential algorithmic approaches would go to far

% In this paper we therefore compare the efficiency of valid adjustment sets directly, rather than considering the behavior of individual covariates. 

%Especially in the graphical framework, where it is straightforward to decide whether a given covariate set is valid, it is thus a more fruitful approach to compare the efficiency of valid adjustment sets directly, rather than considering the behavior of individual covariates. 

\citet{kuroki2003covariate} and \citet{kuroki2004selection} circumvent these difficulties by comparing the efficiency of certain pairs of valid adjustment sets, rather than considering the behavior of individual covariates. 
Both introduce graphical criteria that identify which of two valid adjustment sets provides the smaller asymptotic variance in causal linear models. The criterion from \citet{kuroki2003covariate} compares adjustment sets of size two and the criterion from \citet{kuroki2004selection} compares disjoint adjustment sets. Furthermore, both criteria require a directed acyclic graph (DAG) and a multivariate Gaussian distribution. We extend these results in various directions.

% In the presence of hidden variables, the counterparts of DAGs and CPDAGs are maximal ancestral graphs (MAGs) and partial ancestral graphs (PAGs). The back-door criterion was generalized to CPDAGs, MAGs and PAGs by  \citet{maathuis2013generalized}. 

%
%
%To the best of our knowledge, the only theoretical results for causal linear models are given by \citet[Theorem 5]{kuroki2003covariate} and \citet[Lemma 3]{kuroki2004selection}. Both provide graphical criteria in DAGs to identify which of two valid adjustment sets provides the smaller asymptotic variance. These are very useful, but somewhat limited in scope. The criterion from \citet{kuroki2003covariate} compares adjustment sets of size two and the criterion from \citet{kuroki2004selection} compares disjoint adjustment sets. Furthermore, both criteria are for DAGs and require the underlying probability distribution to be multivariate Gaussian. We extend these results in various directions. 

Our first result is a new graphical criterion (see Theorem \ref{cor:bignew12}) that can compare more pairs of valid adjustment sets than the existing criteria \citep{kuroki2003covariate,kuroki2004selection}. Our result holds for causal linear models with arbitrary error distributions as well as single and joint interventions. They can also be applied to graph types other than DAGs. We note, however, that we still cannot compare all pairs of valid adjustment sets. This is in fact impossible with the graph alone (see Example \ref{Apply-3.1}).

Building on Theorem \ref{cor:bignew12}, we introduce two further results. First, we provide a simple order invariant pruning procedure that, given a valid adjustment set, returns a subset that is also valid and provides equal or smaller asymptotic variance (see Algorithm \ref{algo: pruning} and Theorem \ref{prop:order-indep}). Our procedure is similar to that of \citet[Propositions 1 and 2]{vanderweele2011new}, who conjectured a resulting efficiency gain. Our main contribution is that we formally establish this efficiency gain for causal linear models and show the order invariance of such a procedure. 

Second, we define a valid adjustment set that provides the smallest possible asymptotic variance among all valid adjustment sets relative to $(X,Y)$ in the underlying causal graph $\g$ (see Theorem \ref{thm:optimalSetCPDAG}). 
We denote this adjustment set by $\optb{\g}$ and refer to it as \textit{ asymptotically optimal} . The fact that such an asymptotically optimal set can be defined is perhaps surprising, considering that Theorem \ref{cor:bignew12} only allows for the comparison of certain pairs of valid adjustment sets. 
%Using the results by \citet{ZanderL19} it can be shown that $\optb{\g}$ can be computed in polynomial time. 
Our results depend only on the structure of the causal graph and not on the specific edge weights or error distributions of the underlying causal linear model. 
%Except for the existence of $\optb{\g}$, our results extend to settings with hidden variables. 
We also discuss the particulars of how our results can be applied to cases with unmeasured covariates in the Discussion (Section \ref{sec:dis}).

We also provide numerical experiments to quantify the efficiency that is gained by using $\optb{\g}$ in finite samples (see Section \ref{sim-MSE}), and also apply our methods to single cell data from \citet{sachs2005causal} (see Section \ref{sec-realData}). All proofs can be found in the Supplement \citep{supplement}. We have also made our code available at \url{https://github.com/henckell/CodeEfficientVAS}.

Independent follow-up research has already expanded upon our results in various directions. 
%\citet{ZanderL19} propose a linear run time alogorithm  
In particular, \citet{rotnitzky2019efficient} show that our results on the asymptotic optimality of $\optb{\g}$ extend to a broad class of non-parametric estimators. Building on this, \citet{smucler2020efficient} consider even more general settings and construct adjustment sets with efficiency guarantees other than asymptotic optimality.
\citet{ZanderL19} provide a polynomial time algorithm to compute $\optb{\g}$. 
%$\mathcal{O}(|\mathbf{V}|+|\mathbf{E}|)$, where  $|\mathbf{V}|$  is the number of nodes and $|\mathbf{E}|$ the number of edges for the given DAG $\g=(\mathbf{V},\mathbf{E})$. 
\citet{outcomeIDA} provide an alternative characterization of $\optb{\g}$ and also integrate $\optb{\g}$ into the IDA algorithm by \citet{maathuis2009estimating}. \citet{kuipers2020variance} investigate the theoretical finite sample performance of $\opt{\g}$ in a specific non-linear example and discuss how $\opt{\g}$'s performance relates to causal discovery considerations.

\section{Preliminaries}
\label{prelim}

In this paper we use graphs where nodes represent random variables, and edges represent conditional dependencies and direct causal effects. We now give an overview of the main graphical objects used in this paper. We give the usual graphical definitions and define these objects more formally in Section \ref{graph-supp} of the Supplement. 

We consider three classes of acyclic graphs: directed acyclic graphs (DAGs), completed partially directed acyclic graphs (CPDAGs) and maximally oriented partially directed acyclic graphs (maximal PDAGs) (see Example \ref{ex: CPDAGs, PDAGs and possibly directed} in the Supplement). DAGs are directed graphs, i.e. graphs with all edges of the form $\rightarrow$ and without directed cycles. They arise naturally to describe causal relationships under the assumption of no feedback loops \citep[cf.][]{pearl2009causality}. Generally it is not possible to learn the causal DAG from observational data alone. Under the assumptions of causal sufficiency and faithfulness, one can, however, learn a Markov equivalence class of DAGs, which can be uniquely represented by a CPDAG \citep{meek1995causal,andersson1997characterization,spirtes2000causation, Chickering02}. Given explicit knowledge of some causal relationships between variables, access to interventional data, or some model restrictions, one can obtain a refinement of this class, uniquely represented by a \mpdag{} \citep{meek1995causal,tetrad1998,hoyer08,hauserBuehlmann12,eigenmann17,wang2017permutation}. 
All three graph types encode conditional independence relationships that can be read off the graph by applying the well known d-separation criterion (see Definition 1.2.3 in \cite{pearl2009causality} for DAGs, Definition 3.5 in \cite{maathuis2013generalized} for CPDAGs and Lemma \ref{lemma:dsepp1} of the Supplement for \mpdag{}s). We use the notation $\mathbf{X} \perp_{\g} \mathbf{Y}|\mathbf{Z}$ to denote that $\mathbf{Z}$  \textit{d-separates} $\mathbf{X}$ from $\mathbf{Y}$ in $\g$, with $\mathbf{X},\mathbf{Y}$ and $\mathbf{Z}$ pairwise disjoint nodes sets in a graph $\g$.

\begin{Remark}
	DAGs and CPDAGs are special cases of \mpdag{}s. In the remainder of the paper results are generally stated in terms of \mpdag{}s. Readers unfamiliar with CPDAGs and \mpdag{}s may also disregard this and simply think of all results as being with respect to DAGs. 
\end{Remark}

%\noindent{}\textbf{\Mpdag{}s.} {\color{blue-violet} A PDAG $\g$ is \textit{maximally oriented} (\mpdag{}) if and only if the graphs in Figure \ref{fig:orientationRules} are \textbf{not} induced subgraphs of $\g$.}
%%A directed edge $X \rightarrow Y$ in a \mpdag{} $\g$ corresponds to $X \rightarrow Y$ in every DAG represented by $\g$. 
%DAGs and CPDAGs are special cases of \mpdag{}s. In general, a \mpdag{} $\g$ describes a subset of a Markov equivalence class of DAGs, denoted by $[\g]$.  
%A \mpdag{} $\g$ has the same adjacencies and v-structures as any DAG in $[\g]$. 
%Moreover, a directed edge $X\to Y$ in $\g$ corresponds to a directed edge $X \to Y$ in every
%DAG in $[\g]$, and for any undirected edge $X -Y$  in $\g$, $[\g]$ contains a DAG with $X \rightarrow Y$ and a DAG with $X \leftarrow Y$.

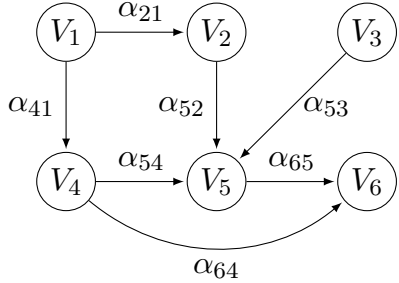
\begin{figure}
	\centering
	\begin{tikzpicture}[>=stealth',->,>=latex,shorten >=1pt,auto,node distance=1.8cm,scale=1.1,transform shape]
	\tikzstyle{vertex}=[circle, draw, inner sep=2pt, minimum size=0.15cm]
	\node[vertex]         (V1)                        {$V_1$};
	\node[vertex]         (V2) [below of= V1]  		{$V_4$};
	\node[vertex]         (V3) [right of = V2] 		{$V_5$};
	\node[vertex]       	 (V4)  [right of= V3] 		{$V_6$};
	\node[vertex]       	 (V5)  [above of= V3] 		{$V_2$};
	\node[vertex]       	 (V6)  [above of= V4] 		{$V_3$};
	
	\draw[->] 	(V1) 	edge  node[left]  {$\alpha_{41}$} (V2);
	\draw[->] 	(V2) 	edge  node[above]  {$\alpha_{54}$}  (V3);
	\draw[->] 	(V3) 	edge   node[above]  {$\alpha_{65}$} (V4);
	\draw[->] 	(V5) 	edge   node[left]  {$\alpha_{52}$} (V3);
	%		\draw[->] 	(V5) 	edge   node[right]  {$\alpha_{53}$} (V3);
	\draw[->] 	(V2) [bend right=40] 	edge   node[below]  {$\alpha_{64}$} (V4);
	\draw[->] 	(V1)	edge   node[above]  {$\alpha_{21}$} (V5);
	\draw[->] 	(V6)	edge   node[right]  {$\alpha_{53}$} (V3);
	
	\end{tikzpicture}
	\caption{DAG  from Examples \ref{examle: linear SEM}, \ref{ex:total effect} and \ref{example-O}}
	%	and \ref{emma:cpdag}.}
	
	\label{figure:super simple}
\end{figure}

%\subsection{Causal graphs and the generalized adjustment criterion}\label{sec: causal graphs and adjustment}

%Throughout, nodes in a graph represent random variables. At times, we use $\mathbf{X}$ to denote both a set of nodes in a graph and the corresponding vector of random variables. The graphical preliminaries are given in the section \ref{S-graph-supp} of the supplement.

%\subsection{Causal liner models and total effects}

We now introduce causal linear models, total effects and defines some notation.

\vsp\noindent\textbf{Causal DAGs, CPDAGs, \mpdag{}s.} 
We consider interventions $do(\mathbf{x})$ (for $\mathbf{X}\subseteq \mathbf{V}$), which represent outside interventions that set $\mathbf{X}$ to $\mathbf{x}$ uniformly for the entire population \citep{pearl1995causal}.
A density $f$  of $\mathbf{V} = \{V_1,\dots,V_p\}$ is \textit{compatible with a causal DAG} $\g =(\mathbf{V},\mathbf{E})$ if all post-intervention densities $f(\mathbf{v}|do(\mathbf{x}))$ factorize as:
\begin{equation}
f(\mathbf{v}|do(\mathbf{x}))=
\begin{cases}
\prod_{V_{i} \in \mathbf{V} \setminus \mathbf{X}}f(v_{i}|\pa(V_{i},\g)), &  \text{if }\mathbf{X} =\mathbf{x}, \\
0, & \text{otherwise.}
\end{cases}
\label{eq11}
\end{equation}
Equation \eqref{eq11} is known as the truncated factorization formula \citep{pearl2009causality}, manipulated density formula \citep{spirtes2000causation} or the g-formula \citep{robins1986new}. 
A density $f$  of $\mathbf{V} = \{V_1,\dots,V_p\}$ is \textit{compatible with a causal \mpdag{} or a causal CPDAG} $\g$ if it is compatible with a causal DAG $\g[D] \in [\g]$, where $[\g]$ is the class of DAGs represented by $\g$.

\vsp\noindent\textbf{Causal linear model.} 
Let $\g = (\mathbf{V,E})$ be a DAG. 
Then $\mathbf{V} =(V_1, \dots, V_p)^T$, $p\ge1$ follows a \textit{causal linear model} compatible with $\g$ if the following two conditions hold:
%Then the distribution $f$ of $\mathbf{V} =(X_1, \dots, X_p)^T$, $p\ge1$ follows a \textit{causal linear model} compatible with $\g$ if the following two conditions hold:
\begin{enumerate}
	\item The distribution $f$ of $\mathbf{V}$ is compatible with the causal DAG $\g$.
	\item $V_1, \dots, V_p$ follows a set of linear equations
	\begin{align}
	V_i = \sum_{V_j \in \pa(V_i,\g)} \alpha_{ij}V_j +\epsilon_{v_i}, \label{sem}
	\end{align}
	where $i,j \in \{1,\dots,p\}, i \neq j$,  $\alpha_{ij} \in \mathbb{R}$ and $\epsilon_{v_1}, \dots, \epsilon_{v_p}$ are jointly independent random variables with mean 0 and finite variance. 
\end{enumerate}
$\mathbf{V}$ follows a \textit{causal linear model} compatible with a \mpdag{} or CPDAG $\g$, if it follows a causal linear model compatible with a DAG $\g[D] \in [\g]$.

We refer to $\epsilon_{v_1}, \dots, \epsilon_{v_p}$ as errors and emphasize that we do not require them to be Gaussian. Furthermore, by  construction $E[\mathbf{V}] = \boldsymbol{0}$. The coefficient $\alpha_{ij}$ corresponding to the edge $V_j \rightarrow V_i$ in the causal DAG $\g$ can be interpreted as the \textit{direct effect} of $V_j$ on $V_i$ with respect to $\mathbf{V}$. 

%Note that due to the underlying DAG being causal, the equations \eqref{sem} describe non-symmetric causal mechanisms. 

\begin{Example}
	\label{examle: linear SEM}
	Consider the DAG $\g=(\mathbf{V},\mathbf{E})$ in Figure \ref{figure:super simple}. Then the generating mechanism 
	\begin{align*}
	&V_1 \leftarrow \epsilon_{v_1} \\
	&V_2 \leftarrow \alpha_{21}V_1 +\epsilon_{v_2} \\
	&V_3 \leftarrow \epsilon_{v_3}\\
	&V_4 \leftarrow\alpha_{41}V_1 + \epsilon_{v_4}\\
	&V_5 \leftarrow\alpha_{52}V_2 + \alpha_{53}V_3 + \alpha_{54}V_4 + \epsilon_{v_5}\\
	&V_6 \leftarrow\alpha_{64}V_4 + \alpha_{65}V_5 + \epsilon_{v_6}
	\end{align*}
	with $\epsilon_{v_1} \sim \mathrm{Uniform}(-1,1),\epsilon_{v_2} \sim \mathrm{Uniform}(-2,2),\epsilon_{v_3} \sim \mathcal{N}(0,1),\epsilon_{v_	4} \sim \mathcal{N}(0,2),\epsilon_{v_5} \sim \mathrm{Uniform}(-3,3)$ and $\epsilon_{v_6} \sim \mathcal{N}(0,3)$ for $\mathbf{V}$ is an example of a causal linear model compatible with $\g$. We use the notation $\leftarrow$ to emphasize that we are considering a generating mechanism and not just an equation.
	
	A do intervention, for example $do(V_4=1)$, then corresponds to replacing the generating mechanism of the intervened on variables with the fixed intervention value, e.g. $V_4 \leftarrow 1$.
	%	An example of a causal linear model compatible with $\g=(\mathbf{V},\mathbf{E})$
	
\end{Example}

%\begin{Remark}
%	Note that, different from the linear structural equation model setting, we do not require the errors from the generating causal equations to be identical between different do-intervention environments. 
%\end{Remark}

%Further, any causal structural equation model gives rise to a unique compatible DAG.

\noindent \textbf{Causal and proper paths} Let $\g$ be a \mpdag{}. A path $(V_1,\dots, V_m)$ in $\g$ is called a \textit{causal path} from $V_1$ to $V_m$ if $V_i \rightarrow V_{i+1}$ for every $i \in \{1,\dots, m-1\}$. 
%It is called a \textit{possibly causal path} from $V_1$ to $V_m$ if for no $i,j \in \{1,\dots, m\},i<j$, such that $V_i \leftarrow V_j$. 
Let $\mathbf{X}$ and $\mathbf{Y}$ be disjoint node sets in a causal \mpdag{} $\g$. A path from $\mathbf{X}$ to $\mathbf{Y}$ is \textit{proper} if only the first node on $p$ is in $\mathbf{X}$

\vsp\noindent\textbf{Total effects}. \citep[][]{pearl2009causality}
Let $\mathbf{V} = (\mathbf{X}^T, \mathbf{Y}^T, \mathbf{Z}^T)^T$ be a random vector, where $\mathbf{X} = (X_1, \dots, X_{k_x})^T$ and $\mathbf{Y}=(Y_1,
\dots, Y_{k_y})^T$. Then the total effect of $\mathbf{X}$ on $\mathbf{Y}$ is defined as the matrix $\tau_{\mathbf{yx}}$, where 
\[
(\tau_{\mathbf{yx}})_{j,i} = \frac{\partial}{\partial x_i} E[Y_j|do(x_1,\dots,x_{k_x})], i \in \{1,\dots,k_x\},i \in \{1,\dots,k_y\}
\]
represents the effect of $X_i$ on $Y_j$ in the joint intervention of $\mathbf{X}$ on $\mathbf{Y}$. 
In general, $\tau_{\mathbf{yx}}$ is a matrix of functions, but 
in causal linear models the partial derivatives do not depend on $x_i$. Hence, $\tau_{\mathbf{yx}}$ reduces to a matrix of numbers, whose values are determined by the coefficients in Equation \eqref{sem} \citep{wright1934method,nandy2017estimating}. We can thus give an equivalent definition of the total effect specific to this setting. 
Consider disjoint node sets $\mathbf{X}=\{X_1,\dots,X_{k_x}\}$ and $\mathbf{Y}=\{Y_1,\dots,Y_{k_y}\}$ in a causal DAG $\g=(\mathbf{V},\mathbf{E})$, such that $\mathbf{V}$ follows a causal linear model compatible with $\g$. 
The \textit{total effect along a causal path} $p$ from $X$ to $Y$ in $\g$ is the product of the edge coefficients along p. %The total effect of $X$ on $Y$, $\tau_{yx}$, is equal to the sum of total effects along all causal paths from $X$ to $Y$. 
The \textit{total effect} of $\mathbf{X}$ on $\mathbf{Y}$ is then the matrix $\tau_{\mathbf{yx}} \in \mathbb{R}^{k_y \times k_x}$ whose $(j,i)$-th value $(\tau_{\mathbf{yx}})_{j,i}$ is equal to the sum of the total effects along all proper causal paths from $\mathbf{X}$ to $Y_j$ starting with $X_i$ in $\g$. 

If $\mathbf{V}$ follows a causal linear model compatible with a causal CPDAG or \mpdag{} $\g$, the total effect of $\mathbf{X}$ on $\mathbf{Y}$ is identifiable if it is the same for every DAG in $[\g]$. 

%By Theorem \ref{S-thm:adjustment-set}, a sufficient (but not necessary) condition for identifiability, is the existence of a valid adjustment set relative to $(\mathbf{X},\mathbf{Y})$ in $\g$. 

\begin{Remark}
	Consider the total effect $\tau_{\mathbf{yx}}$ of $\mathbf{X}$ on $\mathbf{Y}$. If for some $Y_j \in \mathbf{Y}$ and $X_i \in \mathbf{X}$, $Y_j$ is a non-descendant of $X_i$ then $(\tau_{\mathbf{yx}})_{j,i} = 0$.
	Further, the total effect $\tau_{y_jx_i}$ of $X_i$ on $Y_j$ will generally differ from the partial total effect $(\tau_{\mathbf{yx}})_{j,i}$ in the joint intervention on $\mathbf{X}$. This is due to the latter effect not considering causal paths from $X_i$ to $Y_j$ that contain other nodes in $\mathbf{X}\setminus\{X_i\}$. The total effect of $\mathbf{X}$ on any $Y_j$, however, does not depend on the remaining $\mathbf{Y} \setminus \{Y_j\}$.
	\label{remark: joint vs non}
\end{Remark}
%\vsp
%  Finally, we introduce the following notation.

%\noindent\textbf{D-separation:} 
%(Cf.\  Definition 1.2.3 in \cite{pearl2009causality}, Definition 3.5 in \cite{maathuis2013generalized} and Lemma \ref{S-lemma:dsepp1} from the Supplement) 

\vsp
%We use the convention that for any two disjoint node sets $\mathbf{X}$ and $\mathbf{Y}$ it holds that $\emptyset \perp_{\g} \mathbf{X} | \mathbf{Y}$.
%In a causal PDAG this implies that $\mathbf{X} \ci \mathbf{Y}|\mathbf{Z}$
%If  $\mathbf{Z}$ does not block every definite status path between any node in $\mathbf{X}$ and any node in $\mathbf{Y}$ in $\g$, we write $\mathbf{X} \not\perp_{\g} \mathbf{Y}|\mathbf{Z}$. 

\noindent\textbf{Notation for covariance matrices and regression coefficients.}
Consider random vectors $\mathbf{S}, \mathbf{T}$, and $\mathbf{W_1}, \mathbf{W_2}, \dots, \mathbf{W_m}$, and let $\mathbf{W}=(\mathbf{W_1}^T, \dots, \mathbf{W_m}^T)^T$, $k_s=|\mathbf{S}|$ and $k_t=|\mathbf{T}|$. 
We denote the covariance matrix of $\mathbf{S}$ with $\Sigma_{\mathbf{ss}} \in\mathbb{R}^{k_s \times k_s}$ and the covariance matrix between $\mathbf{S}$ and $\mathbf{T}$ with $\Sigma_{\mathbf{st}}  \in \mathbb{R}^{k_s \times k_t}$, where its $(i,j)$th element equals $\Cov(S_i,T_j)$. We further define 
$\Sigma_{\mathbf{ss}.\mathbf{t}} = \Sigma_{\mathbf{ss}} - \Sigma_{\mathbf{st}} \Sigma^{-1}_{\mathbf{tt}} \Sigma_{\mathbf{ts}}$. If $|\mathbf{S}|=1$, we write $\sigma_{ss.\mathbf{t}}$ instead. 
Let $\boldsymbol{\beta}_{\mathbf{s}\mathbf{t}.\mathbf{w}} \in \mathbb{R}^{k_s \times k_t}$ represent the least squares regression coefficient matrix whose $(i,j)$-th element is the regression coefficient of $T_j$ in the regression of $S_i$ on $\mathbf{T}$ and $ \mathbf{W}$, with $\boldsymbol{\hat{\beta}}_{\mathbf{st}.\mathbf{w}}$ denoting the corresponding estimator.  We also use the notation that $\boldsymbol{\beta}_{\mathbf{st}.{\mathbf{w_1w_2}\cdots \mathbf{w_m}}} = \boldsymbol{\beta}_{\mathbf{st}.{\mathbf{w}}}$ and $\Sigma_{\mathbf{s}\mathbf{t}.{\mathbf{w_1w_2}\cdots \mathbf{w_m}}} = \Sigma_{\mathbf{s}\mathbf{t}.\mathbf{w}} $. Given a set $\mathbf{X}=\{X_1,\dots,X_k\}$ we use the notation $\mathbf{X}_{-i}$ to denote $\mathbf{X} \setminus \{X_i\}$.

\section{Main Results} \label{main}

\subsection{Total effect estimation via covariate adjustment}

In causal linear models total effects can be estimated via OLS regression given an appropriate adjustment set. This result is well known in the Gaussian case with \citet{shpitser2010validity} and \citet{perkovic16} having fully characterized the class of valid adjustment sets (see Definition \ref{adjustment}). 

The fact that total effects can be estimated via OLS regression has been shown to generalize to causal linear models with arbitrary error distributions for a singleton $X$ with the adjustment set $\pa(X,\g)$ \citep[Proposition 3.1 from the supplement of][]{nandy2017estimating}. We now extend this property to arbitrary valid adjustment sets and derive the estimator's asymptotic distribution. \\

\begin{Proposition}
	Let $\mathbf{X}=\{X_1,\dots,X_{k_x}\}$ and $\mathbf{Y}=\{Y_1,\dots,Y_{k_y}\}$ be disjoint node sets in a causal DAG $\g=(\mathbf{V},\mathbf{E})$ and let $\mathbf{V}$ follow a causal linear model compatible with $\g$. Let $\mathbf{Z}$ be a valid adjustment set relative to $(\mathbf{X},\mathbf{Y})$ in $\g$. Then 
	\begin{align}
	\sqrt{n} ((\boldsymbol{\hat{\beta}}_{\mathbf{yx}.\mathbf{z}})_{j,i} - (\tau_{\mathbf{yx}})_{j,i}) \xrightarrow{d} \mathcal{N}(0, \frac{\sigma_{y_jy_j.\mathbf{xz}}}{\sigma_{x_ix_i.\mathbf{x}_{-i}\mathbf{z}}}), \nonumber 
	\end{align}
	%	\begin{equation}
	%	a.var(\hat{\tau}^{\mathbf{z}}_{\mathbf{yx}})_{j,i} = a.var(\hat{\beta}_{y_jx_i.\mathbf{x}_{-i}\mathbf{z}}) = \frac{\sigma_{y_jy_j.\mathbf{xz}}}{\sigma_{x_ix_i.\mathbf{x}_{-i}\mathbf{z}}}.
	%	\label{avar.equation}
	%	\end{equation}
	for all $i = 1,\dots,k_x $ and $j = 1, \dots,k_y$, with $\xrightarrow{d}$ denoting convergence in distribution.
	
	\label{avarDAG}
\end{Proposition}

%The proof of Proposition \ref{avarDAG} is given in Section \ref{S-sec:avar} of the Supplement. 
The key aspect to Proposition \ref{avarDAG}, is that it does not require the considered regression of $\mathbf{Y}$ on $\mathbf{X}$ and $\mathbf{Z}$ to be well-specified, in the sense of being linear and having homoskedastic residuals. One may think this generality is not needed, given that we consider causal linear models. However, in a causal linear model with non-Gaussian errors, adjusted regressions other than that of a node on its parents, are not generally well-specified (see Example \ref{section:example}). 
%We make this explicit in Example \ref{S-section:example}. 
We note that for Proposition \ref{avarDAG} to hold for misspecified regresssions, it is essential that $\mathbf{Z}$ is a valid adjustment set. Of course, Proposition \ref{avarDAG} corresponds to what we know for well-specified regressions, in which case the restriction to valid adjustment sets is not needed. For causal linear models with Gaussian errors all regressions are well-specified.

%The asymptotic variance is identical to the well known OLS asymptotic variance formula for the multivariate Gaussian setting or the general homoskedastic residuals case. 
%Nonetheless, this is \emph{not} a trivial result. In a causal linear model with non-Gaussian errors, misspecified OLS regressions, i.e. regressions with non-linearity or heteroskedastic residuals, naturally arise and nonetheless Proposition \ref{avarDAG} holds. We illustrate this behaviour in Example \ref{S-section:example}.

Due to the result in Proposition \ref{avarDAG}, we use the notation $\hat{\tau}^{\mathbf{z}}_{\mathbf{yx}}$ to denote the least squares estimate $\boldsymbol{\hat{\beta}}_{\mathbf{y}\mathbf{x}.\mathbf{z}}$ of $\tau_{\mathbf{y}\mathbf{x}}$, for any valid adjustment set $\mathbf{Z}$ relative to $(\mathbf{X},\mathbf{Y})$. We also write 
\begin{equation}
a.var(\hat{\tau}^{\mathbf{z}}_{\mathbf{yx}})_{j,i} = a.var(\hat{\beta}_{y_jx_i.\mathbf{x}_{-i}\mathbf{z}}) = \frac{\sigma_{y_jy_j.\mathbf{xz}}}{\sigma_{x_ix_i.\mathbf{x}_{-i}\mathbf{z}}}
\label{eq:avar}
\end{equation}
and $a.var(\hat{\tau}^{\mathbf{z}}_{\mathbf{yx}})$ to denote the matrix with entries
\[
a.var(\hat{\tau}^{\mathbf{z}}_{\mathbf{yx}})_{j,i} = a.var((\hat{\tau}^{\mathbf{z}}_{\mathbf{yx}})_{j,i}), \quad i=1, \dots, k_x \ \textit{and} \ j=1, \dots, k_y.
\]

\begin{Remark}
	The terms in Equation \eqref{eq:avar} depend on the distribution of $\mathbf{V}=\{V_1,\dots,V_p\}$ only through the covariance matrix $\Sigma_{\mathbf{vv}}$, which in turn only depends on the underlying causal linear model through the edge coefficients $\alpha_{ij}$ and error variances $\mathrm{var}(\epsilon_{v_i}), i,j \in \{1,\dots,p\}, i \neq j$ \citep[cf.][]{nandy2017estimating}. 
	In particular, this implies that the asymptotic variance $a.var(\hat{\tau}^{\mathbf{z}}_{\mathbf{yx}})$ does not depend on the error
	%		$\epsilon_{v_i},i\in \{1,\dots,p\}$ 
	distribution families.
	\label{remark: var enough}
\end{Remark}

\begin{Example}
	Consider the causal DAG $\g=(\mathbf{V},\mathbf{E})$ in Figure~\ref{figure:super simple} and assume that $\mathbf{V}$ follows the causal linear model from Example \ref{examle: linear SEM}.
	% with edge coefficients $\alpha_{13}$, $\alpha_{21}$ and $\alpha_{23}$ 
	%be the edge coefficient corresponding to $X\leftarrow Y$, $X \rightarrow Z$ and $Y \rightarrow Z$, respectively. 
	The total effect of $V_4$ on $V_6$ in $\g$ is $\tau_{64} = \alpha_{64} + \alpha_{65}\alpha_{54}$. 
	By Proposition \ref{avarDAG}, $\tau_{64}$ also equals the  population level regression coefficient of $V_4$ in the regression of $V_6$ on $V_4$ and any adjustment set of the form $\mathbf{A} \cup \mathbf{B}$, with $\mathbf{A} \subseteq \{V_1,V_2\}$ non-empty and $\mathbf{B} \subseteq \{V_3\}$ possibly empty. (see Definition \ref{adjustment}).
	
	\label{ex:total effect}
\end{Example}

%%%%%%%%%%%%%%%%%%%%%%%%%%%%%%%%%%%%%%%%%%%%%%%%%%%%%%%%%%%%%%%%%%
%%%%%%%%%%%%%%%%%%%%%%%%%%%%%%%%%%%%%%%%%%%%%%%%%%%%%%%%%%%%%%%%%%

%%%%%%%%%%%%%%%%%%%%%%%%%%%%%%%%%%%%%%%%%%%%%%%%%%%%%%%%%%%%%%%%%

\subsection{Comparing valid adjustment sets}

\vsp We now introduce a new graphical criterion for qualitative comparisons between the asymptotic variances resulting from certain pairs of valid adjustment sets, which is more general than the criteria of \citet{kuroki2003covariate} and \cite{kuroki2004selection}.

\begin{Theorem}
	Let $\mathbf{X}$ and $\mathbf{Y}$ be disjoint node sets in a \mpdag{} $\g = (\mathbf{V},\mathbf{E})$, such that $\mathbf{V}$ follows a causal linear model that is compatible with $\g$. Let $\mathbf{Z_1}$ and $\mathbf{Z_2}$ be two valid adjustment sets relative to $(\mathbf{X,Y})$ in $\g$ and let $\mathbf{T}=\mathbf{Z_1} \setminus \mathbf{Z_2}$ and $\mathbf{S}=\mathbf{Z_2} \setminus \mathbf{Z_1}$. 
	If $\mathbf{Y} \perp_{\g} \mathbf{T} | \mathbf{X} \cup \mathbf{Z_2} $ and $\mathbf{X} \perp_{\g} \mathbf{S} | \mathbf{Z_1}$,
	then 
	\[
	a.var(\hat{\tau}_{\mathbf{yx}}^{\mathbf{z_2}}) \leq 
	a.var(\hat{\tau}_{\mathbf{yx}}^{\mathbf{z_1}}),
	\]
	\label{cor:bignew12}
	with the matrix inequality denoting entry wise inequality.
	
\end{Theorem}

The proof of Theorem \ref{cor:bignew12} relies on equation \eqref{eq:avar}. The intuition behind it is that the more  information a conditioning set $\mathbf{B}$ contains on a target variable $A$ the smaller $\sigma_{aa.\mathbf{b}}$. 
%Thus, additional information on $\mathbf{Y}$ decreases the asymptotic variance, while additional information on $\mathbf{X}$ increases it. 
Thus, the assumed conditional independence statements imply that 
\begin{align}
\sigma_{x_ix_i.\mathbf{x_{-i}\mathbf{z_1}}} \leq \sigma_{x_ix_i.\mathbf{x_{-i}\mathbf{z_2}}} \ \textrm{and} \
\sigma_{y_jy_j.\mathbf{xz_2}} \leq \sigma_{y_jy_j.\mathbf{xz_1}}, \nonumber 
\end{align}
for all $X_i \in \mathbf{X}$ and $Y_j \in \mathbf{Y}$.
%$\mathbf{Z_2} \cup \mathbf{X}$ contains more information on $\mathbf{Y}$ than $\mathbf{Z_1} \cup \mathbf{X}$, while $\mathbf{Z_1}$ contains more information on $\mathbf{X}$ than $\mathbf{Z_2}$. 

We stress that when a causal linear model with non-Gaussian errors is considered, Theorem \ref{cor:bignew12} holds only for pairs of \emph{valid adjustment sets}. This is due to Proposition \ref{avarDAG} only holding for misspecified regressions when a valid adjustment set is considered.

\begin{figure}
	\centering
	\subfloat[]{
		\centering
		\begin{tikzpicture}[scale=0.12]
		\tikzstyle{every node}+=[inner sep=0pt]
		\draw [black] (28.9,-14.1) circle (3);
		\draw (28.9,-14.1) node {$A$};
		\draw [black] (41.7,-14.1) circle (3);
		\draw (41.7,-14.1) node {$B$};
		\draw [black] (54.4,-14.1) circle (3);
		\draw (54.4,-14.1) node {$C$};
		\draw [black] (35.2,-23.5) circle (3);
		\draw (35.2,-23.5) node {$X$};
		\draw [black] (47.7,-23.5) circle (3);
		\draw (47.7,-23.5) node {$Y$};
		\draw [black] (21.9,-23.5) circle (3);
		\draw (21.9,-23.5) node {$D$};
		\draw [black] (38.4,-23.5) -- (44.7,-23.5);
		\fill [black] (44.7,-23.5) -- (43.9,-23) -- (43.9,-24);
		\draw [black] (32.2,-23.5) -- (24.9,-23.5);
		\fill [black] (24.9,-23.5) -- (25.7,-24) -- (25.7,-23);
		\draw [black] (30.6,-16.9) -- (33.48,-21.04);
		\fill [black] (33.48,-21.04) -- (33.44,-20.1) -- (32.62,-20.67);
		\draw [black] (39.8,-16.5) -- (36.85,-20.99);
		\fill [black] (36.85,-20.99) -- (37.7,-20.6) -- (36.87,-20.05);
		\draw [black] (43.5,-16.6) -- (46.14,-20.94);
		\fill [black] (46.14,-20.94) -- (46.15,-19.99) -- (45.3,-20.51);
		\draw [black] (52.7,-16.7) -- (49.48,-21.08);
		\fill [black] (49.48,-21.08) -- (50.35,-20.73) -- (49.55,-20.14);
		\end{tikzpicture}
	}
	\unskip
	\hspace{1cm}
	\subfloat[]{
		\centering
		\begin{tikzpicture}[scale=0.12]
		\tikzstyle{every node}+=[inner sep=0pt]
		\draw [black] (28.9,-14.1) circle (3);
		\draw (28.9,-14.1) node {$A$};
		\draw [black] (41.7,-14.1) circle (3);
		\draw (41.7,-14.1) node {$B$};
		\draw [black] (54.4,-14.1) circle (3);
		\draw (54.4,-14.1) node {$C$};
		\draw [black] (35.2,-23.5) circle (3);
		\draw (35.2,-23.5) node {$X$};
		\draw [black] (47.7,-23.5) circle (3);
		\draw (47.7,-23.5) node {$Y$};
		\draw [black] (21.9,-23.5) circle (3);
		\draw (21.9,-23.5) node {$D$};
		\draw [black] (38.4,-23.5) -- (44.7,-23.5);
		\fill [black] (44.7,-23.5) -- (43.9,-23) -- (43.9,-24);
		\draw [black] (24.9,-23.5) -- (32.2,-23.5);
		\fill [black] (32.2,-23.5) -- (31.4,-23) -- (31.4,-24);
		\draw [black] (30.57,-16.59) -- (33.53,-21.01);
		\fill [black] (33.53,-21.01) -- (33.5,-20.07) -- (32.67,-20.62);
		\draw [black] (32,-14.1) -- (38.7,-14.1);
		\fill [black] (38.7,-14.1) -- (37.9,-13.6) -- (37.9,-14.6);
		\draw [black] (52.6,-16.6) -- (49.44,-21.05);
		\fill [black] (49.44,-21.05) -- (50.31,-20.69) -- (49.49,-20.11);
		\draw [black] (51.4,-14.1) -- (44.7,-14.1);
		\fill [black] (44.7,-14.1) -- (45.5,-14.6) -- (45.5,-13.6);
		\end{tikzpicture}
	}
	\caption{(a) Causal DAG from Examples \ref{Apply-3.1}, \ref{ex:pruning} and \ref{example-O}, (b) causal DAG from Examples \ref{hidden} and \ref{example-O}.}
	\label{newgraph}
\end{figure}
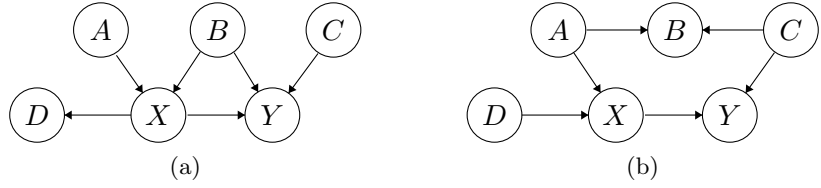

\begin{Example}
	Consider the DAG $\g = (\mathbf{V},\mathbf{E})$ in Figure \ref{newgraph}(a) and assume that $\mathbf{V}$ follows a causal linear model compatible with $\g$. Any valid adjustment set $\mathbf{Z}$ relative to $(X,Y)$ in $\g$ is of the form $\{B\} \cup \mathbf{S}$, where $\mathbf{S} \subseteq \{A,C,D\}$ (see Definition \ref{adjustment}).
	Fix any such set $\mathbf{Z}$. One can check that i) $Y \perp_{\g} A | X \cup \mathbf{Z} \setminus \{A\}$, ii) $Y \perp_{\g} D | X \cup \mathbf{Z} \setminus \{D\}$  and iii) $X \perp_{\g} C | \mathbf{Z} \setminus \{C\}$. 
	We can thus apply Theorem \ref{cor:bignew12} to the following pairs of valid adjustment sets:
	\begin{enumerate}
		\item  $\mathbf{Z_1} = \mathbf{Z}$ and $\mathbf{Z_2}=\mathbf{Z} \setminus \{A\}$,
		\item  $\mathbf{Z_1} = \mathbf{Z}$ and $\mathbf{Z_2}=\mathbf{Z} \setminus \{D\}$ and
		\item $\mathbf{Z_1}=\mathbf{Z} \setminus \{C\}$ and $\mathbf{Z_2} = \mathbf{Z}$.
	\end{enumerate} 
	We can conclude that adding $A$ or $D$ to any conditioning set worsens the asymptotic variance, while a converse statement holds for $C$. Consequently, $\{B,C\}$ provides the best asymptotic variance, while the set $pa(X,\g)=\{A,B\}$ does not fare well. 
	
	\begin{table}
		\caption{Asymptotic variances for six randomly drawn causal linear models compatible with the DAG in Figure \ref{newgraph}(a) \label{asyvar}}
		\centering
		\begin{tabular}{r|rrrrrrrr}
			Adjustment set & Case 1 & Case 2 & Case 3 & Case 4 & Case 5 & Case 6 \\ 
			\hline
			$\{A,B\}$ & 5.38 & 5.47 & 0.85 & 0.57 & 5.42 & 0.64  \\ 
			$\{A,B,C\}$ & 1.44 & 4.44 & 0.51 & 0.39 & 2.02 & 0.61  \\ 
			$\{B\}$ & 3.49 & 4.40 & 0.54 & 0.26 & 2.76 & 0.39  \\ 
			$\{B,C\}$ & 0.94 & 3.58 & 0.32 & 0.18 & 1.03 & 0.37  \\ 
			$\{A,B,D\}$ & 7.20 & 7.39 & 12.65 & 0.65 & 5.72 & 0.69  \\ 
			$\{A,B,C,D\}$ & 1.93 & 6.01 & 7.59 & 0.45 & 2.13 & 0.65  \\ 
			$\{B,D\}$ & 5.31 & 6.33 & 12.34 & 0.35 & 3.05 & 0.44  \\ 
			$\{B,C,D\}$ & 1.42 & 5.15 & 7.41 & 0.24 & 1.14 & 0.41  \\ 
		\end{tabular}
	\end{table}

	In order to empirically verify these results, we randomly drew six causal linear models compatible with $\g$ and computed the asymptotic variances $a.var(\hat{\tau}^\mathbf{z}_{yx})$ for each valid adjustment set $\mathbf{Z}$ relative to $(X,Y)$ in $\g$, across these 6 models. Specifically, we did the following for each model. We drew error variances $\sigma_{vv}$ for each node $V \in \mathbf{V}$ independently from a standard uniform distribution and edge coefficients $\alpha_{vw}$ for each edge $(W,V) \in \mathbf{E}$ independently from a standard normal distribution. From these parameters we computed the causal linear model's covariance matrix and then, in accordance with Proposition \ref{avarDAG}, the asymptotic variances corresponding to each valid adjustment set. We did not consider error properties other than the variance (and mean 0) as they are irrelevant for the asymptotic variances (see Remark \ref{remark: var enough}).
	
	The thus obtained asymptotic variances are given in Table \ref{asyvar}. They show that the three proven trends do in fact hold and that $\{B,C\}$ provides the best asymptotic variance in the considered models. Interestingly, the order of the asymptotic variances corresponding to any two sets that cannot be compared using Theorem \ref{cor:bignew12}, such as $\{A,B,C\}$ and $\{B\}$, or $\{A,B\}$ and $\{B,D\}$, are in fact inconsistent throughout the considered models.
	\label{Apply-3.1}

\end{Example}

We now give two simple corollaries of Theorem \ref{cor:bignew12}. The first one shows that superfluous parents of $X$ are harmful for the asymptotic variance, while the second one shows that parents of $Y$ are beneficial.

\begin{Corollary}
	Let $\mathbf{X}$ and $\mathbf{Y}$ be disjoint node sets in a \mpdag{} $\g = (\mathbf{V},\mathbf{E})$ and let $\mathbf{V}$ follow a causal linear model compatible with $\g$. Let $\mathbf{Z}$ be a valid adjustment set relative to $(\mathbf{X},\mathbf{Y})$ in $\g$ and let $P \in \pa(\mathbf{X},\g)$. If $\mathbf{Z'} = \mathbf{Z} \setminus \{P\}$ is a valid adjustment set relative to $(\mathbf{X},\mathbf{Y})$ in $\g$, then 
	\[
	a.var(\hat{\tau}_{\mathbf{yx}}^{\mathbf{z'}}) \leq 
	a.var(\hat{\tau}_{\mathbf{yx}}^{\mathbf{z}}).
	\]
	\label{cor:badparents}
\end{Corollary}

\begin{Corollary}
	Let $\mathbf{X}$ and $\mathbf{Y}$ be disjoint node sets in a \mpdag{} $\g = (\mathbf{V},\mathbf{E})$ and let $\mathbf{V}$ follow a causal linear model compatible with $\g$. Let $\mathbf{Z}$ be a valid adjustment set relative to $(\mathbf{X},\mathbf{Y})$ in $\g$ and let $R \in 
	\pa(\mathbf{Y},\g)$. If $\mathbf{Z'} = \mathbf{Z} \cup \{R\}$ is a valid adjustment set relative to $(\mathbf{X},\mathbf{Y})$ in $\g$, then
	\[
	a.var(\hat{\tau}_{\mathbf{yx}}^{\mathbf{z'}}) \leq 
	a.var(\hat{\tau}_{\mathbf{yx}}^{\mathbf{z}}).
	\]
	\label{cor:goodparents}
\end{Corollary}

We now give a third corollary of Theorem \ref{cor:bignew12}, especially relevant for randomized trials, where $\pa(\mathbf{X},\g)=\emptyset$. It shows that, when restricting oneself to covariates not in $\de(\mathbf{X},\g)$, enlarging an adjustment set can only be beneficial for the asymptotic variance.  In particular this implies that adjusting for additional pre-treatment covariates in a randomized trial can only be beneficial for the asymptotic variance. 

%In particular, this results can be applied  is particularly relevant for randomized trials, where it can be used even when the causal graph is not known.

% which arises in the context of randomized trials. In this special case

\begin{Corollary}
	Let $\mathbf{X}$ and $\mathbf{Y}$ be disjoint node sets in a DAG $\g = (\mathbf{V},\mathbf{E})$, such that $\pa(\mathbf{X},\g)=\emptyset$ and let $\mathbf{V}$ follow a causal linear model compatible with $\g$. Let $\mathbf{Z}$  and $\mathbf{Z}'$ be two node sets in $\g$, such that $\mathbf{Z} \cap (\de(\mathbf{X},\g)\cup\mathbf{Y}) = \emptyset $ and $\mathbf{Z'} \cap (\de(\mathbf{X},\g)\cup\mathbf{Y}) = \emptyset$. Then $\mathbf{Z}$ and $\mathbf{Z'}$ are valid adjustment sets relative to $(\mathbf{X},\mathbf{Y})$ in $\g$ and if $\mathbf{Z} \subseteq \mathbf{Z}'$,
	\[
	a.var(\hat{\tau}_{\mathbf{yx}}^{\mathbf{z'}}) \leq 
	a.var(\hat{\tau}_{\mathbf{yx}}^{\mathbf{z}}).
	\]
	\label{cor:pretreatment}
\end{Corollary}

\subsection{Pruning procedure}

\label{sec:prune}

The result from Theorem~\ref{cor:bignew12} can be used to prune a valid adjustment set to obtain a subset that is still valid and yields a smaller asymptotic variance. 
Generally, which of the subsets $\mathbf{\tilde{Z}} \subseteq \mathbf{Z}$ provides the optimal asymptotic variance depends on the edge coefficients in the underlying causal linear model (see Example \ref{hidden}).
% as shown in Example \ref{hidden} where both $\{Z_1, Z_2\}$ and the empty set can be the asymptotically optimal subset, depending on the edge coefficients. 
However, we can use Theorem \ref{cor:bignew12} to identify a subset such that there is no other subset for which Theorem \ref{cor:bignew12} guarantees a better asymptotic variance. 
This is formalized in Algorithm \ref{Algorithm} whose soundness is stated in Theorem \ref{prop:order-indep}.

In practice, such pruning is advisable as it reduces the number of variables that need to be measured while also improving precision. Although similar pruning procedures exist \citep{hahn2004functional,vanderweele2011new}, we believe Theorem \ref{cor:bignew12} and Theorem \ref{prop:order-indep} to be the first theoretical guarantees for such pruning in causal linear models. 

\begin{algorithm}[h!]
	\label{algo: pruning}
	\LinesNumbered 
	\SetKwData{Set}{Set}
	\SetKwFunction{Foreach}{Foreach}
	\SetKwInOut{Input}{input}\SetKwInOut{Output}{output}
	\Input{Causal \mpdag{} $\g$ and disjoint node sets $\mathbf{X},\mathbf{Y}$ and $\mathbf{Z}$ in $\g$, such that $\mathbf{Z}$ is  a valid adjustment set relative to $(\mathbf{X},\mathbf{Y})$ in $\g$}
	\Output{Valid adjustment set $\mathbf{Z'} \subseteq \mathbf{Z}$ relative to $(\mathbf{X},\mathbf{Y})$ in $\g$, such that $a.var(\hat{\tau}^{\mathbf{z'}}_{\mathbf{yx}}) \leq a.var(\hat{\tau}^{\mathbf{z}}_{\mathbf{yx}})$} 		
	\BlankLine
	\Begin{
		$\mathbf{Z'} = \mathbf{Z}$; \\
		\ForEach{ $Z \in \mathbf{Z'}$}{
			\If{ $\mathbf{Y} \perp_{\g} Z | \mathbf{X} \cup (\mathbf{Z'} \setminus \{Z\})$}
			{ $\mathbf{Z'} = \mathbf{Z'} \setminus \{Z\}$;}} 
		\Return{$\mathbf{Z'}$;}
		\caption{{\bf Pruning procedure} \label{Algorithm}}}
	
\end{algorithm}

%Pruning $\mathbf{Z}$ into $\mathbf{Z'}$ according to Theorem \ref{cor:bignew12} involves removing nodes from $\mathbf{Z}$ that are d-separated from $\mathbf{Y}$ given $\mathbf{X}$ and the remaining nodes in $\mathbf{Z}$. 

%	We now give a simple and efficient algorithm for pruning a given valid adjustment set to obtain a possibly asymptotically optimal subset that is still a valid adjustment set. 
%It follows from Lemma \ref{S-lemma:order-indep} in the Supplement.

\begin{Theorem}\label{prop:order-indep}
	Let $\mathbf{X}$ and $\mathbf{Y}$ be disjoint node sets in a \mpdag{} $\g = (\mathbf{V},\mathbf{E})$ and let $\mathbf{V}$ follow a causal linear model compatible with $\g$. Let $\mathbf{Z}$ be a valid adjustment set relative to $(\mathbf{X,Y})$ in $\g$. Applying Algorithm \ref{Algorithm} then yields a valid adjustment set $\mathbf{Z'} \subseteq \mathbf{Z}$, such that $a.var(\hat{\tau}^{\mathbf{z'}}_{\mathbf{yx}}) \leq a.var(\hat{\tau}^{\mathbf{z}}_{\mathbf{yx}})$ and there is no other subset of $\mathbf{Z}$ for which Theorem \ref{cor:bignew12} guarantees a better asymptotic variance than $\mathbf{Z'}$. Further, Algorithm \ref{Algorithm} outputs the same set $\mathbf{Z'}$, regardless of the order in which the nodes in $\mathbf{Z}$ are considered. 
	%and any intermediate set of the pruning algorithm is also a valid adjustment set with respect to $(\mathbf{X,Y})$ in $\g$.  
\end{Theorem}

\begin{Example}
	
	We now return to Example \ref{Apply-3.1} and the DAG $\g = (\mathbf{V},\mathbf{E})$ in Figure \ref{newgraph}(a) to illustrate Algorithm \ref{Algorithm}. Fix some valid adjustment $\mathbf{Z}$ relative to $(X,Y)$ in $\g$. As i) $Y \perp_{\g} A | X \cup \mathbf{Z} \setminus \{A\}$, ii) $Y \perp_{\g} D | X \cup \mathbf{Z} \setminus \{D\}$  and iii) $X \perp_{\g} C | \mathbf{Z} \setminus \{C\}$, Algorithm \ref{Algorithm} will discard the nodes $A$ and $D$, while keeping the nodes $B$ and $C$ whenever these nodes are in $\mathbf{Z}$. This is done independently of the order in which the nodes are considered. Hence, $\mathbf{Z}$ will either be pruned to $\{B\}$ or $\{B,C\}$. Both these sets are valid adjustment sets relative to $(X,Y)$ in $\g$ and $\{B,C\}$ yields the optimal asymptotic of all valid adjustment set, while $\{B\}$ yields the optimal asymptotic variance of all valid adjustment sets that do not contain $C$.
	
	\label{ex:pruning}
	
\end{Example}

\begin{Example}
	\label{hidden}

	We now give an example in which one cannot use Theorem \ref{cor:bignew12} to decide which subset $\mathbf{\tilde{Z}} \subseteq \mathbf{Z}$ of a valid adjustment set $\mathbf{Z}$ provides the optimal asymptotic variance. Instead, the optimal subset depends on the edge coefficients and error variances of the underlying causal linear model.
	
	Consider the DAG $\g$ in Figure \ref{newgraph}(b) and two sets of possible edge coefficients for $\g$. Let all edge coefficients that are not explicitly mentioned be $1$ and let $\alpha_{ba} = 0.5, \alpha_{xa} = 0.25$ and $\alpha_{yx} = 2$ in Case i), while $\alpha_{xa} = 0.7$ and $\alpha_{yc}=0.5$ in Case ii). With all error variances equal to 1 in both cases, one obtains the asymptotic variances shown in Table \ref{non-exhaustive}, where we ignore error properties other than variance (and mean 0) in accordance with Remark \ref{remark: var enough}.

	%							\begin{table}
	%			\caption[]{Non-exhaustive table of the asymptotic variances corresponding to valid adjustment sets in the two cases considered in Example \ref{hidden} \label{non-exhaustive}}
	%			\centering
	%			\begin{tabular}{c | c  c}
	%				Adjustment set & Case i) & Case ii) \\
	%				\hline	
	%				$\{C\}$ & 0.48 & 0.4 \\
	%				$\{B,C\}$ & 0.49 & 0.45 \\
	%				$\{A,C\}$ & 0.5 & 0.5 \\
	%				$\emptyset$ & 0.97 & 0.5 \\
	%				$\{A,B\}$ & 0.75 & 0.56 \\
	%				$\{A\}$ & 1 & 0.62 \\
	%			\end{tabular}
	%		\end{table}
	
	\begin{table}
		\caption[]{Non-exhaustive table of the asymptotic variances corresponding to valid adjustment sets in Example \ref{hidden} \label{non-exhaustive}}
		\centering
		\begin{tabular}{c | c c c c c c}
			Adjustment set & $\{C\}$ & $\{B,C\}$ & $\{A,C\}$ & $\emptyset$ & $\{A,B\}$ & $\{A\}$  \\
			\hline
			Case i) &  0.48 & 0.49 & 0.5 & 0.97 & 0.75 &  1  \\ 
			Case ii) & 0.4 & 0.45 & 0.5 & 0.5 & 0.56 & 0.62\\
			%		\hline	
			%		 0.48 & 0.4 \\
			%		 0.49 & 0.45 \\
			%		 0.5 & 0.5 \\
			%		0.97 & 0.5 \\
			%		 0.75 & 0.56 \\
			%		 1 & 0.62 \\
		\end{tabular}
	\end{table}
	% in the two cases considered

	The set $\{C\}$ provides the smallest asymptotic variance in both cases and will also be the output of Algorithm \ref{Algorithm} applied to any valid adjustment set containing $C$. If we instead consider valid adjustment sets that do not contain $C$ the situation is more complex. If, for example, we apply Algorithm \ref{Algorithm} to $\{A,B,D\}$, the output is $\{A,B\}$, which is the subset that yields the optimal asymptotic variance in Case i), but is bested by the empty set in Case ii). These two sets cannot be compared with Theorem \ref{cor:bignew12}. However, Theorem \ref{cor:bignew12} still implies that the valid adjustment sets $\{A\},\{D\},\{A,D\}$ and $\{A,B,D\}$ provide worse asymptotic variances than both the empty set and $\{A,B\}$. Algorithm \ref{Algorithm} will prune these sets to either $\{A,B\}$ or the empty set, depending on whether they originally included $\{A,B\}$.

	%	
	%	%Hence, even when $\opt{\g}$ cannot be used, Theorem \ref{cor:bignew12} can still be used to make a well informed choice when picking an alternative valid adjustment set. 
	
	%	%	Note that this counter example extends to PAGs, as in the PAG corresponding to Figure \ref{Mbias} with $S_3$ marginalized out  $\{S_1\},\{S_1,S_2\}$ and the empty set are still valid adjustment sets relative to $(X,Y)$ and obviously provide the same asymptotic variances.
\end{Example}

\subsection{The optimal adjustment set} \label{sec:cpdag}

\label{sec:optimal-set}
We will now define a set that provides the optimal asymptotic variance among all valid adjustment sets. 
This is remarkable since Theorem \ref{cor:bignew12} can only compare the asymptotic variance provided by certain valid adjustment sets (see Examples \ref{Apply-3.1} and \ref{hidden}). Nevertheless, it allows us to define this optimal set, whose optimality only depends on the underlying causal graph. We first give some preparatory definitions which for simplicity we restrict to the DAG setting. The general definitions for \mpdag{}s are given in Section \ref{graph-supp}. This includes the definition of the set of possible descendants $\possde(\mathbf{X},\g)$, which in the case that $\g$ is a DAG, reduces to the set of descendants $\de(\mathbf{X},\g)$.

\noindent \textbf{Causal and forbidden nodes} 
Consider a DAG $\g$ and two disjoint node sets $\mathbf{X},\mathbf{Y}$ in $\g$. We define \textit{causal nodes} relative to $(\mathbf{X,Y})$ in $\g$, denoted $\cnb{\g}$, as all nodes on proper causal paths from $\mathbf{X}$ to $\mathbf{Y}$, excluding nodes in $\mathbf{X}$. 
For singleton $X$, causal nodes are also called mediating nodes. 
%	Analogously, we define
%	\textit{possible causal nodes} relative to $(\mathbf{X,Y})$ in $\g$, denoted $\posscnb{\g}$, as all nodes on proper possibly causal paths from $\mathbf{X}$ to $\mathbf{Y}$, excluding nodes in $\mathbf{X}$.
%	 (see Section \ref{S-graph-supp} in the Supplement \citep{supplement}). 
% \textit{Possible causal nodes} relative to $(\mathbf{X}$, $\mathbf{Y})$ in $\g$ are all nodes on proper possibly causal paths from $\mathbf{X}$ to $\mathbf{Y}$, excluding nodes in $\mathbf{X}$.
% Causal and possibly causal nodes relative to $(\mathbf{X,Y})$ and $\g$ are denoted by $\cnb{\g}$ and $\posscnb{\g}$ respectively.
We then define the \textit{forbidden nodes} relative to $(\mathbf{X},\mathbf{Y})$ in $\g$ as
\[
\fb{\g} = \de(\cnb{\g}, \g) \cup \mathbf{X}.
\]
Note that, differently from \citet{perkovic16}, we also include $\mathbf{X}$ in $\fb{\g}$ to simplify notation. The forbidden set characterizes those covariates that may never be included into a valid adjustment set (see Definition \ref{adjustment}).

\begin{Definition}
	Let $\mathbf{X}$ and $\mathbf{Y}$ be disjoint node sets in a \mpdag{} $\g$.
	We define $\optb{\g}$ as:
	\[
	\optb{\g} =  \pa( \cnb{\g},\g)\setminus \fb{\g}.    
	\]
	\label{def:optset2}
\end{Definition}

%In this broader sense $\pa(X,\g)$ and $\opt{\g}$ are therefore symmetric; mirroring Corollaries \ref{cor:badparents} and \ref{cor:goodparents}.
%	This definition surprisingly does not depend on possible parents of possibly causal nodes as any such node is in fact forbidden.

%\begin{figure}
%	\centering
%	\begin{tikzpicture}[scale=0.15]
%	\tikzstyle{every node}+=[inner sep=0pt]
%	\draw [black] (21.1,-24.1) circle (3);
%	\draw (21.1,-24.1) node {$X$};
%	\draw [black] (34,-24.1) circle (3);
%	\draw (34,-24.1) node {$Z$};
%	\draw [black] (46.1,-24.1) circle (3);
%	\draw (46.1,-24.1) node {$Y$};
%	\draw [black] (27,-13.7) circle (3);
%	\draw (27,-13.7) node {$W$};
%	\draw [black] (37.1,-24.1) -- (43.1,-24.1);
%	\fill [black] (43.1,-24.1) -- (42.3,-23.6) -- (42.3,-24.6);
%	\draw [black] (24.2,-24.1) -- (31,-24.1);
%	\fill [black] (31,-24.1) -- (30.2,-23.6) -- (30.2,-24.6);
%	\draw [black] (25.3,-16.2) -- (22.51,-21.45);
%	\fill [black] (22.51,-21.45) -- (23.33,-20.98) -- (22.44,-20.51);
%	\draw [black] (28.8,-16.2) -- (32.35,-21.59);
%	\fill [black] (32.35,-21.59) -- (32.33,-20.65) -- (31.49,-21.2);
%	\end{tikzpicture}
%	\caption{Causal DAG $\g$ within which $pa(Y,\g) \setminus \f{\g}$ is not a valid adjustment set showing}
%	\label{paY}
%\end{figure}	

\begin{Theorem} \label{thm:optimalSetCPDAG}
	Let $\mathbf{X}$ and $\mathbf{Y}$ be disjoint node sets in a causal \mpdag{} $\g = (\mathbf{V},\mathbf{E})$, such that $\mathbf{Y} \subseteq \possde(\mathbf{X},\g)$. Let the density $f$ of $\mathbf{V}$ be compatible with $\g$ and let $\mathbf{O}=\optb{\g}$. Then the following three statements hold:
	\begin{enumerate}[label = (\roman*)]
		\item\label{res1-cpdag} The set $\mathbf{O}$ is a valid adjustment set relative to $(\mathbf{X,Y})$ in $\g$ if and only if there exists a valid adjustment set relative to $(\mathbf{X,Y})$ in $\g$. 
		\item\label{res2-cpdag} Let $\mathbf{Z}$ be a valid adjustment set relative to $(\mathbf{X},\mathbf{Y})$ in $\g$. If $\mathbf{V}$ follows a causal linear model compatible with $\g$ then
		\[
		a.var(\hat{\tau}^{\mathbf{o}}_{\mathbf{yx}}) \leq a.var(\hat{\tau}^{\mathbf{z}}_{\mathbf{yx}}).
		\]
		\item\label{res3-cpdag} Let $\mathbf{Z}$ be a valid adjustment set relative to $(\mathbf{X},\mathbf{Y})$ in $\g$, such that 
		\[
		a.var(\hat{\tau}^{\mathbf{o}}_{\mathbf{yx}}) = a.var(\hat{\tau}^{\mathbf{z}}_{\mathbf{yx}}).
		\]
		If $\mathbf{V}$ follows a causal linear model compatible with $\g$ and $f$ is faithful to $\g$ then $\mathbf{O}\subseteq \mathbf{Z}$. 
	\end{enumerate}
\end{Theorem}

\begin{Remark}
	In Theorem \ref{thm:optimalSetCPDAG} we assume that $\mathbf{Y} \subseteq \possde(\mathbf{X},\g)$. If 
	$\mathbf{Y} \not\subseteq \possde(\mathbf{X},\g)$ we can instead consider the total effect of $\mathbf{X}$ on $\tilde{\mathbf{Y}} = \mathbf{Y}\cap \possde(\mathbf{X},\g)$, since the total effect of $\mathbf{X}$ on $\mathbf{Y} \setminus \tilde{\mathbf{Y}}$ is $\mathbf{0}$ (see Remark \ref{remark: joint vs non}). Hence, this restriction only limits us from superfluously estimating some zero values. 
	\label{remark:pruned}
\end{Remark}

%As the total effect of $\mathbf{X}$ on one $Y_j \in \mathbf{Y}$ is independent of $\mathbf{Y}_{-j}$ in the linear setting considered here, this restriction only limits us from superfluously estimating some zero values. %\citep{pearl2009causality}. 

Statement \ref{res1-cpdag} implies that our optimal set, similarly to the Adjust$(\mathbf{X},\mathbf{Y},\g)$ set from Definition 12 in \citet{perkovic16}, can be used to check if there exists a valid adjustment set, albeit with the added qualifier that $\mathbf{Y}$ has to be appropriately pruned in advance (see Remark \ref{remark:pruned}). 
%Interestingly, $\optb{\g}$ does not contain any nodes in $\possde(\mathbf{X},\g)$.
Due to statement \ref{res2-cpdag} in Theorem \ref{thm:optimalSetCPDAG} we call $\optb{\g}$ \emph{asymptotically optimal}. Statement \ref{res3-cpdag} implies that in case of faithfulness no other asymptotically optimal set is of smaller or equal size than $\optb{\g}$. 

%In case of faithlessness however, $\optb{\g}$ may contain nodes that do not affect the asymptotic variance, while being also superfluous for unbiasedness. 

%and hence negatively affect the finite sample performance.

%The reason $\optb{\g}$ is asymptotically optimal, is due to it containing as much information on each $Y_j \in \mathbf{Y}$ as any valid adjustment set may, while simultaneously doing the converse for each $X_i \in \mathbf{X}$. 

As a corollary to Theorem \ref{thm:optimalSetCPDAG} jointly with Theorem \ref{prop:order-indep}, the output of Algorithm \ref{Algorithm} is $\optb{\g}$, whenever the starting valid adjustment set $\mathbf{Z}$ is a superset of $\optb{\g}$. It is of course simpler to compute $\mathbf{O}$ directly rather than via pruning.

Remarkably, given a \mpdag{} $\g$ amenable relative to some tuple of node sets $(\mathbf{X},\mathbf{Y})$, such that $\mathbf{Y} \subseteq \possde(\mathbf{X},\g)$, $\optb{\g}$ is not only the optimal set amongst all valid adjustment sets in $\g$ but also among all valid adjustment sets in any DAG $\g[D] \in [\g]$. In fact $\f{\g[D']}=\f{\g}$ and $\opt{\g[D']} = \opt{\g}$ for all DAGs $\g[D'] \in [\g]$ (see Lemmas \ref{lemma:optimalSetEquiv} and \ref{lemma:forbIsSame}). \\

Intuitively, $\optb{\g}$ is constructed to maximize information on $\mathbf{Y}$, while minimizing information on $\mathbf{X}$ and preserving validity. Although one may think that a simpler set, such as $\pa(\mathbf{Y},\g) \setminus (\mathbf{X} \cup \mathbf{Y})$ would suffice for this purpose, this is not generally the case.
%In the special case $\cnb{\g}=\mathbf{Y}$, $\optb{\g}=\pa(\mathbf{Y},\g) \setminus \mathbf{X}$. However, $\pa(\mathbf{Y},\g) \setminus (\mathbf{X} \cup \mathbf{Y})$ is not generally a valid adjustment set and $\optb{\g}$ may contain covariates not in $\pa(\mathbf{Y},\g)\setminus\mathbf{X}$. 
We illustrate this in Example \ref{example-O}. Interestingly, \citet{outcomeIDA} have shown that $\optb{\g}$ can indeed be characterized as $\pa(\mathbf{Y},\tilde{\g}) \setminus (\mathbf{X} \cup \mathbf{Y})$, in a specific latent projection graph $\tilde{\g}$ of $\g$. 

\begin{Example}
	\label{example-O}
	Consider the DAG in Figure \ref{figure:super simple} and the two DAGs from Figure \ref{newgraph}, denoted, respectively, as $\g_1, \g_{2.a}$ and $\g_{2.b}$. We now illustrate how to construct $\optb{\g}$ and the results from Theorem \ref{thm:optimalSetCPDAG}. 
	
	Consider first $\g_1$ and suppose we are interested in the total effect of $V_4$ on $V_6$ as in Example \ref{ex:total effect}. Here $\pa(\mathrm{cn}(V_4,V_6,\g_1),\g_1)=\{V_2,V_3,V_4,V_5\}$ and $\mathrm{forb}(V_4,V_6,\g_1)=\{V_4,V_5,V_6\}$. Therefore, 
	\[\mathbf{O}(V_4,V_6,\g)= \{V_2,V_3\}.\]
	As shown in Example \ref{ex:total effect}, $\{V_2,V_3\}$ is a valid adjustment set relative to $(V_4,V_6)$ in $\g_1$.
	%		 by Theorem \ref{thm:optimalSetCPDAG}, $\{V_2,V_3\}$ is a valid adjustment set relative to $(V_4,V_6)$ in $\g_1$ (see Example \ref{ex:total effect}). 
	It can also easily be verified with Theorem \ref{cor:bignew12} that $\{V_2,V_3\}$ provides a smaller asymptotic variance than any of the alternative valid adjustment sets, as $V_1 \perp_{\g_1} V_6 | \{V_2,V_3,V_4\}$ and $V_3 \perp_{\g_1} V_4 | \{V_2\}$. Hence, it is asymptotically optimal as claimed in Theorem \ref{thm:optimalSetCPDAG}.

	Consider now $\g_{2.a}$. Here $\pa(\cn{\g_{2.a}},\g_{2.a})=\{X,B,C\}$ and $\f{\g_{2.a}}=\{X,Y\}$. Therefore, \[\opt{\g_{2.a}}=\{B,C\}.\] 
	Consider now $\g_{2.b}$. Here $\pa(\cn{\g_{2.b}},\g_{2.b})=\{X,C\}$ and $\f{\g_{2.b}}=\{X,Y\}$. Therefore, \[\opt{\g_{2.b}}=\{C\}.\]
	In these two cases the results from Theorem \ref{thm:optimalSetCPDAG} are corroborated by Example \ref{Apply-3.1} and Example \ref{ex:pruning}, respectively.		
	
	We now discuss why $\optb{\g}$ takes its distinctive form by considering these three examples. By the result from Theorem \ref{cor:bignew12}, an asymptotically optimal valid adjustment set with respect to $(\mathbf{X},\mathbf{Y})$ in $\g$ must contain less or equal information on $\mathbf{X}$ and more or equal information on $\mathbf{Y}$ than any other valid adjustment set.
	%		 while blocking all proper non-causal paths from $\mathbf{X}$ to $\mathbf{Y}$ and not containing nodes in $\fb{\g}$ (see Definition \ref{S-adjustment}).
	
	%		Since  $\de(\mathbf{Y},\g) \subseteq \fb{\g}$ if $\mathbf{Y} \subseteq \de(\mathbf{X},\g)$ (see Remark \ref{remark:pruned}), 
	One might intuitively expect $\pa(\mathbf{Y},\g) \setminus (\mathbf{X} \cup \mathbf{Y})$ satisfies these properties. This is indeed the case for two of the three examples considered here, with \[\opt{\g_{2.a}}=\{B,C\}=\pa(Y,\g_{2.a})\setminus \{X,Y\}\] and 
	\[\opt{\g_{2.b}}=\{C\}=\pa(Y,\g_{2.b})\setminus \{X,Y\}.\]
	This pattern, however, fails to hold for $\g_{1}$. Here, $\pa(V_6,\g_1) \setminus \{V_4,V_6\}=\{V_5\}$ is not a valid adjustment set relative to $(V_4,V_6)$ in $\g_1$, as the mediator $V_5 \in \mathrm{forb}(V_4,V_6,\g_1)$. 
	%		In particular, this means that $V_5$ cannot be used to block the non-causal path $(V_4,V_1,V_2,V_5,V_6)$, even though it is the node closest to $V_6$ on this path. 
	
	The construction of $\boldsymbol{O}(V_4,V_6,\g_1)$ solves this problem by using the next-closest non-forbidden nodes instead, that is, the non-forbidden parents $\{V_2,V_3\}$ of the causal node $V_5$. This ensures validity, while maximizing information on $V_6$ and not providing unnecessary information on $V_4$. 
	Specifically, $V_2$, as the non-forbidden node closest to $V_6$ (and furthest from $V_4$) on the non-causal path $(V_4,V_1,V_2,V_5,V_6)$, is the most efficient choice to block this path. Moreover, $V_3$, although superfluous for validity, contains only information on $V_6$ and therefore improves precision. Interestingly, it does so even though $V_3 \notin \pa(V_6,\g_1)$.
	%		The latter case also illustrates how nodes not in $\pa(\mathbf{Y},\g)$ may be beneficial for precision.
	
	%		n the specific case of $\boldsymbol{O}(V_4,V_6,\g_1)$ these are $\{V_2,V_3\}$.

	%		 $\pa(V_5,\g) \setminus \mathrm{forb}(V_4,V_6,\g_1)$ instead; maximizing information on $V_6$ while ensuring validity. In addition $\boldsymbol{O}(V_4,V_6,\g_1)$ also contains $V_3 \in \pa(V_5,\g) \setminus \mathrm{forb}(V_4,V_6,\g_1)$, which similarly to $C \in \opt{\g_{2.a}}$ and $C \in \opt{\g_{2.b}}$, while superfluous for validity 
	%		In fact $\g_1$ serves as good ill
	
	%		 $\boldsymbol{O}(V_4,V_6,\g_1) \neq \pa(V_6,\g_1) \setminus \{V_4\}$.
\end{Example}

\section{Simulation study}

\label{sim-MSE}

\begin{figure}[p]		
	
	\centering
	\captionsetup[subfloat]{position=top,labelformat=empty}
	\subfloat[\textbf{Single interventions}]{
		\centering
		\includegraphics[height=6.8cm, width=7.05cm]{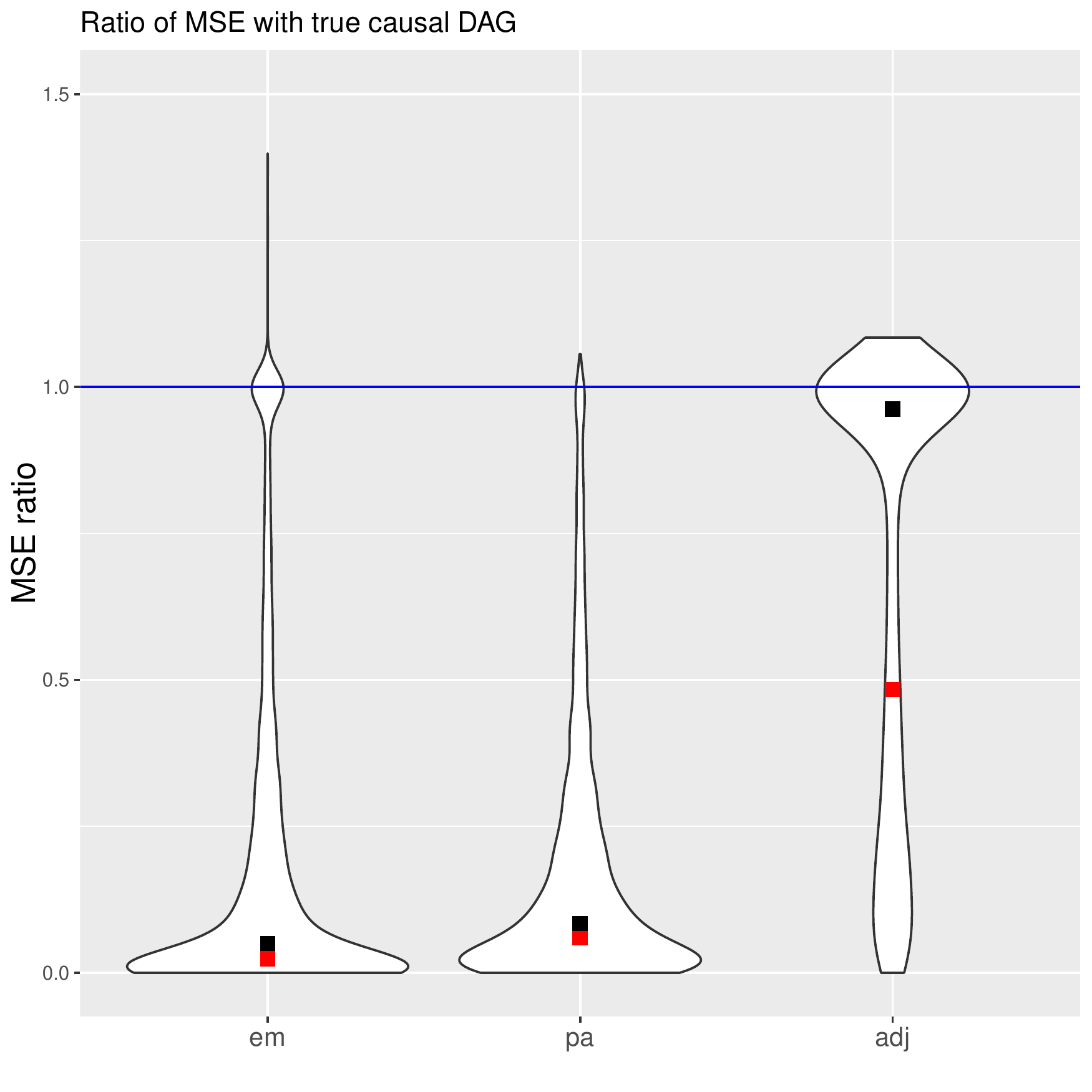}
		\includegraphics[height=6.8cm, width=7.05cm]{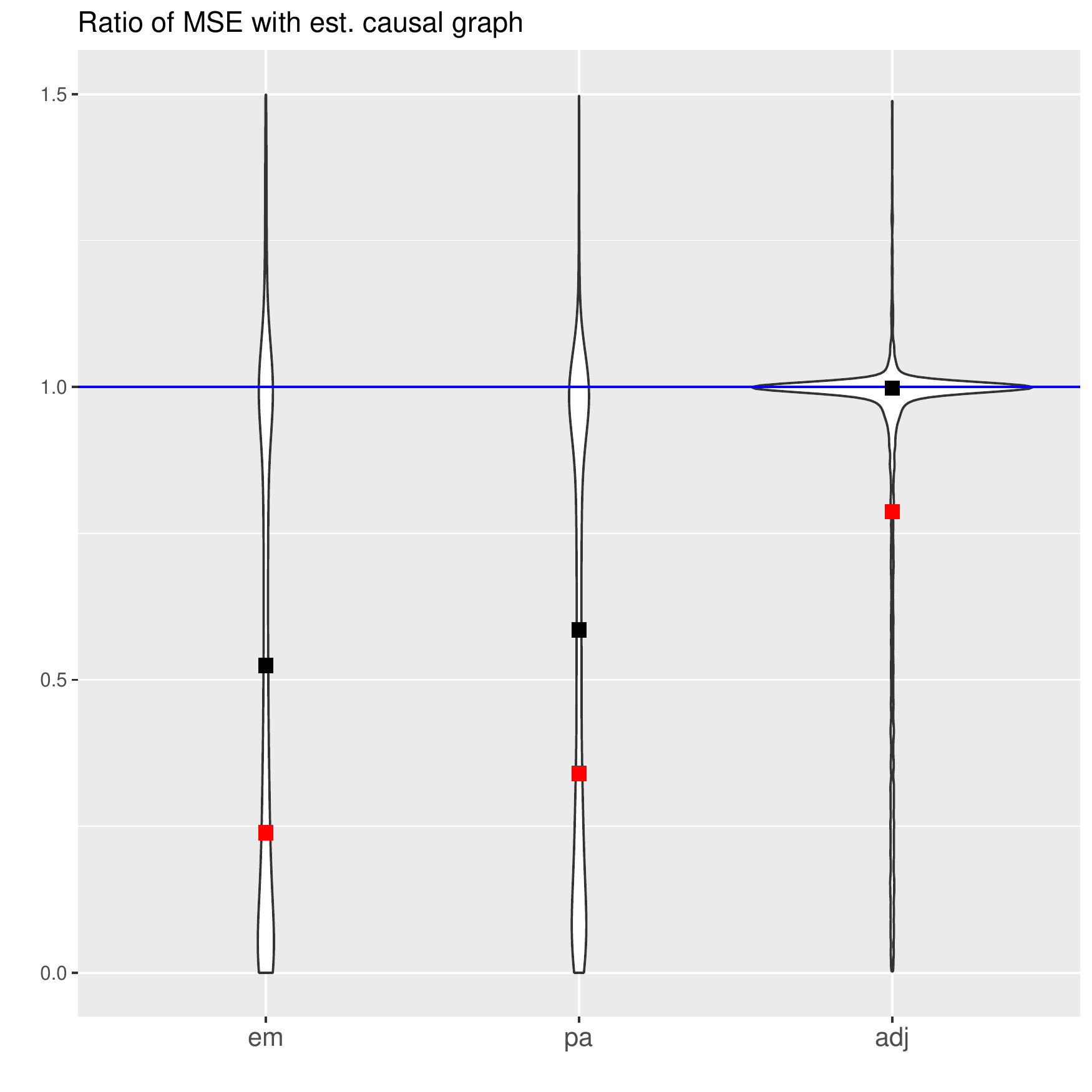}
	} \par \bigskip 
	\subfloat[\textbf{Joint interventions}]{
		\centering
		\includegraphics[height=6.8cm, width=7.05cm]{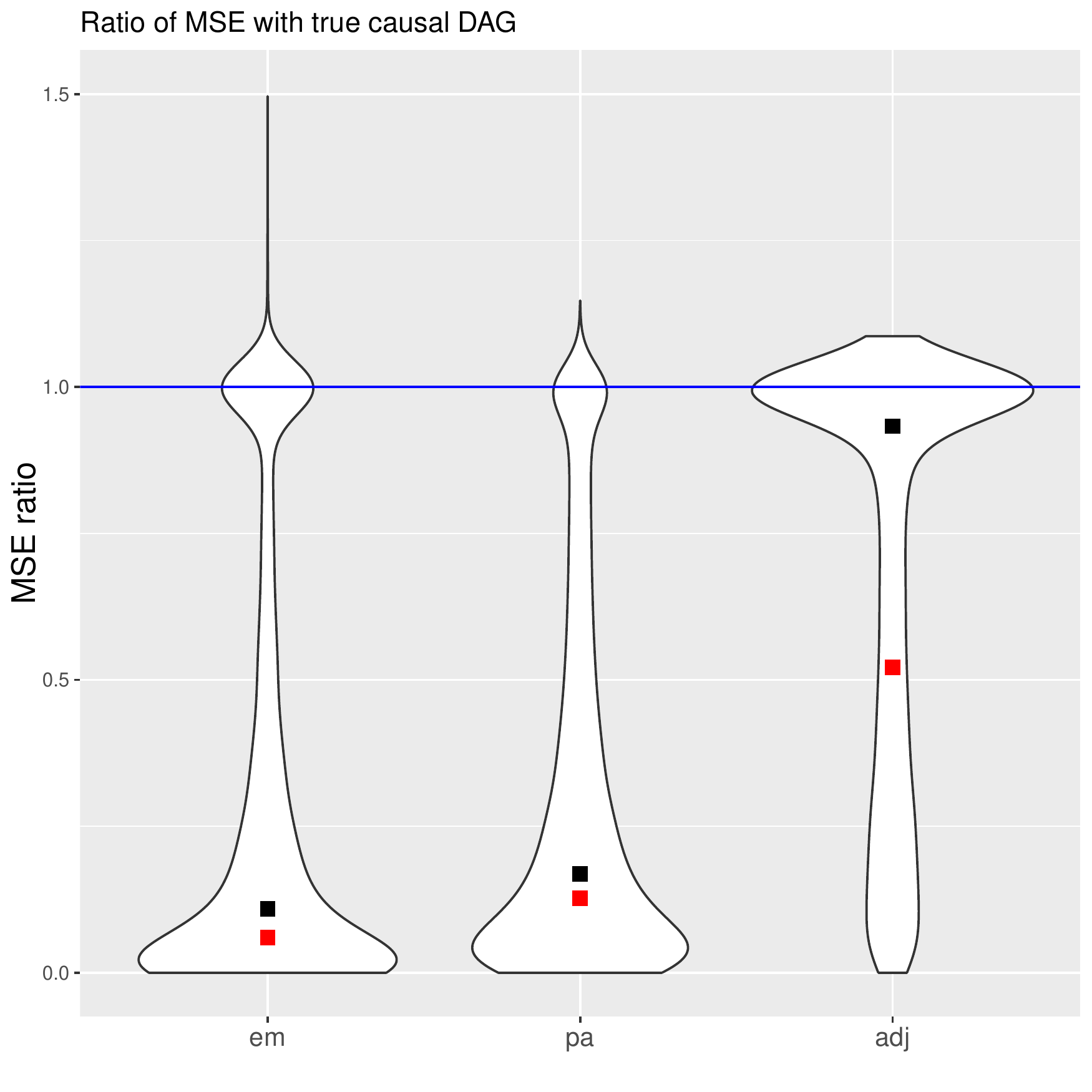}
		\includegraphics[height=6.8cm, width=7.05cm]{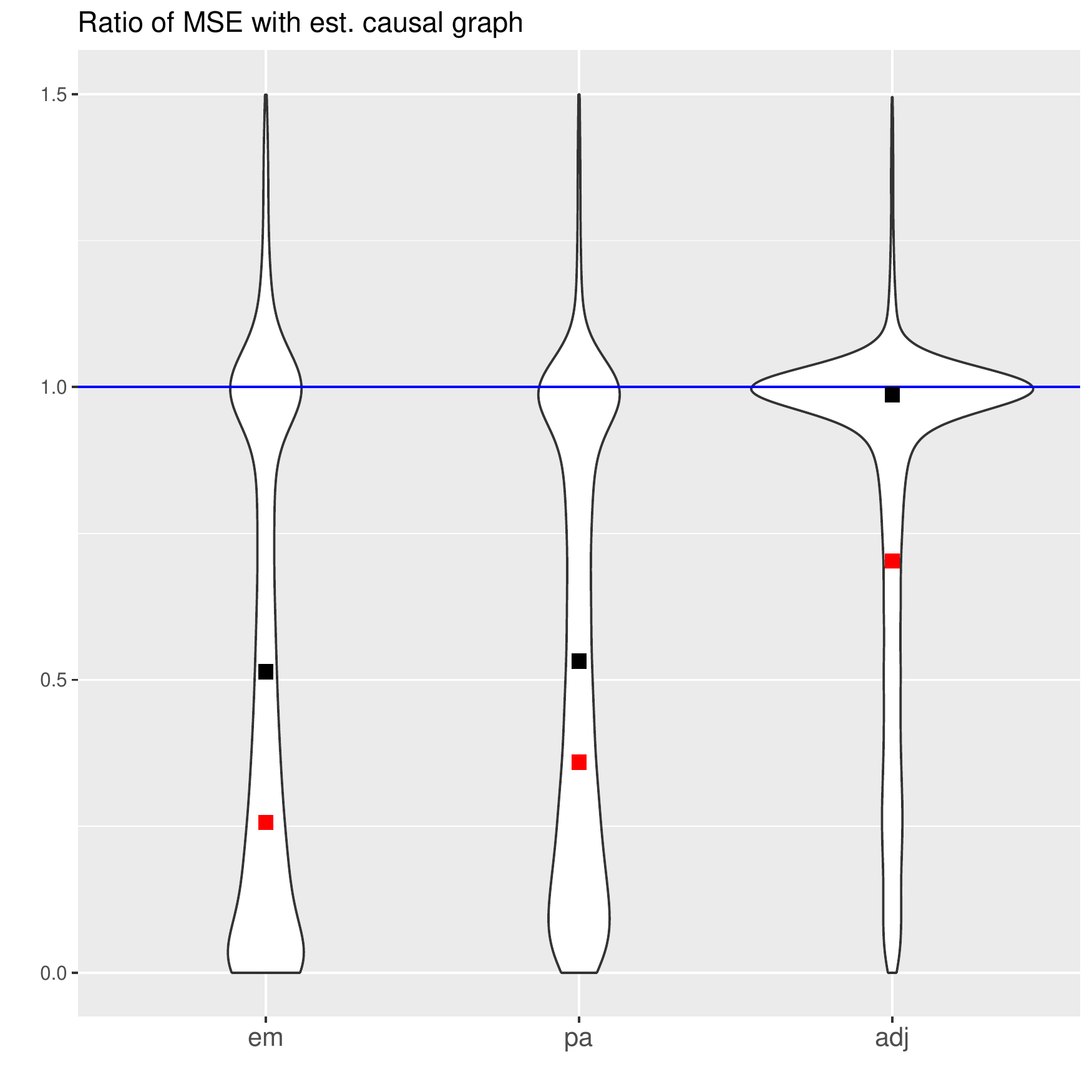}
	}
	\caption{ 
		Violin plots of the ratios of the mean squared errors provided by $\text{\textbf{O}}(\mathbf{X},Y,\g)$ and the three alternative adjustment sets \ref{AS-em} - \ref{AS-adj} from Section \ref{sim-MSE}, respectively. 
		The true causal DAG cases are on the left and the graph estimate cases on the right. The single intervention cases are at the top and the joint intervention cases are at the bottom.  The red squares show the geometric average of the ratios, the black squares the median. While the summary measures are computed with all ratios, ratios larger than $1.5$ are not shown in the plots.  From left to right the percentage of ratios larger than $1.5$ in the single intervention setting is $0.2, 0.0, 0.0, 7.1, 1.1,$ and  $1.5$ percent, respectively. For the joint interventions setting the percentages are $0.2, 0.0, 0.0, 8.8, 2.2$ and $1.8$ percent, respectively.} 
	%		We removed the smallest and largest 2.5 percent of ratios to draw the violin plot but not to compute the summary measures.
	\label{violin}
\end{figure}

We investigate the finite sample performance of adjusting for $\optb{\g}$ by sampling data from randomly generated causal linear models and comparing the empirical mean squared error provided by $\optb{\g}$ to three alternative adjustment sets. 
A detailed explanation of our simulation setup is given in Section \ref{setupMSE} of the Supplement \citep{supplement}.

We randomly generate a total of 10'000 DAGs, with the number of nodes chosen from $\{10,20,50,100\}$ and the expected neighborhood size from $\{2,3,4,5\}$. Each graph is associated with a causal linear model. The edge coefficients of the model are drawn independently from a uniform distribution on $[-2,-0.1]\cup [0.1, 2]$, and the errors are either drawn from a Gaussian distribution, a t-distribution with 5 degrees of freedom, a logistic distribution or a uniform distribution, with variances in the range of $[0.5, 1.5]$.

For each DAG $\g[D]$, we randomly draw $(\mathbf{X},Y)$ such that $|\mathbf{X}| \in \{1,2,3\}$ and $Y \in \cap_{X_i\in \mathbf{X}} \de(X_i,\g[D])$. We do this for the following two reasons. First, the restriction to a singleton $Y$ is sensible by Remark \ref{remark: joint vs non}. Secondly, if $Y \notin \de(X_i,\g[D])$ for some $X_i\in \mathbf{X}$ then the corresponding entry of the total effect $(\tau_{y\mathbf{x}})_i =0$ (see Remark \ref{remark: joint vs non}). 
%	Deciding whether $Y \notin \de(X_i,\g[D])$ is, however, a causal discovery problem and does not depend on the choice of adjustment set. 
We then verify whether there exists a valid adjustment set with respect to $(\mathbf{X},Y)$ in both the DAG $\g[D]$ and its CPDAG $\g[C]$. If not, we resample $\mathbf{X}$ and $Y$. 

For each causal linear model we generate 100 data sets with sample sizes $n \in \{125,500,2000,10000\}$. We then consider two settings: i) We suppose knowledge of the true causal DAG $\g[D]$ and ii) we estimate a graph $\widehat{\g}$ from the data. If the errors are drawn from a Gaussian distribution $\widehat{\g}$  is estimated with the Greedy Equivalence Search (GES) algorithm \citep{Chickering02}, otherwise with the Linear Non-Gaussian Acyclic Models (LiNGAM) algorithm \citep{shimizu2006linear}. In both cases we use the algorithms as implemented in the {\em pcalg} R-package \citep{Kalisch2012}.

We then compute total effect estimates, by adjusting for $\mathbf{O}(\mathbf{X},Y,\g)$ and three alternative adjustment sets. This is done with respect to both the true causal DAG $\g[D]$ and the estimated causal graph $\g[\widehat G]$, with two special cases for the estimates with respect to $\g[\widehat G]$. Firstly, no estimate is returned if there was no valid adjustment set relative to  $(\mathbf{X},Y)$ in $\widehat{\g}$, i.e., these cases are discarded for the mean squared error computation. In such cases, we recommend the use of alternative total effect estimators such as the {\em IDA} algorithm by \citet{maathuis2009estimating} and the {\em jointIDA} algorithm by \citet{nandy2017estimating}. 
Secondly, $\boldsymbol{0}$ is returned as the estimate whenever $Y \notin \possde(\mathbf{X},\widehat{\g})$, since the total effect on a non-descendant is $\boldsymbol{0}$. The pair $(\mathbf{X},Y)$ is sampled, ensuring that these two special cases do not occur in either the true DAG or its corresponding CPDAG. 

The three alternative adjustments sets are:
\begin{enumerate}[label=(\roman*)]
	\item The empty set, representing a non-causal baseline. It is generally not a valid adjustment set and is denoted by ``em".  \label{AS-em}
	\item The set $\pa(\mathbf{X},\g) \setminus \text{forb}(\mathbf{X},Y,\g)$, which in the setting $|\mathbf{X}|=1$ is the valid adjustment set $\pa(X,\g)$. If $|\mathbf{X}| > 1$, it is not generally a valid adjustment set. It is denoted by ``pa". \label{AS-pa}
	\item The valid adjustment set $\text{Adjust}(\mathbf{X},Y,\g)$ from \citet{perkovic16}. It is denoted by ``adj". \label{AS-adj}
	%		\item The valid adjustment set of all non-forbidden nodes, denoted ``nforb".
	%		\item The set of all nodes except for $\mathbf{X}$ and $Y$, i.e., adjusting for all available covariates.  It is generally not a valid adjustment set and is denoted by ``all". \label{AS-all}
\end{enumerate}

For each causal linear model, we thus have four adjustment sets in two graphical  settings. In each of these cases, we compute the empirical mean squared error of our respective estimates with respect to the true total effect. We emphasize that we do \emph{not} consider the estimated standard errors or residuals from the regression analyses. To quantify the advantage of $\text{\textbf{O}}(\mathbf{X},Y,\g)$, we compute the ratio of the mean squared error corresponding to $\text{\textbf{O}}(\mathbf{X},Y,\g)$ and each of the three alternative adjustment sets. This is done separately for the two graphical settings.

Figure \ref{violin} is a violin plot of these ratios.
We see that $\text{\textbf{O}}(\mathbf{X},Y,\g)$ provides consistently smaller mean squared errors than any of the considered alternatives. Except for $\text{Adjust}(\mathbf{X},Y,\g)$, all alternative sets are clearly outperformed by $\text{\textbf{O}}(\mathbf{X},Y,\g)$, with geometric averages below $0.5$. As might be expected, the gain becomes smaller when the underlying causal DAG has to be estimated, but it remains respectable. Notably, the proportion of ratios larger than $1.5$ is small, even negligible when the true DAG is known. For a more thorough discussion how the ratios behave depending on the parameters see Section \ref{resultsMSE} of the Supplement.
The bulges at $1$ are due to two reasons. Firstly, cases in which the compared sets are similar or the same. Secondly, cases in which $Y \notin \possde(\mathbf{X},\widehat{\g})$ occurs for a considerable number of the estimated graphs $\widehat{\g}$ (see Section \ref{graphIssues}  of the Supplement).

The only true contender to $\text{\textbf{O}}(\mathbf{X},Y,\g)$ in terms of performance is $\text{Adjust}(\mathbf{X},Y,\g)$. It should be noted, however, that $\text{Adjust}(\mathbf{X},Y,\g)$ is a superset of $\text{\textbf{O}}(\mathbf{X},Y,\g)$ and hence will be more cumbersome to measure in practice (see also Figure \ref{sizediff} of the Supplement). 

%	As it still tends to provide a larger mean squared error than $\text{\textbf{O}}(\mathbf{X},Y,\g)$, there appears to be no practical advantage to using it.

Another point worth noting is the bad performance of $\textrm{pa}(\mathbf{X},\g) \setminus \textrm{forb}(\mathbf{X},Y,\g)$. Even though this set is a valid adjustment set if $|\mathbf{X}|=1$, it only provides a small gain compared to the empty set, especially when the graphical structure has to be estimated. This aptly illustrates the importance of taking efficiency considerations into account when choosing a valid adjustment set.

%	One drawback of the optimal set is its typically increased size. As can be seen in Figure \ref{S-sizediff} from the Supplement, the optimal set is in most cases larger than the parent set. 

%	
In summary, these results indicate that there are benefits to using $\text{\textbf{O}}(\mathbf{X},Y,\g)$. These benefits decrease when the underlying causal structure is not known in advance, but do remain respectable. 

\section{Real data example}
\label{sec-realData}

\begin{figure}[t!]
	\centering
	\includegraphics[scale=.5]{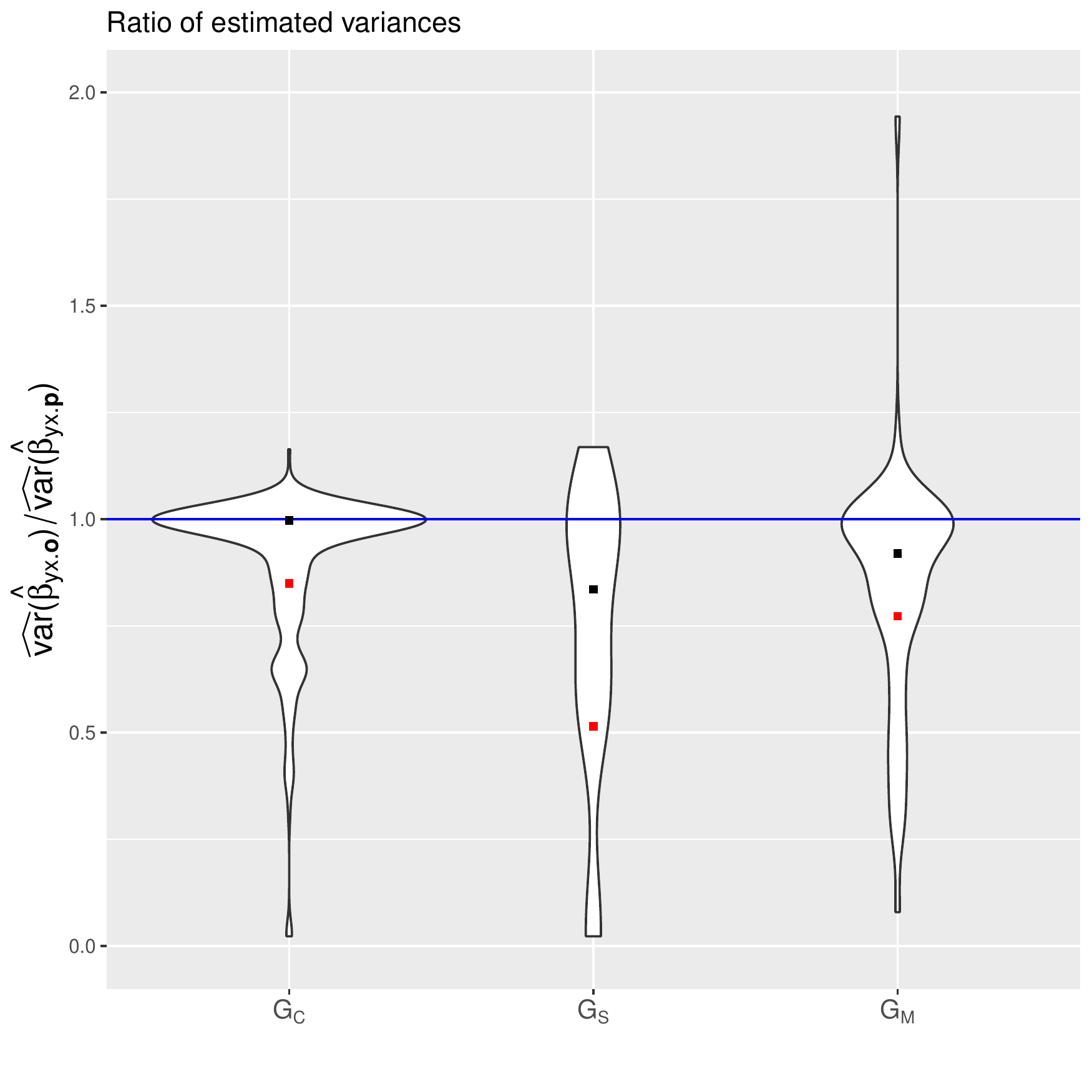}
	\caption{Violin plots of the ratios $\widehat{\mathrm{var}}(\hat{\beta}_{yx.\mathbf{o}})/\widehat{\mathrm{var}}(\hat{\beta}_{yx.\mathbf{p}})$, with $\mathbf{O}=\mathbf{O}(X,Y,\g)$ and $\mathbf{P}=\pa(X,\g)$, for all pairs of nodes $(X,Y)$ such that $Y \in \de(X,\g)$ and $\opt{\g}\neq\pa(X,\g)$ in 8 experimental conditions from \citet{sachs2005causal}, obtained under the assumption that $\g_C,\g_S$ or $\g_M$, respectively, is the true underlying causal graph. Here, $\g_C$ denotes the consensus graph, $\g_S$ the acyclic graph estimated by \citet{sachs2005causal} and $\g_M$ the acyclic graph estimated by \citet{mooij2013cyclic}. The red squares show the geometric average of the ratios, the black squares the median.}	
	\label{CIlength}
\end{figure} 

Our result can easily be integrated into existing approaches to covariate adjustment. Using $\opt{\g}$ to estimate the total effect $\tau_{yx}$  instead of, for example, $\pa(X,\g)$ only requires the minimal additional effort of computing $\opt{\g}$ from the causal graph $\g$. And yet, replacing $\pa(X,\g)$ with $\opt{\g}$ can only improve the asymptotic variance when the true causal graph $\g$ is used. 
%an efficiency gain at essentially no cost. 
Of course, errors in the used graph $\g$ and finite sample considerations might lead to cases where the use of $\opt{\g}$ actually leads to a loss of efficiency in practice, but our simulations (Section \ref{sim-MSE}) indicate that this risk is manageable and that an overall efficiency gain, at essentially no cost, is the norm. 

%the ease with which they can be integrated into existing approaches to covariate adjustment while . 
%They maintain the general framework \cite{ZanderL19} 

To investigate this further, we apply our results to the single cell data of \citet{sachs2005causal}. This data set consists of flow cytometry measurements of 11 phosphorylated proteins and phospholipids in human T-cells, collected under 14 different experimental conditions. Each experimental condition corresponds to a different intervention on the abundance or activity of the proteins. We chose this data set due to the the availability of a consensus graph \citep[see Figure Figure 5.a) in][]{mooij2013cyclic} and the large sample size.
%	 and the availability of interventional data.

Given that there is some uncertainty regarding the consensus graph, we apply our results using the following three different graphs: the consensus graph, the DAG estimated by \citet{sachs2005causal} and the DAG estimated under the restriction to at most 17 edges by \citet{mooij2013cyclic}. These three DAGs are given in Figure 5 of \citet{mooij2013cyclic} and we 
%	refer to them as "Consensus", "Sachs et al." and "Mooij and Heskes", respectively.
denote them by $\g_C,\g_S$ and $\g_M$, respectively.

%	We had three gold standard graphs available: the consensus graph
%	We used the graph given by \citet{mooij2013cyclic}[Figure 5.\ a)] as our gold standard graph 
Our data analysis is as follows. We first log transform the data as it is heavily right skewed. For each of the three DAGs we then do the following: Restricting ourselves to the 8 experimental conditions for which \citet{mooij2013cyclic} provide a graphical interpretation of the condition's effect, we adjust our starting DAG accordingly. For each such adjusted DAG $\g$, we then compute all pairs $(X,Y)$ of nodes in $\g$, such that $Y \in \de(X,\g)$, to ensure that there is a non-trivial total effect to estimate, and $\opt{\g} \neq \pa(X,\g)$, to ensure that we compare different estimators. For each such pair, we compute the least squares regressions of $Y$ on $X$ and $\mathbf{O}$ as well as of $Y$ on $X$ and $\mathbf{P}$. We note that $\hat{\beta}_{yx.\mathbf{o}}$ and $\hat{\beta}_{yx.\mathbf{p}}$ are estimators for the total effect of $X$ on $Y$ in the considered data regime; not necessarily in the observational regime.

As the true total effects are unknown, we compare the least squares coefficient variance estimates $\widehat{var}(\hat{\beta}_{yx.\mathbf{o}})$ and $\widehat{var}(\hat{\beta}_{yx.\mathbf{p}})$ by considering their ratio 
$\widehat{var}(\hat{\beta}_{yx.\mathbf{o}})/\widehat{var}(\hat{\beta}_{yx.\mathbf{p}}).$ 
%	of these two least squares regressions
%by considering their ratio. 
Note that, differently from Section \ref{sim-MSE}, where we are able to compute the empirical mean squared error with respect to the known true total effect, this approach raises some concerns regarding post-selection inference for the two estimated graphs $\g_{S}$ and $\g_{M}$, which we disregard here. 
%	We ignore these concerns
%	Assuming an underlying causal linear model with finite fourth moments this is a consistent estimator of the corresponding asymptotic variance ratio by Proposition \ref{avarDAG} and the continuous mapping theorem. 

%This explain the difference between the bulge of ratios close to $1$.
\begin{figure}
	\centering
	\begin{tikzpicture}[>=stealth',->,>=latex,shorten >=1pt,auto,node distance=1.8cm,scale=1.1,transform shape]
	\tikzstyle{vertex}=[circle, draw, inner sep=2pt, minimum size=0.15cm]
	\node[vertex]         (V1)                        {$X$};
	\node[vertex]         (V2) [right of= V1]  		{$M$};
	\node[vertex]         (V3) [right of = V2] 		{$Y$};
	\node[vertex]         (U1)  [above of =V1]	{$P$};
	\node[vertex]       	(U2) [above of =V3]	{$O$};
	
	\draw[->] 	(V1) 	edge  node[below]  {$1$} (V2);
	\draw[->] 	(V2) 	edge  node[below]  {$1$}  (V3);
	\draw[->] 	(U1) 	edge   node[left]  {$1$} (V1);
	\draw[->] 	(U1) 	edge   node[right]  {$1$} (V2);
	\draw[->] 	(U2) 	edge   node[left]  {$1$} (V2);
	\draw[->] 	(U2) 	edge   node[right]  {$-1$} (V3);
	
	\end{tikzpicture}
	\caption{Causal DAG  from Section \ref{sec-realData}.}
	%			, with the edge weight denoting edge coefficients for the corresponding causal linear model.}
	%	and \ref{emma:cpdag}.}
	
	\label{figure:faithless}
\end{figure}
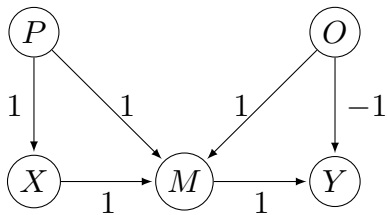 

Figure $\ref{CIlength}$ shows violin plots of these ratios, aggregated over all considered $(X,Y)$ pairs in the 8 experimental settings, with one plot for each of the graphs $\g_C,\g_S$ and $\g_M$.  
%	Since all three considered gold standard graphs are small (p=11) and not particularly dense, $\mathbf{O}(X,Y,\g) = \pa(X,\g)$ for many pairs $(X,Y)$. In this case, the ratio $\widehat{var}(\hat{\beta}_{yx.\mathbf{o}})/\widehat{var}(\hat{\beta}_{yx.\mathbf{p}})$ will always be $1$. For ease of presentation we removed these cases before plotting Figure $\ref{CIlength}$. 
The plots in Figure $\ref{CIlength}$ show that using $\opt{\g}$ instead of $\pa(X,\g)$ results in geometric means smaller than $1$ for all three graphs. In particular, only few of the ratios are larger than $1$, with only one larger than $1.2$, showing that this gain is obtained at little risk of a potential downside. The gains are rather modest, however. This is likely due to the small size and sparsity of the considered graphs. In such settings, even when $\opt{\g}\neq\pa(X,\g)$, the two sets will often share nodes and only differ in minor ways. As a result they provide similar asymptotic variances. In fact, this behaviour can also be seen in our simulations where small graph sizes and small expected neighborhood sizes result in mean squared error ratios closer to $1$ (see Figure \ref{moreBox} of the Supplement).

It is also interesting to consider how the gain in efficiency differs between the three graphs. It is smallest for $\g_C$ and largest for $\g_S$. 
%	This difference is likely due to the large proportion of ratios close to $1$, a behaviour most pronounced for $\g_C$, smaller but still notable for $\g_M$ and almost absent for $\g_S$. Considering that we removed the cases where $\opt{\g}=\pa(X,\g)$, this behaviour might be surprising, but there is a likely explanation. 
As discussed by \citet{mooij2013cyclic}, several strong faithfulness violations seem to be present in the considered data set. While Theorem \ref{thm:optimalSetCPDAG} does not require faithfulness to hold, faithfulness violations can lead to cases where $\opt{\g}$ and $\pa(X,\g)$ differ and yet provide the same asymptotic variance, making them both asymptotically optimal. For an example consider the DAG $\g$ in Figure \ref{figure:faithless} and the corresponding causal linear model with the edge weights indicated on the edges (and arbitrary error variances). Here, $\opt{\g}=\{P,O\}$ and $\pa(X,\g)=\{P\}$ but due to the non-faithfulness of the causal linear model, $O \ci Y | \{X,P\}$, even though $O \in \pa(Y,\g)$. Since $O \perp_{\g} X | P$, we can thus conclude with Lemma \ref{close} that for the considered causal linear model.
\[a.var(\hat{\tau}_{\mathbf{yx}}^{\mathbf{o}}) =
a.var(\hat{\tau}_{\mathbf{yx}}^{p}).\]
%	As $\pa(X,\g) \subset \opt{\g}$, we can even expect $\pa(X,\g)$ to perform slightly better than $\opt{\g}$ in finite samples.
%	We stress that Theorem \ref{thm:optimalSetCPDAG} nonetheless holds and $\opt{\g}$ is still asymptotically optimal, but so is $\{P\}$. 
It appears that this issue is most prevalent for $\g_C$, as it is the densest of the three considered graphs and faithfulness violations require multiple connecting paths between nodes. In $\g_S$,  the least dense of the three considered graphs, this issue appears to be less prominent, leading to a larger gain. 

%While the faithfulness violations did result in many variance ratios close to $1$ and thus lead to an increase in the observed geometric means, we emphasize that there were still very few ratios larger than $1$.  The faithfulness violations' effect was to make $\opt{\g}$ perform more similarly to $\pa(X,\g)$, not worse. This illustrates the robustness of our results.

%	Nonetheless, the variance estimate provided by $\mathbf{O}$ is often smaller than the one provided by  $\mathbf{P}$; in some instances substantially so. 

%	The converse never occurs by a respectable margin, with no ratio being larger than 1.05. 
% Considering the difficulties inherent to this data set (faithlessness, cyclicity) this is a promising results.

%	$\pa(X,\g)$ only provides the smaller standard error in one single case. 
%	
%	Nonetheless, $\mathbf{O}(X,Y,\g)$ does regularly provides smaller standard errors and hence also smaller confidence intervals than $\pa(X,\g)$ more consistently and by larger margins than the converse occurs. We  thus obtain encouraging  results.
%	 even though the data almost certainly violates some of our assumptions (acyclicity,  linearity,  hidden confounders). 
	\section{Discussion}

\label{sec:dis}

We provide a series of results on graphical criteria for efficient total effect estimation via adjustment in causal linear models. Specifically, we present a graphical criterion to qualitatively compare the asymptotic variance that many pairs of valid adjustment sets provide. Further, supposing the existence of a valid adjustment set, we provide a variance reducing pruning procedure as well as an asymptotically optimal valid adjustment set. These results formalize and strengthen existing intuition regarding efficiency. They form a versatile tool set for choosing among valid adjustment sets, a choice that can have a significant impact on the mean squared error.

We do, however, only consider total effect estimation via covariate adjustment. Other estimators, such as ensemble estimators or the front-door criterion \citep{hayashi2014estimating} may be more efficient.

Our results require an in-depth understanding of the causal structure in the form of a causal DAG or an amenable \mpdag{}. However, our results are not considerably more affected by this difficulty than covariate adjustment as a whole. For example, suppose that we consider singleton $X$ and $Y$, such that $Y \subseteq \possde(X,\g)$ (see Remark \ref{remark:pruned}). Then knowledge of an amenable \mpdag{} $\g$ is required both for $\pa(X,\g)$ to be a valid adjustment set and for $\opt{\g}$ to be identifiable. In practice, $\opt{\g}$ may be more sensitive to graph estimation errors than $\pa(X,\g)$, since $\pa(X,\g)$ only relies on estimating the local neighborhood of $X$ accurately. Nonetheless, our simulations indicate that even when the underlying causal DAG has to be estimated, $\optb{\hat{\g}}$ typically provides a smaller mean squared error than $\pa(X,\g)$ or $\emptyset$. 

Since our results cover DAGs, CPDAGs and maximal PDAGs, it is natural to ask: can they be extended to settings with latent variables? The answer is: partially. Theorem \ref{cor:bignew12} extends to settings with latent variables and without selection bias, by simply changing d-separation to (definite) m-separation \citep{RichardsonSpirtes02,zhang2008causal} in the latent variable graph (MAG or PAG) and then using Theorem 4.18 from \citet{RichardsonSpirtes02} (see also Lemma 20 in \citet{zhang2008causal}) and Lemma 26 from \citet{zhang2008causal}. However, Theorem \ref{thm:optimalSetCPDAG} does not extend to latent variable models as can be seen in Example \ref{hidden}. If we suppose here that $C$ is latent and only consider the valid adjustment sets that do not contain $C$, then the valid adjustment set providing the optimal asymptotic variance depends on the edge coefficients and error variances.

There is a caveat to this partial extension. Since all our results are with respect to valid adjustment sets, they do not apply if there is unmeasured confounding, i.e., no valid adjustment set is fully observed. The one partial exception to this is Theorem \ref{cor:bignew12}, which, under the assumption of Gaussian errors for the causal linear model, does hold for arbitrary adjustment sets. Interestingly, there is research indicating that in the presence of unmeasured confounding, the broad guideline of choosing the adjustment set $\mathbf{Z}$ in way that minimizes information on $\mathbf{X}$ while maximizing information on $\mathbf{Y}$ should still be followed to minimize bias amplification \citep{pearl10,wooldridge2016should,ding2017instrumental}.  As result, a pruning procedure along the lines of Algorithm \ref{algo: pruning}  might still be warranted in this setting.

%Add Pearl z-bias citation

%	Since our results cover DAGs, CPDAGs and \mpdag{}s, it is natural to ask whether they could be extended to settings with latent variables? The answer is: not entirely. Theorem \ref{cor:bignew12} extends to settings with unmeasured variables and without selection bias, by simply changing d-separation to (definite) m-separation \citep{RichardsonSpirtes02,zhang2008causal} in the latent variable graph and then using Theorem 4.18 from \citet{RichardsonSpirtes02} (see also Lemma 20 in \citealp{zhang2008causal}) and Lemma 26 from \citet{zhang2008causal}. 
%	However, Theorem \ref{cor:bignew12} does not extend to the latent variable models as can be seen in Example \ref{hidden}.
%	In the PAG corresponding to Figure \ref{newgraph}(b) with $S_3$ marginalized out  $\{S_1\},\{S_1,S_2\}$ and the empty set are still valid adjustment sets relative to $(X,Y)$ and obviously provide the same asymptotic variances. Hence, there is no graphically identifiable asymptotically optimal valid adjustment set in this PAG.

Our results do not apply to non-amenable CPDAGs and non-amenable \mpdag{}s. In this setting one can use the IDA algorithm from \citet{maathuis2009estimating} and \cite{maathuis2010predicting} for CPADGs or the modified version by \citet{perkovic17} for \mpdag{}s. Both output a list of possible total effect estimates by adjusting for the possible parent sets of $\mathbf{X}$, one for each DAG compatible with the considered causal graph. As the parents are often an inefficient valid adjustment set, one may wonder whether it is possible to apply our results to improve the IDA algorithm's efficiency. This is indeed possible, as shown by \citet{outcomeIDA}.

Finally, another possible generalization is to consider settings with selection bias. \citet{correa2018generalized} give a necessary and sufficient graphical criterion for causal effect estimation under confounding \textit{and} selection bias. It remains to be investigated whether the results presented in this paper generalize to this setting.

%Other avenues for future research are graphical criteria for efficient total effect estimation with approaches other than covariate adjustment, such as settings with factorial target variables and instrumental variable estimators. 

%Here are some natbib examples. You can cite examples using the citation key \citep{TM83} in your .bib file. (On Overleaf, you can access the .bib file via the Project menu.) There are commands for in-text citations, like \citet{GMP81}. And you can pass an option to specify additional details, such as a page or chapter number, as an option \citep[p. 130]{Ful83}.
%
%We hope you find Overleaf useful, and please let us know if you have any feedback using the help menu above. 

\bibliography{mybib}%
\bibliographystyle{apa}%

\newpage
\renewcommand\appendixpagename{Supplement}
\renewcommand\appendixtocname{Supplement}
\appendix
\appendixpage

This is the Supplement to ``Graphical criteria for efficient total effect estimation via adjustment in causal linear models'' which we will refer to as the Main paper. Results from the Main paper are referenced to by their original numbering (e.g., Proposition \ref{avarDAG}) whereas references to results in the Supplement begin with a letter (e.g., Lemma \ref{avar.equ}).

 \section{Graphical preliminaries and existing results}
 
%  This section introduces some additional notation used throughout the proofs presented in this supplement. Important existing results are given for added comprehensiveness. Notations and results used in only one section are presented in that section. 

\subsection{Graphical preliminaries and examples}
\label{graph-supp}

\noindent{\bf Graphs}. 
We consider simple graphs with a finite node set $\mathbf{V}$ and an edge set $\mathbf{E}$, where
edges can be either directed ($\to$) or undirected ($-$).  
If all edges in $\mathbf{E}$ are directed ($\to$), then $\g = (\mathbf{V},\mathbf{E})$ is a \emph{directed graph}. If all edges in $\mathbf{E}$ are directed or undirected, then $\g =  (\mathbf{V},\mathbf{E})$ is a \textit{partially directed graph}. A graph $\g' = (\mathbf{V'},\mathbf{E'})$ is an induced subgraph of $\g = (\mathbf{V,E})$ if $\mathbf{V'} \subseteq \mathbf{V}$ and $\mathbf{E'}$ contains all edges in $\mathbf{E}$ between nodes in $\mathbf{V'}$. 
% Throughout this paper, any graph is a partially directed graph. 

\vsp \noindent{\bf Paths}. 
Two nodes are \emph{adjacent} if there exists an edge between them.
A \textit{path} $p$ from a node $X$ to a node $Y$ in a graph $\g$ is a sequence of distinct nodes $(X = Z_1, \dots, Z_m=Y)$ such that 
$Z_i$ and $Z_{i+1}$ are adjacent in $\g$ for all $i \in \{1, \dots, m-1\}$. Then $X$ and $Y$ are called \emph{endpoints} of $p$. We use $p(Z_i, Z_j)$ to denote the \textit{subpath} $(Z_i, Z_{i+1}, \dots, Z_j)$ of $p$, with possibly $Z_i =Z_j$ in which case the subpath is simply a node. A path from a set $\mathbf{X}$ to a set $\mathbf{Y}$ is a path from some $X \in \mathbf{X}$ to some $Y \in \mathbf{Y}$. A path from a set $\mathbf{X}$ to a set $\mathbf{Y}$ is \textit{proper} if only the first node is in $\mathbf{X}$ (cf.\ \citet{shpitser2010validity}). 
A path $p = (Z_1, \dots Z_m)$ is called \emph{directed} from $Z_1$ to  $Z_m$ if $Z_i \rightarrow Z_{i+1}$ for all $ i \in \{1,\dots,m-1\}$. It is called \emph{possibly directed} from $Z_1$ to $Z_m$ if there are no $i, j \in \{1,\dots,m\}$, $i <j$, such that $Z_{i} \leftarrow Z_{j}$ (cf.\ \citet{perkovic17}).

\begin{Remark}
	Our definition of a possibly directed path is non-standard, as $j-i>1$ is allowed. This is required for its use in \mpdag{}s, as we show in Example \ref{ex: CPDAGs, PDAGs and possibly directed}. 
\end{Remark}

\noindent{\bf Ancestry}.  If $X \to Z$, then $X$ is a \textit{parent} of $Z$ and $Z$ is a \textit{child} of $X$. 
If there is a  directed path from $X$ to $Y$, then $X$ is an \textit{ancestor} of $Y$ and $Y$ a  \textit{descendant} of $X$. If there is a  possibly directed path from $X$ to $Y$, then $X$ is a \textit{possible ancestor} of $Y$ and $Y$ a  \textit{possible descendant} of $X$. 
We use the convention that every node is an ancestor, possible ancestor, descendant and possible descendant of itself. 
The sets of parents, ancestors, descendants, possible ancestors  and possible descendants of $X$ in $\g$ are denoted by $\pa(X,\g)$, $\an(X,\g)$, $\de(X,\g)$, $\possan(X,\g)$ and $\possde(X,\g)$, respectively. 
For sets $\mathbf{X}$, we let $\pa(\mathbf{X},\g)=\bigcup_{X_i \in \mathbf{X}} \pa(X_i,\g)$, with analogous definitions for $\an(\mathbf{X},\g)$, $\de(\mathbf{X},\g),\possan(\mathbf{X},\g)$, $\possde(\mathbf{X},\g)$.

\vsp\noindent{\bf Colliders, definite status paths and v-structures}. A node $V$ is a  \emph{collider} on a path $p$ if $p$ contains a subpath $(U,V,W)$ such that $U \rightarrow V \leftarrow W$. A node $V$ is called a \emph{non-collider} on $p$ if $p$ contains a subpath $(U,V,W)$ such that 
(i) $U \leftarrow V$, or (ii) $V\to W$, or (iii) $U-V-W$ and $U$ and $W$ are not adjacent in $\g$.
A path $p$ is of \textit{definite status} if every non-endpoint node on $p$ is either a collider or a non-collider. 
%A definite status path $p = (Z_1, \dots Z_m)$ is possibly directed if $Z_i \leftarrow Z_{i+1}$ is not in $\g$ for any $i \in \{1,m-1\}$ (Lemma 3.5 in \citealp{perkovic17}).
If $U \rightarrow V \leftarrow W$ is in $\g$ and $U$ and $W$ are not adjacent in $\g$, then $(U,V,W)$ is called a v-structure in $\g$.

\vsp\noindent\textbf{Directed cycles, DAGs and PDAGs.} A directed path from $X$ to $Y$, together with the edge $Y\to X$ forms a \textit{directed cycle}. A directed graph without directed cycles is called a \textit{directed acyclic graph} (DAG) and a partially directed graph without directed cycles is called a \textit{partially directed acyclic graph} (PDAG). 
%For PDAG $\g$, we let $[\g]$ denote every DAG represented by $\g$.

\vsp\noindent \textbf{Blocking and d-separation in PDAGs.} (Cf.\  Definition 1.2.3 in \cite{pearl2009causality} and Definition 3.5 in \cite{maathuis2013generalized}). Let $\mathbf{Z}$ be a set of nodes in a PDAG. A definite status path $p$ is  \textit{blocked} by $\mathbf{Z}$ if (i) $p$ contains a non-collider that is in $\mathbf{Z}$, or (ii) $p$ contains a collider $C$ such that no descendant of $C$ is in $\mathbf{Z}$.
A definite status path that is not blocked by a set $\mathbf{Z}$ is \textit{open} given $\mathbf{Z}$. If $\mathbf{X},\mathbf{Y}$ and $\mathbf{Z}$ are three pairwise disjoint sets of nodes in a PDAG $\g$, then $\mathbf{Z}$  \textit{d-separates} $\mathbf{X}$ from $\mathbf{Y}$ in $\g$ if $\mathbf{Z}$ blocks every definite status path between any node in $\mathbf{X}$ and any node in $\mathbf{Y}$ in $\g$. 
We then write $\mathbf{X} \perp_{\g} \mathbf{Y}|\mathbf{Z}$. If  $\mathbf{Z}$ does not block every definite status path between any node in $\mathbf{X}$ and any node in $\mathbf{Y}$ in $\g$, we write $\mathbf{X} \not\perp_{\g} \mathbf{Y}|\mathbf{Z}$. 

\begin{Remark}
	We use the convention that for any two disjoint node sets $\mathbf{X}$ and $\mathbf{Y}$ it holds that $\emptyset \perp_{\g} \mathbf{X} | \mathbf{Y}$.
\end{Remark}

\vsp\noindent\textbf{Markov property and faithfulness.} (Cf.\ \citep[Definition 1.2.2][]{pearl2009causality}) Let $\mathbf{X}$, $\mathbf{Y}$ and $\mathbf{Z}$ be disjoint sets of random variables. We use the notation $\mathbf{X} \ci \mathbf{Y}|\mathbf{Z}$ to denote that $\mathbf{X}$ is conditionally independent of $\mathbf{Y}$ given $\mathbf{Z}$. A density $f$ is called \emph{Markov} with respect to a DAG $\g$ if  $\mathbf{X} \perp_{\g} \mathbf{Y} | \mathbf{Z}$ implies $\mathbf{X} \ci \mathbf{Y} | \mathbf{Z}$ in $f$. If this implication also holds in the other direction, then $f$ is \emph{faithful} with respect to $\g$.

\begin{figure}[t]
	
	\centering
	\subfloat[]
	{
		%{.22\textwidth}
		\centering
		\captionsetup[subfigure]{width=80pt}%
		\begin{tikzpicture}[->,>=latex,shorten >=2pt,auto,node distance=1.2cm,scale=1.3,transform shape]
		\tikzstyle{state}=[inner sep=0.5pt, minimum size=5pt]
		
		% rule 1
		\node[state] (Xia) at (0,0) {};
		\node[state] (Xka) at (0,1.2) {};
		\node[state] (Xja) at (1.2,0) {};
		
		\path (Xka) edge (Xia);
		\draw[-]
		(Xia) edge (Xja);
		\end{tikzpicture}
		%			\caption{}
		\label{subfig:forb1}
	}
\hspace{0.5cm}
	\subfloat[]{
		\centering
		\begin{tikzpicture}[->,>=latex,shorten >=1pt,auto,node distance=.2cm,scale=1.3,transform shape]
		\tikzstyle{state}=[inner sep=0.5pt, minimum size=5pt]
		
		% rule 2
		\node[state] (Xic) at (2.8,0) {};
		\node[state] (Xkc) at (2.8,1.2) {};
		\node[state] (Xjc) at (4,0) {};
		
		\path (Xic) edge (Xkc)
		(Xkc) edge (Xjc);
		\draw[-]
		(Xic) edge (Xjc);
		
		\end{tikzpicture}
		%			\caption{}
		%			\label{fig:forb2}
	}
\hspace{0.5cm}
	\subfloat[]{
		\centering
		\begin{tikzpicture}[->,>=latex,shorten >=1pt,auto,node distance=.2cm,scale=1.3,transform shape]
		\tikzstyle{state}=[inner sep=0.5pt, minimum size=5pt]

		% rule 3
		\node[state] (Xie) at (0,-1) {};
		\node[state] (Xke) at (1.2,-1) {};
		\node[state] (Xle) at (0,-2.2) {};
		\node[state] (Xje) at (1.2,-2.2) {};
		
		\path (Xke) edge (Xje)
		(Xle) edge (Xje);
		
		\draw[-]  (Xie) edge (Xje);
		\path[-]
		(Xke) edge (Xie)
		(Xle) edge (Xie);
		\end{tikzpicture}
	}
\hspace{0.5cm}
	\subfloat[]{
		\centering
		\begin{tikzpicture}[->,>=latex,shorten >=1pt,auto,node distance=.2cm,scale=1.3,transform shape]
		\tikzstyle{state}=[inner sep=0.5pt, minimum size=5pt]
		
		% rule 4
		\node[state] (Xig) at (2.8,-1) {};
		\node[state] (Xjg) at (4,-1) {};
		\node[state] (Xkg) at (2.8,-2.2) {};
		\node[state] (Xlg) at (4,-2.2) {};
		
		\path (Xlg) edge (Xkg)
		(Xjg) edge (Xlg);
		\draw[-]
		(Xig) edge (Xkg);
		\path[-]
		(Xig) edge (Xjg)
		(Xig) edge (Xlg);
		\end{tikzpicture}
		%	\subcaption{}
		\label{subfig:forb4}
	}
	\caption{Figures (a) - (d) each show one of the four forbidden induced subgraphs of a maximal PDAG \citep[see orientation rules in][]{meek1995causal}.}
	\label{fig:orientationRules}
\end{figure}
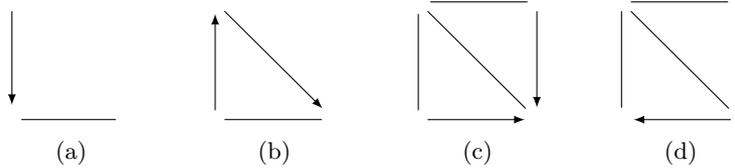

\vsp\noindent\textbf{Markov equivalence and CPDAGs.}  \citep[Cf.\ ][]{meek1995causal,andersson1997characterization} All DAGs that encode the same d-separation relationships are called \textit{Markov equivalent} and form a \textit{Markov equivalence class} of DAGs, which can be represented by a \textit{completed partially directed acyclic graph} (CPDAG). For any CPDAG $\g[C]$, we let $[\g[C]]$ denote the Markov equivalence class of DAGs represented by $\g[C]$. Conversely, if a DAG $\g[D]$ is in $[\g[C]]$, then $\g[C]$ is called the CPDAG of $\g[D]$.
A CPDAG $\g[C]$ has the same adjacencies and v-structures as any DAG in $[\g[C]]$. 
Moreover, a directed edge $X\to Y$ in $\g[C]$ corresponds to a directed edge $X \to Y$ in every
DAG in $[\g[C]]$, and for any undirected edge $X -Y$  in $\g[C]$, $[\g[C]]$ contains a DAG with $X \rightarrow Y$ and a DAG with $X \leftarrow Y$. 

\vsp\noindent\textbf{\Mpdag{}s.} A PDAG $\g$ is \textit{maximally oriented} (\mpdag{}) if and only if the graphs in Figure \ref{fig:orientationRules} are \textbf{not} induced subgraphs of $\g$. 

In general, a \mpdag{} $\g$ describes a subset of a Markov equivalence class of DAGs, denoted by $[\g]$.  
A \mpdag{} $\g$ has the same adjacencies and v-structures as any DAG in $[\g]$. 
Moreover, a directed edge $X\to Y$ in $\g$ corresponds to a directed edge $X \to Y$ in every
DAG in $[\g]$, and for any undirected edge $X -Y$  in $\g$, $[\g]$ contains a DAG with $X \rightarrow Y$ and a DAG with $X \leftarrow Y$.

%\begin{Example}
%	Consider the DAG $\g[D]$ in Figure~\ref{fig:exintropdag}(d) from the main paper. The CPDAG $\g[C]$ of $\g[D]$ is given in Figure~\ref{fig:exintropdag}(a) and two \mpdag{}s $\g$ and $\g'$ of $\g[D]$ are given in  Figures \ref{fig:exintropdag}(b) and \ref{fig:exintropdag}(c) also from the main paper.
%	
%	
%	%In general, for a path $p$ to be possibly directed from $X$ to $Y$ in a CPDAG, it is enough to check that no edge on $p$ is directed towards $X$ \citep[Lemma~7.5 in][]{maathuis2013generalized}. In \mpdag{}s, this reasoning only holds for definite status paths \citep[see Lemma~3.5 in][]{perkovic17}.
%	\label{ex: possibly directed}
%\end{Example}
 
 \begin{Example}\label{ex: CPDAGs, PDAGs and possibly directed}    
 	Consider the DAG $\g[D]$ in Figure~\ref{fig:exintropdag}(d). The CPDAG $\g[C]$ of $\g[D]$ is given in Figure~\ref{fig:exintropdag}(a) and two \mpdag{}s $\g$ and $\g'$ of $\g[D]$ are given in  Figures \ref{fig:exintropdag}(b) and \ref{fig:exintropdag}(c). The CPDAG $\g[C]$ represents 8 DAGs, the \mpdag{} $\g$ represents five DAGs and the \mpdag{} $\g'$ represents two DAGs (see Figure \ref{fig:allDAGs}). Moreover, $\g[D] \in [\g'] \subseteq [\g] \subseteq [\g[C]]$, illustrating that $\g'$ and $\g$ represent refinements of the Markov equivalence class $[\g[C]]$.
 	
 	We consider now possibly directed paths in $\g[C]$ and $\g$. According to our definition the path $V_3-V_4-V_1$ is possibly directed from $V_3$ to $V_1$ in $\g[C]$, but not in $\g$, since $\g$ contains the edge $V_1\rightarrow V_3$. As a result, $V_1 \in \possde(V_3,\g[C])$ but $V_1 \notin \possde(V_3,\g)$. The rationale behind these definitions is that there is a DAG in $[\g[C]]$ containing $V_3\to V_4\to V_1$, but there is no such DAG in $[\g]$ (see Figure \ref{fig:allDAGs}).

 	%	We now consider possibly directed paths in $\g[C]$ and $\g$. According to our definition, as given in Section \ref{S-graph-supp} of the Supplement \cite{supplement} the path $V_3-V_4-V_1$ is possibly directed from $V_3$ to $V_1$ in $\g[C]$, but not in $\g$, since $\g$ contains the edge $V_1\rightarrow V_3$. As a result, $V_1 \in \possde(V_3,\g[C])$ but $V_1 \notin \possde(V_3,\g)$. The rationale behind these definitions is that there is a DAG in $[\g[C]]$ containing $V_3\to V_4\to V_1$, but there is no such DAG in $[\g]$ (see Figure \ref{S-fig:allDAGs} in the Supplement \citep{supplement}).
 	
 	%In general, for a path $p$ to be possibly directed from $X$ to $Y$ in a CPDAG, it is enough to check that no edge on $p$ is directed towards $X$ \citep[Lemma~7.5 in][]{maathuis2013generalized}. In \mpdag{}s, this reasoning only holds for definite status paths \citep[see Lemma~3.5 in][]{perkovic17}.
 \end{Example}

 \begin{figure}[t]
 	
 	\begin{center}
  		\begin{tikzpicture}[>=stealth',shorten >=1pt,auto,node distance=0.8cm,scale=0.8, transform shape,align=center,minimum size=3em]
 \node[state] (x) at (0,0) {$X$};
 \node[state] (z) [right =of x] {$M$};
 \node[state] (y) [right =of z] {$Y$};
 
 \node[state] (a1) at ($(x)!0.5!(z)+(0,+1.5)$) {$A_1$}; 
 \node[state] (a2) [right =of a1] {$A_2$};
 \node[state] (r) [right =of a2] {$R$};
 \node[state] (v) [left =of a1] {$V$};
 
 \node[state] (b1) at ($(x)!0.5!(z)+(0,-1.5)$) {$B_1$}; 
 \node[state] (b2) [right =of b1] {$B_2$};
 \node[state] (f) [right =of b2] {$F$};
 \node[state] (d) [left =of b1] {$D$}; 
 %	{\{x, y\}};

 %	\path (x) edge (y);
 \path[->]   (x) edge    (z);
 \path[->]   (z) edge    (y);
 %\path[<->]   (x) edge  [dashed, bend left]  (y);
 \path[->]   (a1) edge     (x);
 \path[->]   (a1) edge     (a2);
 \path[->]   (a2) edge     (y);
 \path[->]   (b1) edge     (x);
 \path[->]   (b2) edge     (b1);
 \path[->]   (b2) edge     (y);
 \path[->]   (r) edge     (y);
 \path[->]   (y) edge     (f);
 \path[->]   (v) edge     (x);
 \path[->]   (x) edge     (d);
 
 \end{tikzpicture}
 	\end{center}
 	\caption{Causal DAG from Examples \ref{ex: introduction figure explained} and \ref{section:example}}
 %		Causal DAG $\g$
 	\label{AllinOne}
 \end{figure}
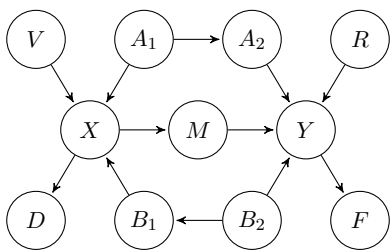
 
 \noindent \textbf{Causal, non-causal and possibly causal paths and nodes} A directed path from $X$ to $Y$ in a causal graph is also called a \textit{causal path} from $X$ to $Y$. Analogously, a possibly directed path from $X$ to $Y$ is called a \textit{possibly causal path}. 
 %is a sequence of distinct nodes $(X=Z_1, \dots, Z_m=Y)$, $m \ge 1$, such that $Z_i \to Z_{i+1}$ is in $\g$ for all $i \in \{1, \dots, m-1\}$.  
 %A \textit{possibly causal path} from $X$ to $Y$ is a path from $X$ to $Y$ where no edge points towards $X$. 
 A \textit{non-causal path} from $X$ to $Y$ is a path that is not possibly directed from $X$ to $Y$.
 Let $\mathbf{X}$ and $\mathbf{Y}$ be disjoint node sets in a causal \mpdag{} $\g$. We define \textit{causal nodes} relative to $(\mathbf{X,Y})$ in $\g$, denoted $\cnb{\g}$, as all nodes on proper causal paths from $\mathbf{X}$ to $\mathbf{Y}$, excluding nodes in $\mathbf{X}$. 
 For singleton $X$ the causal nodes are also called the mediating nodes. 
 Analogously, we define
 \textit{possible causal nodes} relative to $(\mathbf{X,Y})$ in $\g$, denoted $\posscnb{\g}$, as all nodes on proper possibly causal paths from $\mathbf{X}$ to $\mathbf{Y}$, excluding nodes in $\mathbf{X}$. 
 % \textit{Possible causal nodes} relative to $(\mathbf{X}$, $\mathbf{Y})$ in $\g$ are all nodes on proper possibly causal paths from $\mathbf{X}$ to $\mathbf{Y}$, excluding nodes in $\mathbf{X}$.
 % Causal and possibly causal nodes relative to $(\mathbf{X,Y})$ and $\g$ are denoted by $\cnb{\g}$ and $\posscnb{\g}$ respectively.  

 \vsp\noindent{\bf Forbidden nodes}. \citep[Cf.][]{perkovic16} Let $\mathbf{X}$ and $\mathbf{Y}$ be disjoint node sets in a \mpdag{} $\g$. We define \textit{forbidden nodes} relative to $(\mathbf{X},\mathbf{Y})$ in $\g$ as
 \[
 \fb{\g} = \possde(\posscnb{\g}, \g) \cup \mathbf{X}.
 \]

 %\begin{Example}  Consider the causal DAG $\g[D]$ in Figure~\ref{fig:exintropdag}(c). Since $\g[D]$ is a DAG, $\posscn{\g[D]} = \cn{\g[D]}=\{Y\}$ and $\possde(Y,\g[D]) = \de(Y,\g[D]) = \{Y,Z\}$. Hence, 
 %and all paths in $\g[D]$ are of definite status. 
 %$\f{\g[D]} = \{X,Y,Z\}$.
 
 %Now consider the causal CPDAG $\g[C]$ of $\g[D]$ in Figure~\ref{fig:exintropdag}(a) and the causal \mpdag{} $\g$ of $\g[D]$ in  Figure~\ref{fig:exintropdag}(b). In $\g[C]$,  $\posscn{\g[C]} = \{Z,Y\}$ and $\f{\g[C]} =  \{X,Y,Z\}$. 
 %Also, $\de(X,\g[C]) = \{X\}$ and $\possde(X,\g[C]) = \{X,Z,Y\}$.
 %Similarly, one can easily see that $\fb{\g} = \{X,Y,Z\}. 
 %In $\g$, $\posscn{\g[C]}  = \{Y\}$, since $X - Z- Y $ is a non-causal path in $\g$. However, $\f{\g} =  \{X,Y,Z\}$, since $Z$ is a possible descendant of $Y$ in $\g$. 
 %Additionally, $\de(X,\g) = \{X,Y\}$ and $\possde(X,\g) = \{X,Z,Y\}$.
 %\end{Example}
 
 \vsp\noindent{\bf Amenability}. \citep{perkovic16}  Let $\mathbf{X}$ and $\mathbf{Y}$ be disjoint node sets in a \mpdag{} $\g$. If all proper possibly causal paths from $\mathbf{X}$ to $\mathbf{Y}$ start with a directed edge out of $\mathbf{X}$, then we call $\g$ \textit{amenable} relative to $(\mathbf{X},\mathbf{Y})$.

 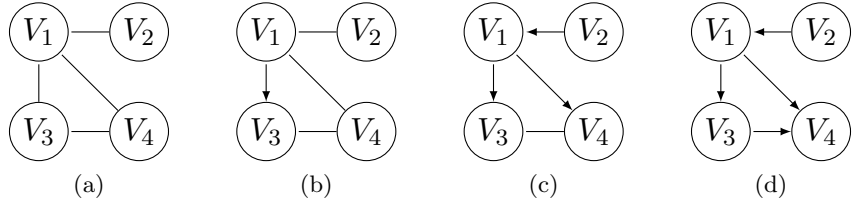
\begin{figure}
 	\centering
 	\subfloat[]{
 		\centering
 		\begin{tikzpicture}[->,>=latex,shorten >=1pt,auto,node distance=1.2cm,scale=1.1,transform shape]
 		\tikzstyle{vertex}=[circle, draw, inner sep=2pt, minimum size=0.15cm]
 		\node[vertex]         (X)                        {$V_1$};
 		\node[vertex]         (V2) [right of= X]  		{$V_2$};
 		\node[vertex]         (V1) [below of = V2] 		{$V_4$};
 		\node[vertex]       	 (Y)  [below of= X] 		{$V_3$};
 		\draw[-] 	(V1)	edge    (X);
 		\draw[-] 	(V2) 	edge    (X);
 		\draw[-] 	(Y) 	edge    (X);
 		\draw[-] 	(V1) 	edge    (Y);
 		\end{tikzpicture}
 		%		\caption{}
 		%		\label{fig:cpdagintro}
 	}
 	%	\unskip
 	\hspace{.5cm}
 	\subfloat[]{
 		\centering
 		\begin{tikzpicture}[->,>=latex,shorten >=1pt,auto,node distance=1.2cm,scale=1.1,transform shape]
 		\tikzstyle{vertex}=[circle, draw, inner sep=2pt, minimum size=0.15cm]
 		\node[vertex]         (X)                        {$V_1$};
 		\node[vertex]         (V2) [right of= X]  		{$V_2$};
 		\node[vertex]         (V1) [below of = V2] 		{$V_4$};
 		\node[vertex]       	 (Y)  [below of= X] 		{$V_3$};
 		\draw[-] 	(V1) 	edge    (X);
 		\draw[-] 	(V2) 	edge    (X);
 		\draw[<-] 	(Y) 	edge    (X);
 		\draw[-] 	(V1) 	edge    (Y);
 		\end{tikzpicture}
 		%		\caption{}
 		%		\label{fig:pdagintro}
 	}
 	\hspace{.5cm}
 	%	\unskip 
 	%	\vrule
 	\subfloat[]{
 		\centering
 		\begin{tikzpicture}[->,>=latex,shorten >=1pt,auto,node distance=1.2cm,scale=1.1,transform shape]
 		\tikzstyle{vertex}=[circle, draw, inner sep=2pt, minimum size=0.15cm]
 		\node[vertex]         (X)                        {$V_1$};
 		\node[vertex]         (V2) [right of= X]  		{$V_2$};
 		\node[vertex]         (V1) [below of = V2] 		{$V_4$};
 		\node[vertex]       	 (Y)  [below of= X] 		{$V_3$};
 		\draw[<-] 	(V1) 	edge    (X);
 		\draw[->] 	(V2) 	edge    (X);
 		\draw[<-] 	(Y) 	edge    (X);
 		\draw[-] 	(V1) 	edge    (Y);
 		\end{tikzpicture}
 		%		\caption{}
 		%		\label{fig:pdag2}
 	}
 	\hspace{.5cm}
 	%	\unskip
 	%	\vrule
 	\subfloat[]{
 		\centering
 		\begin{tikzpicture}[->,>=latex,shorten >=1pt,auto,node distance=1.2cm,scale=1.1,transform shape]
 		\tikzstyle{vertex}=[circle, draw, inner sep=2pt, minimum size=0.15cm]
 		\node[vertex]         (X)                        {$V_1$};
 		\node[vertex]         (V2) [right of= X]  		{$V_2$};
 		\node[vertex]         (V1) [below of = V2] 		{$V_4$};
 		\node[vertex]       	 (Y)  [below of= X] 		{$V_3$};
 		\draw[<-] 	(V1) 	edge    (X);
 		\draw[->] 	(V2) 	edge    (X);
 		\draw[<-] 	(Y) 	edge    (X);
 		\draw[<-] 	(V1) 	edge    (Y);
 		\end{tikzpicture}
 		%		\caption{}
 		%		\label{fig:dagintro}
 	}
 	\caption{(a) CPDAG, (b) \mpdag{}, (c) \mpdag{} and  (d) DAG  from Example \ref{ex: CPDAGs, PDAGs and possibly directed}}
% 	, \ref{ex:total effect} and \ref{emma:cpdag}	.}
 	
 	\label{fig:exintropdag}
 \end{figure} 

 \begin{Definition}(\textbf{Generalized adjustment criterion}; \citealp[Def.\ 4.3 of][]{perkovic17})
 	\label{def: generalized adjustment criterion}
 	Let $\mathbf{X},\mathbf{Y}$ and $\mathbf{Z}$ be pairwise disjoint node sets in a causal \mpdag{} $\g$. Then $\mathbf{Z}$ satisfies the generalized adjustment criterion relative to $(\mathbf{X,Y})$ in~$\g$ if the following three conditions hold:
 	\begin{enumerate}
 		%		[label = (\cctext*), leftmargin=0.5cm,align=left]
 		\item\label{cond0} The graph $\g$ is amenable relative to $(\mathbf{X},\mathbf{Y})$. 
 		\item\label{cond1} $\mathbf{Z} \cap \fb{\g} = \emptyset$, and
 		\item\label{cond2} all proper non-causal definite status paths from $\mathbf{X}$ to $\mathbf{Y}$ are blocked by $\mathbf{Z}$.
 	\end{enumerate}  \label{adjustment}
 \end{Definition}

\begin{Theorem} (\citealp[Theorem~4.4 in][]{perkovic17})
	Let $\mathbf{X},\mathbf{Y}$ and $\mathbf{Z}$ be pairwise disjoint node sets in a causal \mpdag{} $\g$. Then 
	\begin{equation}
	f(\mathbf{y}|do(\mathbf{x}))=
	\begin{cases}
	f(\mathbf{y}|\mathbf{x}) & \text{if }\mathbf{Z} = \emptyset,\\
	\int_{\mathbf{z}}f(\mathbf{y}|\mathbf{x,z})f(\mathbf{z})d\mathbf{z}  & \text{otherwise.}
	\end{cases}
	\label{defadjustmentmpdag}
	\end{equation}
	for any density $f$ compatible with $\g$ if and only if $\mathbf{Z}$ satisfies the generalized adjustment criterion relative to $(\mathbf{X,Y})$ in $\g$ (Def.\ \ref{def: generalized adjustment criterion}).
	\label{thm:adjustment-set}
\end{Theorem}

	Theorem \ref{thm:adjustment-set} establishes that Definition \ref{def: generalized adjustment criterion} characterizes all covariate sets that can be used for causal effect estimation via adjustment. 
	%	There are analogues versions of this Theorem for MAGs and PAGs in \cite{perkovic2015complete, perkovic16}. 
	It is as a consequence of this Theorem that we refer to sets that satisfy the generalized adjustment criterion as \emph{valid adjustment sets}.

\begin{Example}
	\label{ex: introduction figure explained}
	Consider the DAG $\g$ in Figure \ref{AllinOne}. Since $\g$ is a DAG, it is trivially amenable and 
	\[
	\f{\g} = \de(\cn{\g},\g) \cup \{X\} = \de(\{M,Y\},\g) \cup \{X\}= \{X,M,Y,F\}.
	\]
	Further, any valid adjustment set needs to contain at least one node from $\{A_1,A_2\}$ and one node from $\{B_1,B_2\}$ to satisfy the blocking criterion. The remaining nodes $V,D,R$ are neither required nor forbidden. This shows that any valid adjustment set has to be of the following form:  
	$\mathbf{Z} = \mathbf{A} \cup \mathbf{B} \cup \mathbf{C}$, where $\mathbf{A} \subseteq \{A_1,A_2\}$ and $\mathbf{B} \subseteq \{B_1, B_2\}$ are non empty and $\mathbf{C} \subseteq \{V,D,R \}$ is possibly empty. 
	
	As an example of a joint intervention let $\mathbf{X} = \{X,A_2\}$ and $\mathbf{Y}=\{Y,F\}$. The amenability follows trivially from the fact that $\g$ is a DAG and
	\begin{align*}
	\fb{\g} &= \de(\cnb{\g},\g) \cup \{\mathbf{X}\} \\ 
	&= \de(\{M,Y,F\},\g) \cup \{\mathbf{X}\}= \{X,A_2,M,Y,F\}. 
	\end{align*}
	Further,  any valid adjustment set here must block the two proper non-causal paths $(X,B_1,B_2,Y)$ and $(X,B_1,B_2,Y,F)$ from $\mathbf{X}$ to $\mathbf{Y}$. Thus, any valid adjustment set has to be of the following form:  
	$\mathbf{Z} = \mathbf{B} \cup \mathbf{C}$, where $\mathbf{B} \subseteq \{B_1, B_2\}$ is non empty and $\mathbf{C} \subseteq \{A_1,V,D,R \}$ is possibly empty. 
\end{Example}

\vsp

\noindent{}\textbf{Unshielded paths, corresponding paths and path concatenation.} A path \linebreak $(V_{i},V_{j},V_{k})$ in a partially directed graph $\g$ is an \emph{unshielded triple} if $ V_{i} $ and $ V_{k}$ are not adjacent in $\g$. A path is \textit{unshielded} if all successive triples on the path are unshielded. 
%Every unshielded path is a path of definite status. 
If $\g$ and $\g^*$ are two graphs with identical adjacencies and $p$ is a path in $\g$, then the \textit{corresponding path} $\pstar$ in $\g^*$ consists of the same node sequence as $p$.  We denote the concatenation of paths by $\oplus$. For example given a path $p=(X_1, \dots, X_k)$ it holds that $p = p(X_1,X_{m}) \oplus p(X_{m},X_{k})$ for $1 \le m \le k$. In general, concatenating paths does not result in a path, but we only use the symbol $\oplus$ if the result is in fact a path. 

\vsp

\noindent{}\textbf{Partial total effects in a causal linear model.} Consider a causal DAG $\g=(\mathbf{V},\mathbf{E})$, such that $\mathbf{V}$ follows a causal linear model compatible with $\g$. Let $\mathbf{Z}$ be a node set and $X,Y \notin \mathbf{Z}$ be two nodes in $\g$. Then the \emph{partial total effect} $\tau_{yx.\mathbf{z}}$ of $X$ on $Y$ given $\mathbf{Z}$ is defined as the sum of the total effect along all causal paths from $X$ to $Y$ that do not contain nodes from $\mathbf{Z}$.

\vsp

We now give a small simulation study to illustrate that Proposition \ref{avarDAG} also covers regressions with heteroskedastic residuals and how these arise in causal linear models with non-Gaussian errors. For an example of non-linearity we refer to Example 1 from the supplement of \citet[][]{nandy2017estimating}.

\begin{Example}
	\label{section:example}
	
	\begin{figure}[t]
		\subfloat{
			\includegraphics[scale=.38]{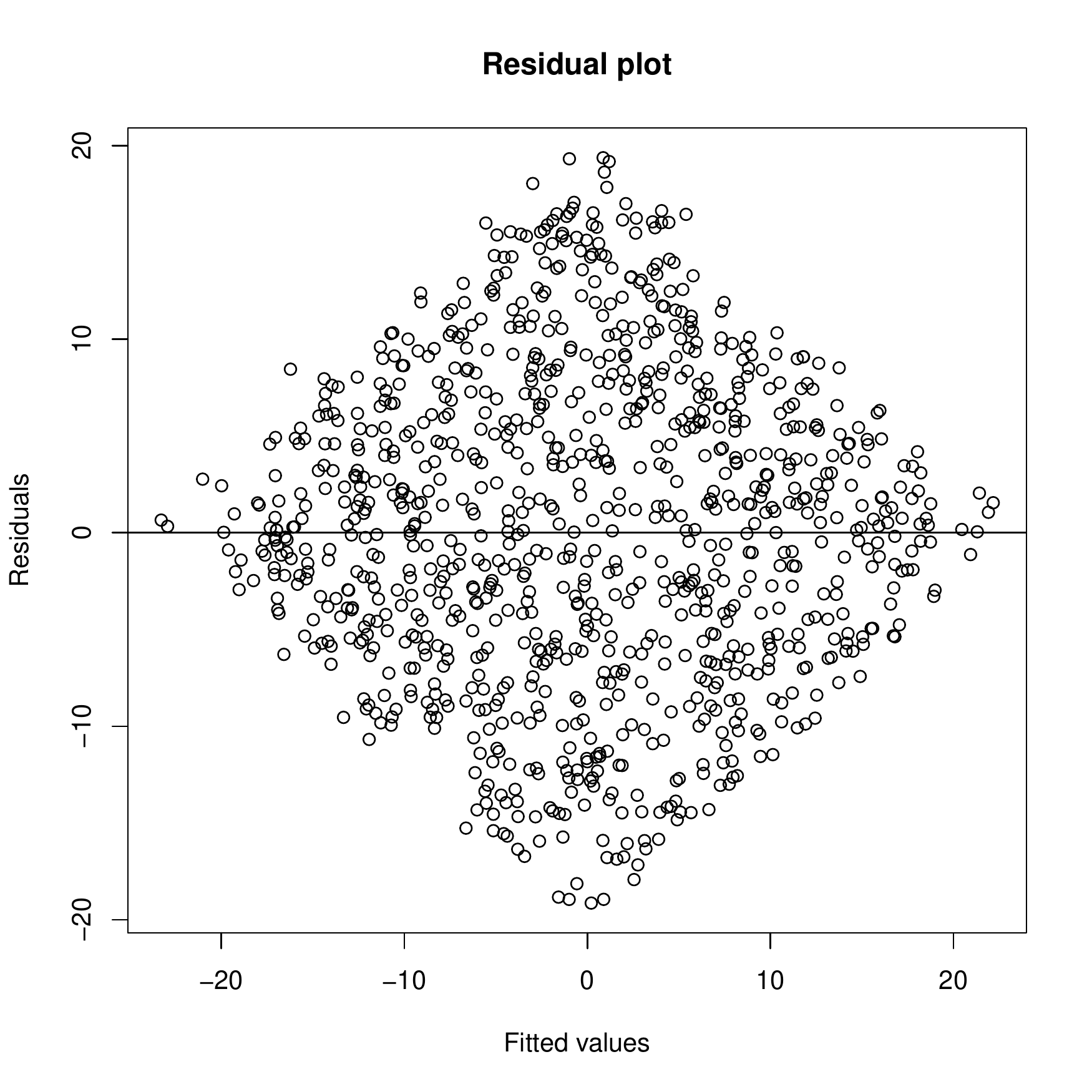}
		}
		\subfloat{
			\includegraphics[scale=.38]{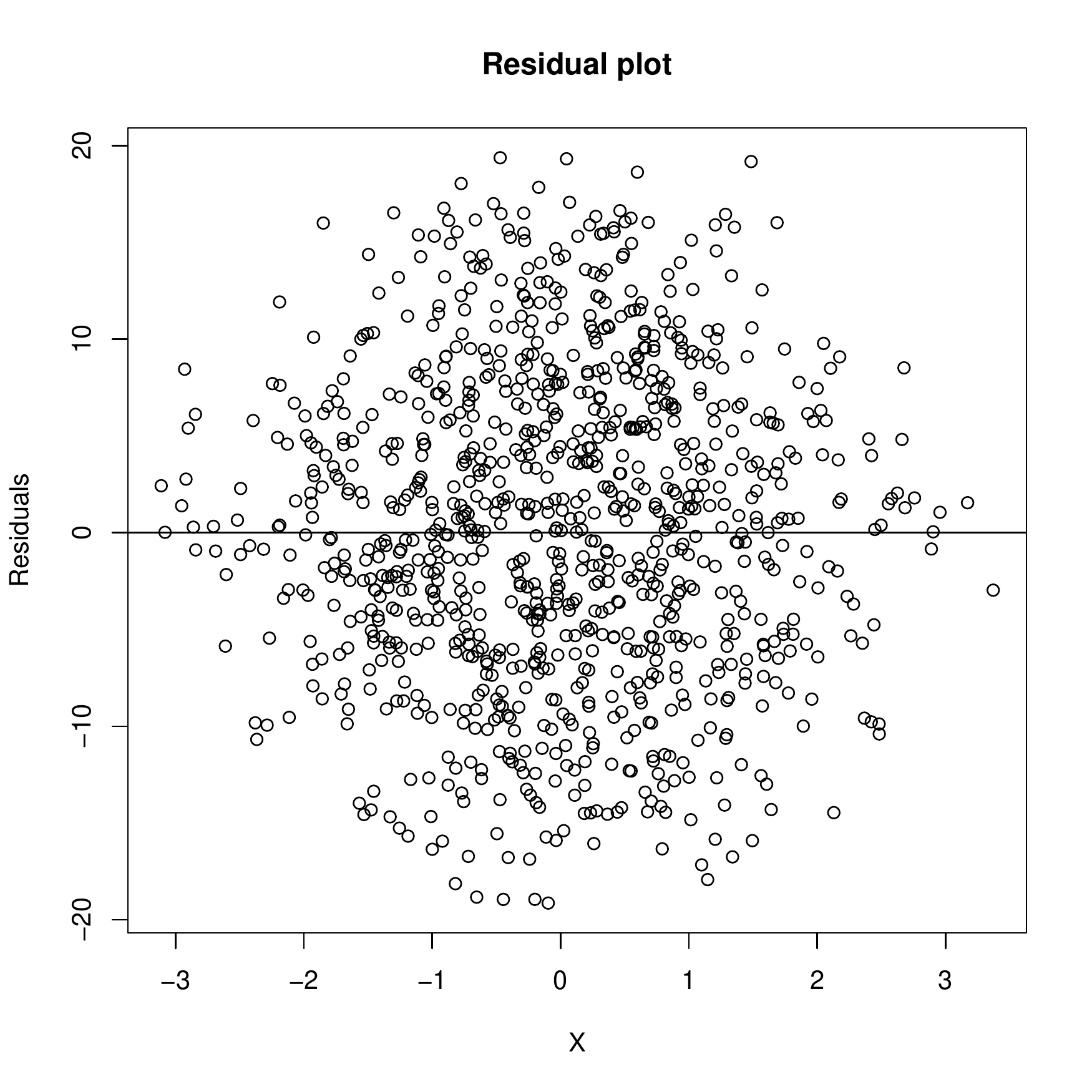}
		}
		\label{figure:residuals}
		\caption{Residual vs Fitted and Residual vs X plot for the example OLS regression considered in Example \ref{section:example}.}
	\end{figure}
	
	Consider the DAG $\g$ given in Figure \ref{AllinOne}, with for simplicity the nodes $V,D,R,F$ dropped. We sampled data from the causal linear model compatible with $\g$, such that all errors are uniformly distributed on $[-1,1]$ and all edge coefficients except $\alpha_{yb_2}$ are 1 with $\alpha_{yb_2}=20$. We used this data to estimate the total effect of $X$ on $Y$, i.e. $1$, by adjusting for the valid adjustment set $\{A_2,B_1\}$. More precisely, we computed the ordinary least squares (OLS) regression of $Y$ on $X,A_2,B_2$ for one generated data set of 1000 points. 
	Figure \ref{figure:residuals} shows the resulting Residuals vs Fitted plot as well as the Residuals vs X plot. As we can see the residuals are clearly heteroskedastic with their distribution depending on both the fitted values and X.

	To verify that Proposition \ref{avarDAG} nonetheless holds, we then repeated this procedure 1000 times collecting the estimated coefficients. Based on these we computed empirical means and variances for all three coefficients in the regression of $Y$ on $X,A_2,B_1$. We also computed theoretical asymptotic variances in accordance with the formula given in Proposition \ref{avarDAG} from the known true underlying causal linear model. The results are given in Table \ref{table:empirical}.
	
	\begin{table}
		\caption{Empirical mean and scaled to sample size variance from the simulation study in Section \ref{section:example} as well as the theoretical asymptotic variance according to the formula from Proposition \ref{avarDAG} for all three coefficients in the considered OLS regression of $Y$ on $X,A_2,B_1$.}
		\centering
		\begin{tabular}{r|rrr}
			& Emp. mean & Emp. variance & Asy. variance \\ 
			\hline
			$\hat{\beta}_{yx.a_2b_1}$ & 1.01 & 0.13 & 0.13 \\ 
			$\hat{\beta}_{ya_2.xb_1}$ & 1.00 & 0.13 & 0.13 \\ 
			$\hat{\beta}_{yb_1.xa_2}$ & 9.98 & 0.16 & 0.23 \\ 
		\end{tabular}
		\label{table:empirical}
	\end{table}
	
	We see that in fact $\hat{\beta}_{yx.a_2b_1}$ appears to estimate the total effect $1$ and that it does so with the claimed asymptotic variance from Proposition \ref{avarDAG}. Interestingly the same is true for $\hat{\beta}_{ya_2.xb_1}$, which is due to the fact that $B_1$ is a valid adjustment set for the joint effect of $\{X,A_2\}$ on $Y$ (see Example \ref{ex: introduction figure explained}) and therefore Proposition \ref{avarDAG} also holds with respect to $\hat{\beta}_{ya_2.xb_1}$. Interestingly, the theoretical and empirical variances do not match for $\hat{\beta}_{yb_1.xa_2}$. This is due to the fact that given $\{X,A_2\}$ the non-causal path $B_1 \leftarrow B_2 \rightarrow Y$ remains open and therefore $\hat{\beta}_{yx.a_2b_1}$ is not covered by Proposition \ref{avarDAG}. This illustrates that the asymptotic variance result from Proposition \ref{avarDAG} does not necessarily cover all coefficients in the considered OLS regression.
	
\end{Example}

\begin{figure}
	\centering
	\subfloat[]{
		\centering
		\begin{tikzpicture}[->,>=latex,shorten >=1pt,auto,node distance=1.2cm,scale=1.1,transform shape]
		\tikzstyle{vertex}=[circle, draw, inner sep=2pt, minimum size=0.15cm]
		\node[vertex]         (V1)                        {$V_1$};
		\node[vertex]         (V2) [right of= V1]  		{$V_2$};
		\node[vertex]         (V3) [below of = V2] 		{$V_4$};
		\node[vertex]       	 (V4)  [below of= V1] 		{$V_3$};
		\draw[<-] 	(V1)	edge    (V2);
		\draw[->] 	(V1) 	edge    (V3);
		\draw[->] 	(V1) 	edge    (V4);
		\draw[->] 	(V3) 	edge    (V4);
		\end{tikzpicture}
		%		\subcaption{}
		%	\end{subfigure}
	}	\unskip 
	%  \vrule
	\subfloat[]{
		\centering
		\begin{tikzpicture}[->,>=latex,shorten >=1pt,auto,node distance=1.2cm,scale=1.1,transform shape]
		\tikzstyle{vertex}=[circle, draw, inner sep=2pt, minimum size=0.15cm]
		\node[vertex]         (V1)                        {$V_1$};
		\node[vertex]         (V2) [right of= V1]  		{$V_2$};
		\node[vertex]         (V3) [below of = V2] 		{$V_4$};
		\node[vertex]       	 (V4)  [below of= V1] 		{$V_3$};
		\draw[<-] 	(V1)	edge    (V2);
		\draw[->] 	(V1) 	edge    (V3);
		\draw[->] 	(V1) 	edge    (V4);
		\draw[<-] 	(V3) 	edge    (V4);
		\end{tikzpicture}
		%		\subcaption{}
	}	\unskip 
	%\vrule
	\subfloat[]{
		\centering
		\begin{tikzpicture}[->,>=latex,shorten >=1pt,auto,node distance=1.2cm,scale=1.1,transform shape]
		\tikzstyle{vertex}=[circle, draw, inner sep=2pt, minimum size=0.15cm]
		\node[vertex]         (V1)                        {$V_1$};
		\node[vertex]         (V2) [right of= V1]  		{$V_2$};
		\node[vertex]         (V3) [below of = V2] 		{$V_4$};
		\node[vertex]       	 (V4)  [below of= V1] 		{$V_3$};
		\draw[->] 	(V1)	edge    (V2);
		\draw[->] 	(V1) 	edge    (V3);
		\draw[->] 	(V1) 	edge    (V4);
		\draw[->] 	(V3) 	edge    (V4);
		\end{tikzpicture}
		%		\subcaption{}
	} \unskip	
	%\vrule
	\subfloat[]{
		\centering
		\begin{tikzpicture}[->,>=latex,shorten >=1pt,auto,node distance=1.2cm,scale=1.1,transform shape]
		\tikzstyle{vertex}=[circle, draw, inner sep=2pt, minimum size=0.15cm]
		\node[vertex]         (V1)                        {$V_1$};
		\node[vertex]         (V2) [right of= V1]  		{$V_2$};
		\node[vertex]         (V3) [below of = V2] 		{$V_4$};
		\node[vertex]       	 (V4)  [below of= V1] 		{$V_3$};
		\draw[->] 	(V1)	edge    (V2);
		\draw[->] 	(V1) 	edge    (V3);
		\draw[->] 	(V1) 	edge    (V4);
		\draw[<-] 	(V3) 	edge    (V4);
		\end{tikzpicture}
		%		\subcaption{}
	} \par
	%\bigskip
	% #5
	\subfloat[]{
		\centering
		\begin{tikzpicture}[->,>=latex,shorten >=1pt,auto,node distance=1.2cm,scale=1.1,transform shape]
		\tikzstyle{vertex}=[circle, draw, inner sep=2pt, minimum size=0.15cm]
		\node[vertex]         (V1)                        {$V_1$};
		\node[vertex]         (V2) [right of= V1]  		{$V_2$};
		\node[vertex]         (V3) [below of = V1] 		{$V_3$};
		\node[vertex]       	 (V4)  [below of= V2] 		{$V_4$};
		\draw[->] 	(V1)	edge    (V2);
		\draw[->] 	(V1) 	edge    (V3);
		\draw[<-] 	(V1) 	edge    (V4);
		\draw[<-] 	(V3) 	edge    (V4);
		\end{tikzpicture}
		%		\subcaption{}
	}	\unskip	
	%\vrule
	\subfloat[]{
		\centering
		\begin{tikzpicture}[->,>=latex,shorten >=1pt,auto,node distance=1.2cm,scale=1.1,transform shape]
		\tikzstyle{vertex}=[circle, draw, inner sep=2pt, minimum size=0.15cm]
		\node[vertex]         (V1)                        {$V_1$};
		\node[vertex]         (V2) [right of= V1]  		{$V_2$};
		\node[vertex]         (V3) [below of = V1] 		{$V_3$};
		\node[vertex]       	 (V4)  [below of= V2] 		{$V_4$};
		\draw[->] 	(V1)	edge    (V2);
		\draw[<-] 	(V1) 	edge    (V3);
		\draw[<-] 	(V1) 	edge    (V4);
		\draw[<-] 	(V3) 	edge    (V4);
		\end{tikzpicture}
		%		\subcaption{}
	}	\unskip	
	%\vrule
	\subfloat[]{
		\centering
		\begin{tikzpicture}[->,>=latex,shorten >=1pt,auto,node distance=1.2cm,scale=1.1,transform shape]
		\tikzstyle{vertex}=[circle, draw, inner sep=2pt, minimum size=0.15cm]
		\node[vertex]         (V1)                        {$V_1$};
		\node[vertex]         (V2) [right of= V1]  		{$V_2$};
		\node[vertex]         (V3) [below of = V1] 		{$V_3$};
		\node[vertex]       	 (V4)  [below of= V2] 		{$V_4$};
		\draw[->] 	(V1)	edge    (V2);
		\draw[<-] 	(V1) 	edge    (V3);
		\draw[->] 	(V1) 	edge    (V4);
		\draw[->] 	(V3) 	edge    (V4);
		\end{tikzpicture}
		%		\subcaption{}
	}	\unskip	
	%\vrule
	\subfloat[]{
		\centering
		\begin{tikzpicture}[->,>=latex,shorten >=1pt,auto,node distance=1.2cm,scale=1.1,transform shape]
		\tikzstyle{vertex}=[circle, draw, inner sep=2pt, minimum size=0.15cm]
		\node[vertex]         (V1)                        {$V_1$};
		\node[vertex]         (V2) [right of= V1]  		{$V_2$};
		\node[vertex]         (V3) [below of = V1] 		{$V_3$};
		\node[vertex]       	 (V4)  [below of= V2] 		{$V_4$};
		\draw[->] 	(V1)	edge    (V2);
		\draw[<-] 	(V1) 	edge    (V3);
		\draw[<-] 	(V1) 	edge    (V4);
		\draw[->] 	(V3) 	edge    (V4);
		\end{tikzpicture}
		%		\subcaption{}
	}
	\caption{These are the 8 DAGs that form the Markov equivalence class represented by the CPDAG in Figure \ref{fig:exintropdag}(a). The DAGs (a)-(e) form the equivalence class represented by the \mpdag{} in Figure \ref{fig:exintropdag}(b).
		The DAGs (a)-(b) form the equivalence class represented by the \mpdag{} in Figure \ref{fig:exintropdag}(c).}
	\label{fig:allDAGs}
\end{figure}
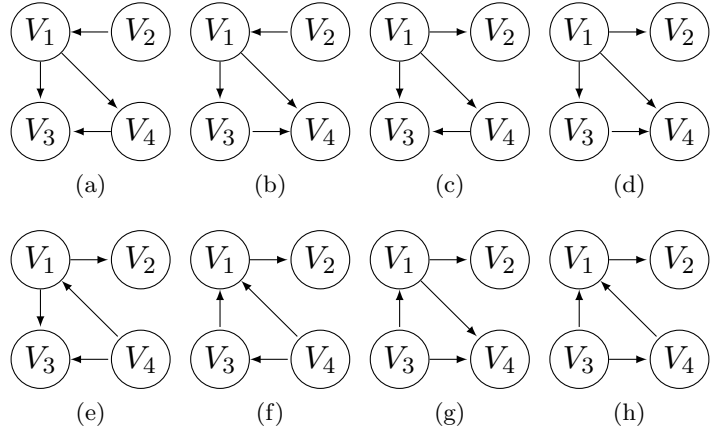 

\subsection{Existing results}

For completeness we first give the well known asymptotic behavior of the least squares estimator in the Gaussian setting. 

%If one restricts oneself to causal linear model with Gaussian errors, it alone would suffice to proof our results.

\begin{Lemma} \citep[Cf. Section 6.6.4 of][]{mardia1980multivariate}
	Let $(X,Y,\mathbf{Z}^{T})^{T}$ be a multivariate Gaussian vector with mean $\mathbf{0}$.
	Then 
	\[
	\sqrt{n}(\hat{\beta}_{yx.\mathbf{z}} - \beta_{yx.\mathbf{z}}) \xrightarrow{d} 
	\mathbf{N}\left( 0,\frac{\sigma_{yy.x\mathbf{z}}}
	{\sigma_{xx.\mathbf{z}}} \right) \quad \textit{as } n \rightarrow \infty,
	\] 
	where $\xrightarrow{d}$ denotes convergence in distribution.
	%	\begin{equation*}
	%	a.var(\hat{\beta}_{yx.\mathbf{z}}) = \frac{\sigma_{yy.x\mathbf{z}}}
	%	{\sigma_{xx.\mathbf{z}}}. 
	%	\end{equation*}
	%We suppress the dependence of $\hat{\beta}_{yx_i.\mathbf{x}_{-i}}$ on the sample size $n$ for the sake of simplicity. 
	\label{avar.equ}
\end{Lemma}

%In a linear causal model with non-Gaussian errors, a constituent random variables is only required to be linear in its parents and is not necessarily linear given another node set. One will thus often encounter misspecified least squares regressions. 

The next three lemmas give theoretical properties of the possibly misspecified least squares regression. They form the foundation for the extension of our results to causal linear models with non-Gaussian errors in Proposition \ref{avarDAG}.

\begin{Lemma} \citep[Cf. Section 3.1 of][]{buja2014models}
	Let $\mathbf{V} = (\mathbf{X}^T,\mathbf{Y}^T)^T$ be a mean $\mathbf{0}$ random vector with finite variance. Then the population level least squares regression coefficient matrix is $\beta_{\mathbf{yx}} = \Sigma_{\mathbf{yx}} \Sigma^{-1}_{\mathbf{xx}}$. 
%	where $\beta_{y_j\mathbf{x}}$ describes the best linear approximation of $Y_j$ given $\mathbf{X}$ for $j = 1, \dots k_y$.
	\label{lemma: buja}
\end{Lemma}

\begin{Lemma}
	Let $(\mathbf{X}^T,\mathbf{Y}^T,\mathbf{Z}^T)^T$ be a mean $\mathbf{0}$ random vector with finite variance such that $\mathbf{X} \ci \mathbf{Y} | \mathbf{Z}$. Then $\beta_{\mathbf{yx}.\mathbf{z}} = 0$.
	\label{lemma: independence}
\end{Lemma}

\begin{proof}
	This result is well known in the Gaussian setting \citep[cf. Section 2.5 of][]{anderson1958introduction}. It generalizes to random vectors with finite variance by the result from Lemma \ref{lemma: buja}, as $\beta_{\mathbf{yx}.\mathbf{z}} $ is fully determined by the covariance matrix alone.
\end{proof}

\begin{Lemma}{\cite[Cf. Corollary 11.1 of][]{buja2014models}}
	Let $(X,Y,\mathbf{Z}^T)^T$ be a mean $\mathbf{0}$ random vector with  finite variance and let $\mathbf{Z'} = (X,\mathbf{Z}^T)^T$. Then
	\[
	\sqrt{n}(\hat{\beta}_{yx.\mathbf{z}} - \beta_{yx.\mathbf{z}}) \xrightarrow{d} 
	\mathbf{N}\left(0,\frac{E[\delta^2_{x\mathbf{z}} \delta^2_{y\mathbf{z'}}]}{E[\delta^2_{x\mathbf{z}}]^2}\right) \quad \textit{as } n \rightarrow \infty,
	\] 
	where $\delta_{y\mathbf{z'}} = Y - \boldsymbol{\beta}_{y\mathbf{z'}}\mathbf{Z'}$ and $\delta_{x\mathbf{z}} = X - \boldsymbol{\beta}_{x\mathbf{z}}\mathbf{Z}$ denote the population error terms of the least squares regressions specified by their respective subscripts.
	%	where $\hat{\beta}_{yx_i.\mathbf{x}_{-i}}$ is the least squares estimator of the regression coefficient of $X_i$ in the regression of $Y$ on $\mathbf{X}$. 
	\label{buja}
\end{Lemma}

The following corollary gives a graphical criterion that is necessary and given a small restriction (see Remark \ref{remark: joint vs non}) also sufficient for the existence of a valid adjustment set in a DAG. 

%Since we typically assume the existence of a valid adjustment set it is often implicitly assumed to hold.

\begin{Corollary}(\citealp[Cf.\ Corollary 27 of][]{perkovic16})
	Let $\mathbf{X}$ and $\mathbf{Y}$ be disjoint node sets in a DAG $\g$. If there exists a valid adjustment set relative to $(\mathbf{X},\mathbf{Y})$ in $\g$, then $\mathbf{X} \cap \de(\cnb{\g},\g) = \emptyset$.  
	
	If we additionally assume that $\mathbf{Y} \subseteq \de(\mathbf{X},\g)$, then there exists a valid adjustment set relative to $(\mathbf{X},\mathbf{Y})$ in $\g$, if and only if $\mathbf{X} \cap \de(\cnb{\g},\g) = \emptyset$.
%	 if and only if there exists a valid adjustment set relative to $(\mathbf{X,Y})$ in $\g$.
	\label{cor:noforbx}
\end{Corollary}

\section{Proofs for Section \ref{prelim}}

\subsection{Proof of Proposition \ref{avarDAG}}
\label{sec:avar}

%
%\begin{Theorem}
%	Consider a DAG $\g=(\mathbf{V},\mathbf{E})$ such that $\mathbf{V}$ follows a causal linear model compatible with $\g$ and disjoint node sets $\mathbf{X}=\{X_1,\dots,X_{k_x}\}$ and $\mathbf{Y}=\{Y_1,\dots,Y_{k_y}\}$ in $\g$. Let $\mathbf{Z}$ be a valid adjustment set with respect to $(\mathbf{X},\mathbf{Y})$ in $\g$. Then 
%	\begin{equation*}
%	\sqrt{n} (\hat{\beta}_{y_jx_i.\mathbf{x}_{-i}\mathbf{z}} - (\tau^{\mathbf{z}}_{\mathbf{yx}})_{j,i})
%	\xrightarrow{D}
%	\mathcal{N}(0,\frac{\sigma_{y_jy_j.\mathbf{xz}}}{\sigma_{x_ix_i.\mathbf{x}_{-i}\mathbf{z}}}),
%	\end{equation*}
%	for all $i = 1,\dots,k_x $ and $j = 1, \dots,k_y$.
%	% where $\hat{\beta}_{y_jx_i.\mathbf{x}_{-i}\mathbf{z}}$ is the ordinary least squares estimator of $\beta_{y_jx_i.\mathbf{yx}_{-i}\mathbf{z}}$. 
%	Thus,
%	\[
%	a.var((\hat{\tau}^{\mathbf{z}}_{\mathbf{yx}}))_{i,j} = a.var(\hat{\beta}_{y_jx_i.\mathbf{yx}_{-i}\mathbf{w}}) = \frac{\sigma_{y_jy_j.\mathbf{xz}}}{\sigma_{x_ix_i.\mathbf{x}_{-i}\mathbf{z}}}.
%	\]
%	\label{avarDAG} 
%\end{Theorem}
%
%
%Theorem $\ref{avarDAG}$ implies that for total effect estimation via covariate adjustment in a causal linear model compatible with a causal graph, the error types in the generating linear equations are asymptotically insignificant. The consistency and asymptotic variance of the ordinary least squares estimator of the total effect does not depend on the error type and behaves as in the multivariate Gaussian case. This is generally \textit{not} the case when the conditioning set is not a valid adjustment set. 

	Recall that $\delta_{y\mathbf{z}} = Y - \boldsymbol{\beta}_{y\mathbf{z}}\mathbf{Z}$ denotes the population error term of the least squares regression specified by its subscripts. Moreover, we use the notation $\mathbf{Z}_{-x}$ to denote $\mathbf{Z} \setminus \{X\}$, for any set $\mathbf{Z}$.  
	
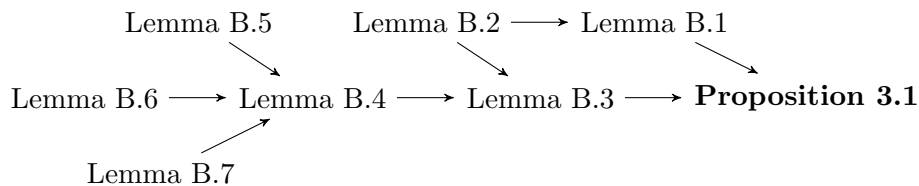
\begin{figure}[ht!]
	\centering
	\begin{tikzpicture}[>=stealth',shorten >=1pt,node distance=3cm, main node/.style={minimum size=0.4cm}]
	\node[main node,yshift=0cm,xshift=0cm]         (M) {\textbf{Proposition \ref{avarDAG}}};

	\node[main node,yshift=1cm,xshift=1cm]  (T1)  [left of= M]   {Lemma \ref{lemma:specialZ}}; 
	\node[main node,yshift=0cm]  (T2) [left of= T1]    {Lemma \ref{lemma:collider}}; 
	\node[main node,yshift=0cm]  (T3) [left of= T2]   {Lemma \ref{lemma:alwaysVAS}}; 
		\node[main node,yshift=0cm,xshift=-0.5cm] (M1) [left of= M] {Lemma \ref{ind-res}};
	\node[main node,yshift=0cm] (M2)  [left of= M1]   {Lemma \ref{open}};   
	\node[main node,yshift=0cm]  (M3) [left of= M2]   {Lemma \ref{lemma:edgeinto}}; 
	\node[main node,yshift=-1cm,xshift=1cm]  (L1) [left of= M2]   {Lemma \ref{lemma: TE-decomposition}}; 
	
%	
	
	%\node[main node,yshift=-1cm,xshift=-1cm]            
	%(C1)  at (T32)  {Corollary \ref{non-increasing}};  
	
	\draw[->] (M1) edge    (M);
    \draw[->] (M2) edge    (M1);
    \draw[->] (M3) edge    (M2);
 	\draw[->] (L1) edge    (M2);
	\draw[->] (T1) edge (M);
    \draw[->] (T2) edge (T1);
    \draw[->] (T2) edge (M1);
    \draw[->] (T3) edge (M2);
	
%	%\draw[->] (L2) edge    (C1);
%	\draw[->] (L2) edge    (T32);
%	\draw[->] (L3) edge    (T32);
	\end{tikzpicture}
	\caption{Proof structure of Proposition \ref{avarDAG}.}
	\label{fig: proof 2.13}
\end{figure}

\begin{proofof}[Proposition~\ref{avarDAG}]
	
	Consider a causal DAG $\g=(\mathbf{V},\mathbf{E})$, such that $\mathbf{V}$ follows a causal linear model compatible with $\g$. Let $\mathbf{X}=\{X_1,\dots,X_{k_x}\},\mathbf{Y}=\{Y_1,\dots,Y_{k_y}\}$ and $\mathbf{Z}$ be pairwise disjoint node sets in $\g$, such that $\mathbf{Z}$ is a valid adjustment set relative to $(\mathbf{X},\mathbf{Y})$ in $\g$.

%	Consider a causal DAG $\g=(\mathbf{V},\mathbf{E})$, such that $\mathbf{V}$ follows a causal linear model compatible with $\g$. Let $\mathbf{X}=\{X_1,\dots,X_{k_x}\}$ and $\mathbf{Y}=\{Y_1,\dots,Y_{k_y}\}$ be disjoint node sets in $\g$.
	
	In the proof of Proposition 3.1 from the Supplement of \citet{nandy2017estimating}, it is shown that $\tau_{yx}=\beta_{yx.\mathbf{z}}$ for singleton $X$ and $Y$, and $\mathbf{Z}=\pa(X,\g)$, independently of the error distributions in the causal linear model. Their argument relies on the fact that both terms of interest depend on the distribution of $\mathbf{V}$ only through $\Sigma_{\mathbf{vv}}$. The general case hence follows from the Gaussian one.
	
	Using Lemma \ref{lemma: buja}, this argument directly extends to a joint intervention of $\mathbf{X}$ on $\mathbf{Y}$ when $\mathbf{Z}$ is valid adjustment set relative to $(\mathbf{X},\mathbf{Y})$ in $\g$, implying  that  $\tau_{\mathbf{yx}} = \boldsymbol{\beta}_{\mathbf{yx}.\mathbf{z}}$.
	Finally, by Lemma \ref{buja}, $\hat{\beta}_{\mathbf{yx}.\mathbf{z}}$ is a consistent estimator of $\beta_{\mathbf{yx}.\mathbf{z}}$ for any conditioning set $\mathbf{Z}$.
	
	 We now show that our asymptotic variance statement holds. Fix $X_i \in \mathbf{X}$ and $Y_j \in \mathbf{Y}$ and let $\mathbf{Z'}= \mathbf{X} \cup \mathbf{Z}$. 
	We first show that the result holds if   $\delta_{y_j\mathbf{z'}} \ci \delta_{x_i\mathbf{z'}_{-x_i}}$. In this case, the statement of Lemma \ref{buja} simplifies to
	\[
	\sqrt{n}(\hat{\beta}_{y_jx_i.\mathbf{z'}_{-x_i}} - \beta_{y_jx_i.\mathbf{z'}_{-x_i}}) \xrightarrow{d} 
	\mathbf{N}\left(0,\frac{E[\delta^2_{y_j\mathbf{z'}}]}{E^[\delta^2_{x_i\mathbf{z'}_{-x_i}}]}\right) \quad \textit{as } n \rightarrow \infty.
	\] 
	Since $\mathbf{Z}$ is a valid adjustment set relative to $(\mathbf{X},\mathbf{Y})$ in $\g$, we have already shown that $(\tau_{\mathbf{yx}})_{j,i} =\beta_{y_jx_i.\mathbf{z'}_{-x_i}}$. Our claim then follows, as 
	\begin{align}
	E[\delta^2_{y_j\mathbf{z'}}] &= E[(Y_j - \boldsymbol{\beta}_{y_j\mathbf{z'}}\mathbf{Z'})^2] \nonumber \\
	&= E[Y_j^2] - 2 E[ \boldsymbol{\beta}_{y_j\mathbf{z'}}\mathbf{Z'} Y_j] + E[(\boldsymbol{\beta}_{y_j\mathbf{z'}}\mathbf{Z'})^2] \nonumber \\
	&= \sigma_{y_jy_j} - 2 \Sigma_{y_j\mathbf{z'}} \Sigma^{-1}_{\mathbf{z'z'}} \Sigma^T_{y_j\mathbf{z'}} + \Sigma_{y_j\mathbf{z'}} \Sigma^{-1}_{\mathbf{z'z'}}  \Sigma_{\mathbf{z'z'}} \Sigma^{-1}_{\mathbf{z'z'}} \Sigma^T_{y_j\mathbf{z'}} \nonumber \\
	&= \sigma_{y_jy_j} - \Sigma_{y_j\mathbf{z'}} \Sigma^{-1}_{\mathbf{z'z'}} \Sigma^T_{y_j\mathbf{z'}} \nonumber \\
	&= \sigma_{y_jy_j.\mathbf{z'}} \nonumber
	\end{align}
	and similarly
	\[
	E[\delta^2_{x_i\mathbf{z'}_{-x_i}}] = \sigma_{x_ix_i.\mathbf{z'}_{-x_i}}.
	\]
	
	It is left to show that $\delta_{y_j\mathbf{z'}} \ci \delta_{x_i\mathbf{z'}_{-x_i}}$. By Lemma \ref{lemma:specialZ}, all non-causal paths from $X_i$ to $Y_j$ are blocked by $ \mathbf{X}_{-i} \cup \mathbf{Z}$ and $\de(\mathbf{D},\g) \cap (\mathbf{X}_{-i} \cup \mathbf{Z}) =\emptyset$, where $\mathbf{D}=\text{cn}(X_i,Y_j,\g) \cap \text{cn}(\mathbf{Z'},Y_j,\g)$. 
	Hence, we can apply Lemma \ref{ind-res} with  $X=X_i,Y=Y_j$ and $\mathbf{Z'}= \mathbf{X} \cup \mathbf{Z}$ to conclude that in fact $\delta_{y_j\mathbf{z'}} \ci \delta_{x_i\mathbf{z'}_{-x_i}}$.

%	The proof of this result is quite complex as it uses typically overlooked properties of both misspecified ordinary least squares regression and causal linear models, so we give a summary of its structure. 
%	By see Lemma \ref{buja}, it suffices to show that $\delta_{y_j\mathbf{z'}} \ci \delta_{x_i\mathbf{z'}_{-x_i}}$, for all $X_i \in \mathbf{X}$ and all $Y_j \in \mathbf{Y}$, where $\mathbf{Z'} = \mathbf{X} \cup \mathbf{Z}$ for our claim to follow. 
%These residuals are each functions in a unique minimally sized subset of the error set from the underlying causal linear model. The key and most innovative part of this proof is Lemma \ref{open}, in which we manage to tie the appearance of an error term in such a minimal subset to a graphical condition on its corresponding node. Lemmas \ref{block} and \ref{ind-res} are then simple technical results that use Lemma \ref{open} to show that when a valid adjustment set is considered, an error term may not appear in both the minimal subsets of error terms corresponding to the two residuals of interest.
%	
%	
\end{proofof}

Lemma \ref{lemma:specialZ} relies on the technical Lemma \ref{lemma:collider} that will be used throughout this Supplement. For any given pair $X \in \mathbf{X}$ and $Y \in \mathbf{Y}$, it gives some properties of the set $\mathbf{Z'}_{-x}=\mathbf{X}_{-x} \cup \mathbf{Z}$. Albeit $\mathbf{Z'}_{-x}$ is not necessarily a valid adjustment set relative to $(X,Y)$, it behaves similarly to one. Lemma \ref{ind-res} then relies on these properties and the technical Lemma \ref{open} to show that our two residuals of interest are in fact independent. A full summary of how the following lemmas relate to each other is given in Figure \ref{fig: proof 2.13}.

\begin{Lemma}
	Let $\mathbf{X},\mathbf{Y}$ and $\mathbf{Z}$ be pairwise disjoint node sets in a causal DAG $\g$, such that $\mathbf{Z}$ is a valid adjustment set with respect to $(\mathbf{X},\mathbf{Y})$ and let $\mathbf{Z'}=\mathbf{X} \cup \mathbf{Z}$. Consider any pair of nodes $X\in \mathbf{X}, Y \in \mathbf{Y}$ and let $\mathbf{D}=\textup{cn}(X,Y,\g) \cap \textup{cn}(\mathbf{Z'},Y,\g)$ be the set of all nodes $N \in \cn{\g}$, such that there exists a directed path from $N$ to $Y$ which contains no nodes from $\mathbf{Z'}$.
	Then the following two statements hold:
	\begin{enumerate}
		\item All non-causal paths from $X$ to $\mathbf{Y}$ are blocked by $\mathbf{Z'}_{-x}$. \label{block}
		\item $\de(\mathbf{D},\g) \cap \mathbf{Z'}_{-x} =\emptyset$. \label{specialforb}
	\end{enumerate}

	\label{lemma:specialZ}
\end{Lemma}

\begin{proof}
	We first prove Statement \ref{block}. Fix $X \in \mathbf{X}$ and $Y \in \mathbf{Y}$. We will prove our claim by contradiction, so suppose that there exists a non-causal path $p$ from $X$ to $\mathbf{Y}$ that is open given $\mathbf{Z'}_{-x}$ and assume that $\mathbf{Z}$ is a valid adjustment set relative to $(\mathbf{X,Y})$ in $\g$. We will show that this implies the existence of a proper non-causal path from $\mathbf{X}$ to $\mathbf{Y}$ that is open given $\mathbf{Z}$, contradicting our assumption that $\mathbf{Z}$ is a valid adjustment set.
	
	Suppose first that no non-endpoint node of $p$ is in $\mathbf{X}$, so that $p$ is a proper non-causal path from $\mathbf{X}$ to $\mathbf{Y}$. By the assumption that $\mathbf{Z}$ is a valid adjustment set, $p$ must then be blocked by $\mathbf{Z}$. Since $p$ is assumed to be open given $\mathbf{Z'}_{-x}$ it is clearly also open given $\mathbf{Z'}=\mathbf{X} \cup \mathbf{Z}$ and hence, Lemma \ref{lemma:collider} with $A=X$ implies that there exists a proper non-causal path from $\mathbf{X}$ to $\mathbf{Y}$ that is open given $\mathbf{Z}$.
	
	Next, suppose that $p$ contains some non-endpoint node in $\mathbf{X}$. Let $X'$ be the node in $\mathbf{X}$ on $p$ that is closest to $Y$. Then $p'=p(X',Y)$ is a proper subpath of $p$. Since $p$ is open given $\mathbf{Z'}_{-x}$, $X'$ must be a collider on $p$. But then $p'$ is both proper and non-causal, and we can repeat the argument from the previous paragraph.

	We now prove Statement \ref{specialforb}. We will first show that $\mathbf{D} \subseteq \cnb{\g}$. Consider a node $D \in \mathbf{D}$. Then $D$ must lie on at least one causal path $p$ from $X$ to $Y$, where we can choose $p$ such that $p(D,Y)$ contains no node in $\mathbf{Z'}$. Let $X' \in \mathbf{X}$ be the node closest to $D$ on $p(X,D)$. Then $p(X',Y)$ is a proper causal path from $\mathbf{X}$ to $\mathbf{Y}$ containing $D$ and hence $D \in \cnb{\g}$.
	
	We will now prove the statement by contradiction, so assume that there exists a node $F \in \de(\mathbf{D},\g) \cap \mathbf{Z'}_{-x}$ and that $\mathbf{Z}$ is a valid adjustment set. Assume first that $F \in \mathbf{X}_{-x}$. But this implies that $\de(\cnb{\g}) \cap \mathbf{X} \neq \emptyset$, which contradicts our assumption that $\mathbf{Z}$ is a valid adjustment set by Corollary \ref{cor:noforbx}. Now assume that $F \in \mathbf{Z}$. In this case $F \in \fb{\g}$, again contradicting our assumption that $\mathbf{Z}$ is a valid adjustment set.
\end{proof}

\begin{Lemma}
	Let $\mathbf{X},\mathbf{Y}$ and $\mathbf{Z}$ be pairwise disjoint node sets in a causal DAG $\g$. Let $A \notin \mathbf{Y}$ be a node and consider a path $p$ from $A$ to $\mathbf{Y}$ in $\g$. If  $p$ is blocked by $\mathbf{Z}$ and open given $\mathbf{X} \cup \mathbf{Z}$, then there exists a proper non-causal path from $\mathbf{X}$ to $\mathbf{Y}$ that is open given $\mathbf{Z}$. 
	\label{lemma:collider}
\end{Lemma}

\begin{proof}
	Let $Y \in \mathbf{Y}$ be the endpoint of $p$.
	By assumption $p$ is blocked by $\mathbf{Z}$, while being open given $\mathbf{X} \cup \mathbf{Z}$. This requires the following three statements to hold:
	\begin{enumerate}
		\item For any non-collider $N$ on $p$, $N \notin \mathbf{X} \cup \mathbf{Z}$.
		\item For any collider $C$ on $p$, $\de(C,\g) \cap (\mathbf{X} \cup \mathbf{Z}) \neq \emptyset$.
		\label{not-special}
		\item There exists at least one collider $C'$ on $p$ such that $\de(C',\g) \cap \mathbf{Z} = \emptyset$.
		\label{special}
		%		\item There exists at least one collider $C'$ on $p$ such that $\de(C',\g) \cap \mathbf{X} \neq \emptyset$. \label{special}
	\end{enumerate}
	Let $C'$ be the node closest to $Y$ fulfilling Statement \eqref{special}. By choice of $C'$ and the assumption that $p$ is open given $\mathbf{X} \cup \mathbf{Z}$, $p'=p(C',Y)$ is open given $\mathbf{Z}$ and $\mathbf{X} \cup \mathbf{Z}$. If $p'$ contains a node in $\mathbf{X}$, any such node must be a collider. Assume that there is such a node and let $X' \in \mathbf{X}$ be the one closest to $Y$ on $p'$. Since $X'$ is a collider, $p'(X',Y)$ is non-causal and by choice of $X'$, $p'(X',Y)$ is also a proper path from $\mathbf{X}$ to $\mathbf{Y}$. Since $p'$ is open given $\mathbf{Z}$ so is $p'(X',Y)$ and hence we are done. Therefore, we will from now on suppose that no node on $p'$ is in $\mathbf{X}$. 
	
	As $C'$ fulfills the requirements of both Statement \eqref{not-special} and \eqref{special}, there exists a path $p''$ of the form 
	$
	C' \rightarrow \dots \rightarrow X'' \in \mathbf{X},
	$
	that is open given $\mathbf{Z}$ and where we choose $X''$ so that $p''$ contains no other node in $\mathbf{X}$.
	Let $I$ be the node closest to $X''$ on $p''$ that is also on $p'$ and consider $q = p''(X'',I) \oplus p'(I,Y)$.  The path $q$ is a path from $\mathbf{X}$ to $\mathbf{Y}$ and we now show that it is proper and non-causal. Since $I$ lies on $p'$ and $p'$ contains no node on $\mathbf{X}$ by assumption, $I \neq X''$. Thus, $q$ is a proper and by nature of $p''$ being directed towards $X''$, non-causal path from $\mathbf{X}$ to $\mathbf{Y}$.
	
	It now only remains to show that $q$ is open given $\mathbf{Z}$. As the two constituent paths are both open given $\mathbf{Z}$ it suffices to consider $I$. By choice of $C'$, no node in $p''$ may be in $\mathbf{Z}$ and hence neither is $I$. Furthermore, by nature of $p''$ being directed towards $X''$, $I$ cannot be a collider on $q$ and it thus follows that $q$ is open given $\mathbf{Z}$.  
\end{proof}

\begin{Lemma}
	 Let $X$ and $Y$ be two nodes in a causal DAG $\g=(\mathbf{V},\mathbf{E})$, let $\mathbf{V}$ follow a causal linear model compatible with $\g$ and let $\mathbf{Z'}$ be a node set, such that $X \in \mathbf{Z'}$ and $Y\notin \mathbf{Z'}$. If $\mathbf{Z'}_{-x}$ is a set that
	 \begin{enumerate}
	 	\item blocks all non-causal paths from $X$ to $Y$ in $\g$ and 
	 	\item fulfills that $\de(\mathbf{D},\g) \cap \mathbf{Z'}_{-x} =\emptyset$, where $\mathbf{D}=\text{cn}(X,Y,\g) \cap \text{cn}(\mathbf{Z'},Y,\g)$,
	 \end{enumerate}
	 then $\delta_{y\mathbf{z'}} \ci \delta_{x\mathbf{z'}_{-x}}$.
	\label{ind-res}
\end{Lemma}

\begin{proof}
	In order to simplify notation we refer to $\mathbf{Z'}_{-x}$ as $\mathbf{Z}$ throughout this proof. Let $\boldsymbol{\epsilon}= \{\epsilon_{v_1},\dots, \epsilon_{v_p}\}$ be the set of errors from the underlying causal linear model and consider 
	$\delta_{x\mathbf{z}}$ and $\delta_{y\mathbf{z'}}$ as functions of $\boldsymbol{\epsilon}$. Then there are minimally sized subsets $\boldsymbol{\epsilon}^{x\mathbf{z}}$ and $\boldsymbol{\epsilon}^{y\mathbf{z'}}$ of $\boldsymbol{\epsilon}$, such that $\delta_{x\mathbf{z}}$ and $\delta_{y\mathbf{z'}}$ are functions of $\boldsymbol{\epsilon}^{x\mathbf{z}}$ and $\boldsymbol{\epsilon}^{y\mathbf{z'}}$ respectively. It suffices to show that $\boldsymbol{\epsilon}^{x\mathbf{z}} \cap \boldsymbol{\epsilon}^{y\mathbf{z'}} = \emptyset$, as then $\delta_{y\mathbf{z'}} \ci \delta_{x\mathbf{z}}$ follows from the joint independence of the errors in $\boldsymbol{\epsilon}$.
	
	We will prove our claim by contraposition, so assume that there exists a node $N$ such that $\epsilon_n \in \boldsymbol{\epsilon}^{x\mathbf{z}} \cap \boldsymbol{\epsilon}^{y\mathbf{z'}}$. We now show that the existence of such a node $N$ implies either the existence of a non-causal path from $X$ to $Y$ that is open given $\mathbf{Z}$ or that $\de(\mathbf{D},\g) \cap \mathbf{Z} \neq\emptyset$. We will do so by going through a series of cases:
	\begin{mylist}
		\item $N =X \in \mathbf{Z'}$, \label{equalX}
		\item $N \in \mathbf{Z}$ \text{ and}\label{inZ} 
		\item $N \notin \mathbf{Z'}$. \label{neither}
	\end{mylist}

	Case \ref{equalX}:  By Lemma \ref{open} with $A=Y$ and $\mathbf{W}=\mathbf{Z'}$ there exists a non-causal path $p$ from $X$ to $Y$ that is open given $\mathbf{Z'}$. Suppose, that $p$ is blocked by $\mathbf{Z}$. Then by Lemma \ref{lemma:collider} with $A=X$,$\mathbf{X}=X,\mathbf{Y}=Y$ and $\mathbf{Z}=\mathbf{Z}$ there exists a non-causal path from $X$ to $Y$ that is open given $\mathbf{Z}$. Otherwise, $p$ is itself such a path.

	Case \ref{inZ}: By applying Lemma \ref{open} twice, once with $A=X$ and $\mathbf{W}=\mathbf{Z}$ and once with $A=Y$ and $\mathbf{W}=\mathbf{Z'}$ we deduce that  there exists i) a path $p$ of the form $X \dots \rightarrow N$ that is open given $\mathbf{Z}$, and ii) a path $p'$ of the form $N \leftarrow \cdots Y$ that is open given $\mathbf{Z'}$. 
	If $p'$ blocked by $\mathbf{Z}$, we can apply Lemma \ref{lemma:collider} with $A=N, \mathbf{X}=X, \mathbf{Y}=Y$ and $\mathbf{Z}=\mathbf{Z}$ to conclude that there exists a non-causal path from $X$ to $Y$ that is open given $\mathbf{Z}$. For the remainder of Case \ref{inZ} we will suppose that $p'$ is open given $\mathbf{Z}$.
	
	Let $I$ be the node closest to $X$ on $p$ that is also on $p'$. We now show that either $q=p(X,I) \oplus p'(I,Y)$ is a non-causal path from $X$ to $Y$ that is open given $\mathbf{Z}$ or $\de(\mathbf{D},\g) \cap \mathbf{Z} \neq \emptyset$.
	Since both $p(X,I)$ and $p'(I,Y)$ are open given $\mathbf{Z}$, it suffices to consider $I$ to decide whether $q$ is open given $\mathbf{Z}$, so we will now sequentially consider the cases that $I $ is $N,X,Y$, a non-collider on $q$ and collider on $q$.
	
	First suppose that $I=N$. Then $q$ is of the form $X \dots \rightarrow N \leftarrow \cdots Y$ and $I=N \in \mathbf{Z}$ is a collider on $q$. Hence, $q$ is open given $\mathbf{Z}$ as well as non-causal.

	Now, consider the case that $I\in \{X,Y\}$. Then $q$ is a subpath of $p'$ or $p$ respectively and hence trivially open given $\mathbf{Z}$.  If $q$ is non-causal, the first possible claim of our contrapositive statement holds true. 
	Hence, suppose that $q$ is a causal path from $X$ to $Y$. If $I=X$ and $q$ is causal, $p'$ is of the form $N \dots X \rightarrow \dots \rightarrow Y$ and hence $X$ is a non-collider on $p'$. But this contradicts that $p'$ is open given $\mathbf{Z'}$. If $I=Y$, then $q$ is a subpath of $p$ and obviously $p(X,Y)=q$ is causal. Since, $p$ is open given $\mathbf{Z}$ and $Y\notin \mathbf{Z}$ it follows that either $Y$ is collider on $p$, with $\de(Y,\g) \cap \mathbf{Z} \neq \emptyset$ or that $Y$ is a non-collider on $p$ and hence $p(Y,N)$ is of the from $Y \rightarrow \dots N$. In the latter case, since $N \in \de(\mathbf{Z},\g)$ and $p(Y,N)$ is open given $\mathbf{Z}$, $p(Y,N)$ either contains a collider $C$, such that $C \in \de(Y,\g)$ and $\de(C,\g) \cap \mathbf{Z} \neq \emptyset$ or $N \in \de(Y,\g)$. Thus, in either case $\de(Y,\g) \cap \mathbf{Z} \neq \emptyset$.
	Further, since $q$ is a causal path from $X$ to $Y$ that is open given $\mathbf{Z}$, $q$ cannot contain nodes from $\mathbf{Z}$ and therefore $Y \in \mathbf{D}$. Thus, $\de(D,\g) \cap \mathbf{Z} \neq \emptyset$. 
	
%	If $I \in \{X,Y\}$, $q$ is clearly open given $\mathbf{Z}$.
	 We now suppose that $I \notin \{N,X,Y\}$ is a non-collider on $q$. Then $I$ cannot be a collider on both $p$ and $p'$. Since $p$ and $p'$ are open given $\mathbf{Z}$ it thus follows that $I \notin \mathbf{Z}$. Therefore, $q$ is also open given $\mathbf{Z}$. 
	 If $q$ is non-causal, the first possible claim of our contrapositive statement holds true. Hence, suppose that $q$ is a causal path from $X$ to $Y$.  
	 As $q$ is causal, $ p(X,I)$ is directed towards $I$ and we have already shown that $I\notin \mathbf{Z}$.
	 By the same argument as in the case $I=Y$, it then follows that $I \in \mathbf{D}$ and $\de(I,\g) \cap \mathbf{Z} \neq \emptyset$. Thus, $\de(\mathbf{D},\g) \cap \mathbf{Z} \neq \emptyset$.

%	 Since $p$ is open given $\mathbf{Z}$ and $I\notin \mathbf{Z}$ it follows that $I$ is either a collider on $p$, with $\de(I,\g)$ or that $I$ is a non-collider on $p$ and hence $p(I,N)$ is of the form $I \rightarrow \dots N$. In the latter case, since $p$ is open given $\mathbf{Z}$ it follows that $p(I,N)$ either contains a collider $C$, such that $C \in \de(I,\g)$ and $\de(C,\g)\cap \mathbf{Z}\neq \emptyset$ or $N \in \de(I,\g)$. In either case, $\de(I,\g) \cap \mathbf{Z} \neq \emptyset$. Further,  since $I$ lies on $q$, a causal path from $X$ to $Y$ that is open given $\mathbf{Z}$ and hence contains no nodes from $\mathbf{Z}$, $I \in \mathbf{D}$.
%	 Therefore, $\de(\mathbf{D},\g) \cap \mathbf{Z} \neq \emptyset$.

	Consider now the case that $I \notin \{N,X,Y\}$ is a collider on $q$. 
	Clearly, $q$ is non-causal. Further, if $I$ is also a collider on either $p$ or $p'$, it follows that $\de(I, \g) \cap \mathbf{Z} \neq \emptyset$ and hence, $q$ is open given $\mathbf{Z}$ in this case. Suppose that $I$ is a collider on $q$, while being a non-collider on both $p$ and $p'$. 
	Then $p$ must be of the form $X \dots 
	\rightarrow I \rightarrow \dots N$ and since $p(I,N)$ is open given $\mathbf{Z}$ and $N \in \de(\mathbf{Z},\g)$, it follows that $\de(I, \g) \cap \mathbf{Z} \neq \emptyset$. Hence, $q$ is non-causal and open given $\mathbf{Z}$. 
	
	Case \ref{neither}: 
	Let us first suppose that $\mathbf{Z} \cap \de(N,\g) = \emptyset$. By Lemma \ref{open} with $A=Y$ and $\mathbf{W} = \mathbf{Z'}$, it follows that $\epsilon_x \in \boldsymbol{\epsilon}^{y\mathbf{z'}}$ or $\epsilon_y \in \boldsymbol{\epsilon}^{y\mathbf{z'}}$.  
	If $\epsilon_x \in \boldsymbol{\epsilon}^{y\mathbf{z'}}$, then we are done by Case \ref{equalX}, so suppose that only the latter statement is true. By Lemma \ref{open} there then exists a directed path $p'$ from $N$ to $Y$ that is open given $\{X,Y\} \cup \mathbf{Z}$ and as a directed path is therefore also open given $\mathbf{Z}$. 
	By Lemma \ref{open} with $A=X$ and $\mathbf{W} = \mathbf{Z}$, it follows that there exists a path $p$ from $X$ to $N$ that is directed towards $X$ and open given $\{X\} \cup \mathbf{Z}$ and hence, is also open given $\mathbf{Z}$. 
	
	Let $I$ be the node closest to $X$ on $p$ that is also on $p'$. We now show that $q=p(X,I) \oplus p'(I,Y)$ is a non-causal path from $X$ to $Y$ that is open given $\mathbf{Z}$.
%	 or $\de(\mathbf{D},\g) \cap \mathbf{Z} \neq \emptyset$.
	Since both $p(X,I)$ and $p'(I,Y)$ are open given $\mathbf{Z}$, it suffices to consider $I$ to decide whether $q$ is open given $\mathbf{Z}$, so we will now sequentially consider the possible properties of $I$. 
	
	If $I=X$, then $q$ is subpath of $p'$ and $X$ lies on $p'$. But as $p'$ is directed this contradicts that it is open given $\{X\} \cup \mathbf{Z}$. If $I=Y$, then $q$ is a subpath of $p$ and hence a non-causal path from $X$ to $Y$ that is open given $\mathbf{Z}$. If $I \notin \{Y,X\}$, then $I$ is a non-collider on $q$ and since no node in $p$ may be in $\mathbf{Z}$, $I \notin \mathbf{Z}$ and it thus follows that $q$ is open given $\mathbf{Z}$. As $p(X,I)$ is directed towards $X$, $q$ is non-causal. Thus, $q$ is a non-causal path from $X$ to $Y$ that is open given $\mathbf{Z}$.

	For the remainder of Case \ref{neither}, we suppose that $\mathbf{Z} \cap \de(N,\g) \neq \emptyset$. By applying Lemma \ref{open} twice, once with $A=X$ and $\mathbf{W}=\mathbf{Z}$ and once with $A=Y$ and $\mathbf{W}=\mathbf{Z'}$, we deduce that there exists i) a path $p$ of from $X$ to $N$ that is open given $\mathbf{Z}$ and ii) a path $p'$ from $N$ to $ Y$ that is open given $\mathbf{Z'}$. 
	If $p'$ is blocked by $\mathbf{Z}$, we can conclude with Lemma \ref{lemma:collider}, as in Case \ref{inZ}, that there exists a non-causal path from $X$ to $Y$ that is open given $\mathbf{Z}$ and are done. For the remainder of Case \ref{neither}, we suppose that $p'$ is open given $\mathbf{Z}$.
	
%	Let us first suppose that $\mathbf{Z} \cap \de(N,\g) = \emptyset$. This implies that there is a directed path $p''$ from $N$ to $X$ that is open given $\mathbf{Z}$. Let $I$ be the node closest to $X$ on $p''$ and let $q=p''(X,I) \oplus p'(I,Y)$. If $I \neq X$  then by nature of $p''(X,I)$, $q$ is a non-causal path and $I$ may not be a collider on $q$ and since $I \notin \mathbf{Z}$ it follows that $q$ is open given $\mathbf{Z}$. If $I=X$, then by 	

Let $I$ be the node closest to $X$ on $p$ that is also on $p'$. As in case \ref{inZ}, we will now show that either $q=p(X,I) \oplus p'(I,Y)$ is a non-causal path from $X$ to $Y$ that is open given $\mathbf{Z}$ or $\de(\mathbf{D},\g) \cap \mathbf{Z} =\emptyset$. We now sequentially consider the possible properties of $I$.
	
	First suppose that $I=N$. Since $\mathbf{Z} \cap \de(N,\g) \neq \emptyset$ and $N \notin \mathbf{Z}$ we can immediately conclude that $q$ is open given $\mathbf{Z}$ independently of whether $N$ is a collider or a non-collider on $q$. Further, if $q$ is causal it follows that $N \in \mathbf{D}$ and hence $\de(\mathbf{D},\g) \cap \mathbf{Z} = \emptyset$. 
	
	Suppose now $I \neq N$. Recall that $\mathbf{Z} \cap \de(N,\g) \neq \emptyset$ and that we have already shown that there exists a path $p$ from $X$ to $N$ and another $p'$ from $N$ to $Y$, such that the former is open given $\mathbf{Z}$ and the latter open given both $\mathbf{Z}$ and $\mathbf{Z'}$. But these are exactly the assumptions required to show our claim in the corresponding case $I \neq N$ in Case \ref{inZ}. Hence, our claim follows with the same argument.

\end{proof}

%\begin{Lemma}
%	Let $X$ and $Y$ be nodes in a causal DAG $\g$. Let $\mathbf{Z}$ be a node set in $\g$, such that $X,Y \notin \mathbf{Z}$. Consider a path $p$ from $X$ to $Y$ that is open given $\{X\} \cup \mathbf{Z}$. If $p$ is blocked by $\mathbf{Z}$, then there exists a non-causal path from $X$ to $Y$, that is open given $\mathbf{W}$.
%	\label{lemma:weird-endpoint}
%\end{Lemma}
%
%\begin{proof}
%	Suppose that $p$ is blocked by $\mathbf{Z}$ and open given $\{X\} \cup \mathbf{Z}$. Then there must exist a collider $C$ on $p$, such that $X \in \de(C,\g)$ and $\mathbf{Z} \cap \de(C,\g)=\emptyset$, where we can choose $C$ ensuring that $p(C,Y)$ is open given $\mathbf{Z}$. Therefore, there exists a path $p'$ from $X$ to $C$ that is open given $\mathbf{Z}$ and directed towards $X$. Let $I$ be the node closest to $Y$ on $p(C,Y)$ that is also in $p'$  and consider $q=p'(X,I) \oplus p(I,Y)$. Since $I$ lies on $p(C,Y), I \neq X$. It follows by nature of $p'$ that $I \notin \mathbf{Z}$,  that $q$ is a non-causal path from $X$ to $Y$ and that either $I=Y$ or that $I$ is a non-collider on $q$. 
%	If $I = Y$,  then $q=p'(X,I)$ is a subpath of $p'$ that is open given $\mathbf{Z}$. If $I \neq Y$, then by the nature of $I$  and the fact $p'(X,I)$ and $p(I,Y)$ are open given $\mathbf{Z}$ by construction it follows that $q$ is open given $\mathbf{Z}$. Hence, in all case $q$ is a non-causal path from $X $ to $Y$ that is open given $\mathbf{Z}$.
%\end{proof}

\begin{Lemma}
	Consider a causal DAG $\g=(\mathbf{V},\mathbf{E})$, let $\mathbf{V}$ follow a causal linear model compatible with $\g$ and let $\boldsymbol{\epsilon}= \{\epsilon_{v_1},\dots, \epsilon_{v_p}\}$ be the set of errors from the underlying causal linear model. Given a node $A$ and a node set $\mathbf{W}=\{W_1,\dots,W_k\}$, such that $A \notin \mathbf{W}$, the residual $\delta_{a\mathbf{w}}$ is a function of some minimally sized subset $\boldsymbol{\epsilon}^{a\mathbf{w}}$ of $\boldsymbol{\epsilon}$. 
	\begin{enumerate}
		\item Let $M \in \mathbf{W}$. If $\epsilon_m \in \boldsymbol{\epsilon}^{a\mathbf{w}}$, then there exists a path of the form 
		$
		A \dots \rightarrow M,
		$
		that is open given $\mathbf{W}$.  \label{in}
		\item Let $M \notin \mathbf{W}$. If $\epsilon_m \in \boldsymbol{\epsilon}^{a\mathbf{w}}$, then there exists  a directed path from $M$ to some node $M' \in (\{A\} \cup \mathbf{W})$ that is open given $\{A\} \cup \mathbf{W}$ and $\epsilon_{m'} \in \epsilon^{a\mathbf{w}}$. Consequentially, there exists a path from $A$ to $M$ that is open given $\mathbf{W}$.
		 \label{out}		
%		
%		of the form 
%		$
%		A \dots \leftarrow C
%		$ 
%		that is open given $\mathbf{W}$.
	\end{enumerate}
	\label{open}
\end{Lemma}

\begin{proof}
	We first give some preparatory thoughts on how to write $\delta_{a\mathbf{w}}$ as a function in the errors from the causal linear model. Consider a random vector $\mathbf{V}$ that follows a causal linear model compatible with a causal DAG $\g=(\mathbf{V},\mathbf{E})$ and fix some node $V_i \in \mathbf{V}$. 
%	
%	Then
%	\begin{align*}
%	V_i = \sum_{V_j \in \pa(V_i,\g)} \alpha_{ij}V_j +\epsilon_{v_i}, 
%	\end{align*}
%	with similar equations holding for the nodes $V_j \in \pa(V_i,\g)$. By sequentially plugging in these equations 
%	
	We can then write $V_i$ as a linear function of the generating errors in the following way:
	\begin{align}
	V_i = \sum_{V_j \in \mathbf{V}} \tau_{v_iv_j} \epsilon_{v_j} = \sum_{V_j \in \an(V_i,\g)} \tau_{v_iv_j} \epsilon_{v_j}, \label{error.poly} 
	\end{align}
	where we use the convention that $\tau_{v_iv_i}=1$ for any $V_i \in \mathbf{V}$ and make use of the fact that $\tau_{v_jv_i}=0$, whenever $V_i \notin \an(V_j,\g)$ \citep{SEM1989}. 
	
	%or alternatively $V_j \notin (de(V_i,\g) \setminus \{V_i\})$. 
	
	%This can be easily verified by starting with the causal linear model equation for $V_j$ and then sequentially plugging in that same equation for the 

	Consider the equation
	\begin{align}
	\delta_{a\mathbf{w}} = A - \sum_{W_i \in \mathbf{W}} \beta_{aw_i.\mathbf{w}_{-i}} W_i. \label{deltas}
	\end{align}
	%	Given a causal linear model of a random vector $\mathbf{V}$, any of the constituent random variables $V_i \in \mathbf{V}$, can be written as a polynomial in the generating errors in the following way:
	%	\begin{align}
	%	V_j = \sum_{V_i \in \mathbf{V}} \tau_{v_jv_i} \epsilon_{V_i}, \label{error.poly}
	%	\end{align}
	%	where we use the convention that $\tau_{v_iv_i}=1$ for any $V_i \in \mathbf{V}$. 
	Applying equation \eqref{error.poly} to the $A$ and the $W_i$ terms in equation \eqref{deltas}, we can write 
	$\delta_{a\mathbf{w}}$ 
	as a linear function in the error terms of the generating causal linear model of the form
	\begin{align}
		\delta_{a\mathbf{w}} = \sum_{V_j \in \mathbf{V}} \gamma_{v_j} \epsilon_{v_j},
		\label{equ-errors}
	\end{align}
	with coefficients $\gamma_{v_j} \in \mathbb{R}$.
	
	We now prove Statement \eqref{in} by contraposition. So assume that $M \in \mathbf{W}$ and that there exists no path $p$ from $A$ to $M$ whose last edge points into $M$ and which is open given $\mathbf{W}$. We will now show that this implies that $\epsilon_m \notin \boldsymbol{\epsilon}^{a\mathbf{w}}$
	
	 It is sufficient to show that the coefficient $\gamma_m$ corresponding to $\epsilon_m$ in equation \eqref{equ-errors} is equal to 0. The value of $\gamma_m$ is
	\begin{align}
	\gamma_m = \tau_{am} - \sum_{W_i \in \mathbf{W}} \beta_{aw_i.\mathbf{w}_{-i}} \tau_{w_im}. 
	\label{eq:gamma}
	\end{align}
	Our claim thus follows, if we show that  
	\begin{align*}
	\tau_{am} = \sum_{W_i \in \de(M,\g) \cap \mathbf{W}} \beta_{aw_i.\mathbf{w}_{-i}} \tau_{w_im}, 
	\end{align*}
	where we have simplified the sum by removing those $W_i \in \mathbf{W}$ with $\tau_{w_im}=0$.
	
	%We now show that equation \eqref{cool} does indeed hold. 
	Let $\mathbf{W}'=\de(M,\g) \cap \mathbf{W},\mathbf{W}''=\mathbf{W} \setminus \de(M,\g)$ and $\mathbf{W}'''=\pa(M,\g) \cup \mathbf{W}''$. 
	By construction, $\mathbf{W}'''$ contains all parents of $M$ while containing no  descendants of $M$. It thus follows with Lemma \ref{lemma:alwaysVAS} that $\mathbf{W}'''$ is a valid adjustment set relative to $M$ and any node that is not in $\mathbf{W}'''$. We note that $\mathbf{W}' \cap \mathbf{W}'''=\emptyset$ by construction. Further, $A \notin \pa(M,\g)$ by assumption and hence $A \notin \mathbf{W}'''$. 
	
	Using the already proven first half of Proposition \ref{avarDAG} to replace the total effects with appropriate regression coefficients and vice versa, it follows that
	\begin{align}
	\tau_{am} &= \beta_{am.\mathbf{w}'''}  \nonumber \\
	&= \beta_{am.\mathbf{w}_{-m}'\mathbf{w}'''} + \sum_{W_i \in \mathbf{W}_{-M}'} \beta_{aw_i.\mathbf{w}_{-w_i}'\mathbf{w}'''} \beta_{w_im.\mathbf{w}'''} \nonumber \\
	&= \sum_{W_i \in \mathbf{W'}} \beta_{aw_i.\mathbf{w}_{-w_i}} \tau_{w_im}.  \nonumber 
	\end{align}
	Here, we use firstly that $\mathbf{W}'''$ is a valid adjustment set with respect to $M$ and any node not in $\mathbf{W}'''$ to conclude that $\beta_{w_im.\mathbf{w'''}}=\tau_{w_im}$ and $\beta_{am.\mathbf{w'''}}=\tau_{am}$.
	Secondly, we use Lemma \ref{Cochran} in the second step, with $\mathbf{T}=A,\mathbf{W}=M,\mathbf{Z}=\mathbf{W}'''$ and $\mathbf{S}=\mathbf{W}'_{-m}$. Lastly, $\mathbf{P} = \pa(M,\g) \setminus \mathbf{W} \perp_{\g} A | \mathbf{W}$ by Lemma \ref{lemma:edgeinto} and we use this result to simplify the conditioning sets in step three by invoking the first statement from Lemma \ref{Wermuth}, with $\mathbf{T}=\mathbf{P}, \mathbf{S}=\emptyset, \mathbf{X}=\mathbf{W}$ and $\mathbf{Y}=\{A\}$, allowing us to drop all nodes in $\mathbf{P}$.

%	We use this result to simplify the conditioning sets in step three by invoking the second statement from Lemma \ref{Wermuth} and dropping all nodes in $\mathbf{P}$. 
	
	We now prove Statement \eqref{out}. For $M=A$ the statement is trivial. Hence consider a node $M \notin \mathbf{W'}=\mathbf{W} \cup \{X\}$ and its corresponding coefficient $\gamma_m$. For ease of notation let $W_{k+1}=A$. We will now show that for $M \notin \mathbf{W'}$, it holds that $\gamma_m = \sum_{W_j \in \mathbf{W'}} \tau_{{w_j}m.\mathbf{w'}_{-j}} \gamma_{w_j}$.
%	, where $\tau_{{w_i}m.\mathbf{w'}_{-i}}$ is the sum of the product of the edge coefficient along all causal paths from $M$ to $W_i$ that do not contain nodes in $\mathbf{W'}_{-i}$; a partial total effect. 
Using equation \eqref{eq:gamma}, this claim is equivalent to
	\begin{align*}
	\tau_{am}  - \sum_{W_i \in \mathbf{W}} \beta_{aw_i.\mathbf{w}_{-i}} \tau_{w_im} 
	&= \sum_{W_j \in \mathbf{W'}} \tau_{{w_j}m.\mathbf{w'}_{-j}} (\tau_{aw_j} -  \sum_{W_i \in \mathbf{W}} \beta_{aw_i.\mathbf{w}_{-i}} \tau_{w_iw_j}). 
	\end{align*}
	By Lemma \ref{lemma: TE-decomposition}, $\tau_{w_im} = \sum_{W_j \in \mathbf{W'}} \tau_{w_jm.\mathbf{w'}_{-j}} \tau_{w_iw_j}$ for any $W_i \in \mathbf{W' }$ and thus, our claim follows. 
	
	The coefficient $\gamma_m$ can therefore only be non-zero, if at least one of the terms $\tau_{{w_j}m.\mathbf{w'}_{-j}} \gamma_{w_j}$ is also non-zero. Let $M' \in \mathbf{W'}$ be a node, such that $\tau_{m'm.\mathbf{w'}_{-m'}} \neq 0$ and $\gamma_{m'}\neq 0$. The first term being non-zero implies the existence of a directed path $p'$ from $M$ to $M'$ that contains no additional nodes from $\mathbf{W'}$ and is hence open given $\mathbf{W'}$ and $\mathbf{W}$. The second term being non-zero, implies that $\epsilon_{m'} \in \epsilon^{a\mathbf{w}}$, which by Statement \eqref{in} requires that there exists a path $p$ of the form
	$
	A \dots \rightarrow M',
	$
	that is open given $\mathbf{W}$, with possibly $A=M'$. Hence the first part of Statement \eqref{out} holds.

%	We now prove statement \ref{out}. By assumption $C \notin \mathbf{W}$ and hence the term $C$ does not appear in equation \eqref{deltas}. 
%	Thus, $\epsilon_{c} \in \epsilon^{a\mathbf{w}}$ can only hold true if $\de(C,\g) \cap (\{A\} \cup \mathbf{W}) \neq \emptyset$. Furthermore, if all the leading coefficients in equation \eqref{deltas} corresponding to the $\epsilon_{c'}$ with $C' \in (\de(C,\g) \cap (\{A\} \cup \mathbf{W}))$ are 0, so is the leading coefficient of $\epsilon_{c}$.  Hence, $\epsilon_{c} \in \epsilon^{a\mathbf{w}}$ can only hold, if there exists a node $C' \in \de(C,\g) \cap ( \{A\} \cup \mathbf{W})$, such that $\epsilon_{c'} \in \epsilon^{a\mathbf{w}}$, with possibly $C' = A$.
%    Then, by statement \ref{in}, there exists a path $p$ of the form
%	$
%	A \dots \rightarrow C',
%	$
%	that is open given $\mathbf{W}$ and
%	by construction there exists a causal path $p'$ from $C$ to $C'$,
%	with possibly $A=C'$. 
%	By choosing $C'$ as close to $C$ as possible, we can assure that the path $p'$ is open given $\mathbf{W}$.
 	 We now prove the second part of Statement \eqref{out}. If $M'=A$, $p'$ is a path of the claimed form, so suppose that $M' \in \mathbf{W}$.
	Let $I$ be the node closest to $A$ at which $p$ and $p'$ intersect, and consider the path $q=p(A,I) \oplus p'(I,M)$. We will now show that $q$ is open given $\mathbf{W}$.  If $I=M$ of $I=A$, $q$ is a subpath of $p$ or $p'$ respectively and as both $p$ and $p'$ are open given $\mathbf{W}$ we are done. 
	
	Hence, suppose that $I\notin \{A,M\}$. 
	As $p$ and $p'$ are open given $\mathbf{W}$, it suffices to consider $I$. Suppose first that $I \in \mathbf{W}$. Since $p'$ is directed and open given $\mathbf{W}$ it thus follows that $I = M'$. Then $M' \in \mathbf{W}$ is a collider on $q$ and it follows that $q$ is open given $\mathbf{W}$. 
	Suppose now that $I \notin \mathbf{W}$. Since $p'$ is directed towards $M'$ and $M' \in \mathbf{W}$ it follows that $\de(I,\g) \cap \mathbf{W} \neq \emptyset$. Hence, $q$ is open given $\mathbf{W}$ independently of whether $I$ is a collider or a non-collider.
	
	%Thus it suffices to show that the subpath $p(C',C)$ is open given $\mathbf{Z}$ to show that $p$ is open given $\mathbf{Z}$. Assume this is not the case. Then there exits a node $C'' \in \mathbf{Z}$ on $p(C',C)$, such that $\epsilon_{c''} \notin \epsilon_{x\mathbf{z}}$ by construction. But then $\epsilon_{c} \notin \epsilon_{x\mathbf{z}}$, unless there is either another directed path from $C$ to $C'$ or an alternative $C'$. Replacing either $C'$ or $C''$ and applying the same logic shows, that there has to exists a path of the form
	%\[
	%X \dots \rightarrow C' \leftarrow \dots \leftarrow C,
	%\]
	%that is open given $\mathbf{Z}$.
	
\end{proof}

\begin{Lemma}
	Let $X$ and $Y$ be nodes in a causal DAG $\g$. Let $\mathbf{Z}$ be a node set in $\g$, such that $Y \notin \mathbf{Z}, \pa(X,\g) \subseteq \mathbf{Z}$ and $\de(X,\g) \cap \mathbf{Z} = \emptyset$. Then $\mathbf{Z}$ is a valid adjustment set relative to $(X,Y)$.
	\label{lemma:alwaysVAS}
\end{Lemma}

\begin{proof}
	As a DAG, $\g$ is trivially amenable relative to $(X,Y)$. Further, $\f{\g} \subseteq \de(X,\g)$ and therefore $\mathbf{Z}$ fulfills the forbidden set condition $\eqref{cond1}$ from Definition $\ref{adjustment}$. 
	
	It only remains to show that $\mathbf{Z}$ blocks all non-causal paths from $X$ to $Y$ so let $p$ be such a path. Assume that $p$ is of the form $X \rightarrow \dots Y$. Then $p$ must contain a collider $C$, such that $C \in \de(X,\g)$. Since, by assumption $\de(X,\g) \cap \mathbf{Z} = \emptyset$ it follows that $p$ is blocked by $\mathbf{Z}$. Now, assume that $p$ is of the form $X \leftarrow \dots Y$. Then $p$ contains a non-collider $N \in \pa(X,\g)$ and is thus blocked by $\mathbf{Z}$.  
\end{proof}

\begin{Lemma}
	Let $A$ and $M$ be two nodes and $\mathbf{W}$ a node set in a DAG $\g$, such that $A \notin \mathbf{W}$ and $M \in \mathbf{W}$. Let $\mathbf{P} = \pa(M,\g) \setminus \mathbf{W}$. If no path from $A$ to $M$, ending with an edge into $M$ that is open given $\mathbf{W}$ exists, then $\mathbf{P} \perp_{\g} A | \mathbf{W}$.
	\label{lemma:edgeinto}
\end{Lemma}

\begin{proof}
We prove the claim by contraposition, so assume that a path $p$ from $A$ to some node $P \in \mathbf{P}$ that is open given $\mathbf{W}$ exists. By choice of $P$ there exists a path $p'$ of the form $P \rightarrow M$. Let $I$ be the node closest to $A$ on $p$ that is also on $p'$ and consider $q=p\oplus p'$. If $I=M$, then $q$ is a subpath of $p$ and hence open given $\mathbf{W}$. Further, $M \in \mathbf{W}$ must be a collider on $p$. But that implies that $q$ is a path from $A$ to $M$, ending with a node into $M$ that is open given $\mathbf{W}$. If $I=P$, then $I\notin \mathbf{W}$ is a non-collider on $q$ and our claim again follows. 
\end{proof}

\begin{Lemma}
	Consider a causal DAG $\g=(\mathbf{V},\mathbf{E})$ and let $\mathbf{V}$ follow a causal linear model compatible with $\g$. Let $N$ be a node and $\mathbf{A}=\{A_1,\dots, A_k\}$ be a node set in $\g$, such that $N \notin \mathbf{A}$. Then 
	\begin{equation}
	\tau_{a_in} = \sum_{A_j \in \mathbf{A}} \tau_{a_jn.\mathbf{a}_{-j}} \tau_{a_ia_j}
	\label{TE-decomposition}
	\end{equation}
	 for any $A_i \in \mathbf{A}$.
	 \label{lemma: TE-decomposition}
\end{Lemma}

\begin{proof}
	We first define two objects. Given two nodes $A$ and $B$ and a node set $\mathbf{C}$, let $\mathbf{P}^{ab.\mathbf{c}}$ denote the set of all directed paths from $B$ to $A$ not containing any nodes in $\mathbf{C}$. Further, given a directed path $p$, let $\tau^p$ denote the total effect along $p$, i.e. the product of the edge coefficients along $p$.
	
	We now prove our claim. Using the definition of the total effect via the path method we can rewrite the left hand term of equation \eqref{TE-decomposition} as
	\begin{align*}
	\tau_{a_in} = \sum_{p \in \mathbf{P}^{a_in}} \tau^p,
	\end{align*}
	and similarly, the right hand term as 
	\begin{align}
	\sum_{A_j \in \mathbf{A}} \tau_{a_jn.\mathbf{a}_{-j}} \tau_{a_ia_j}  &= \sum_{A_j \in \mathbf{A}} \left( \left( \sum_{p \in \mathbf{P}^{a_jn.a_{-j}}} \tau^{p} \right) \sum_{q \in \mathbf{P}^{a_ia_j}} \tau^{q} \right).  
\label{TE-decomposition2}
	\end{align}
	Clearly, the total effect along a directed path $q=p \oplus p'$ is equal to the the product of the total effect along $p$ and the total effect along $p'$. Using this we can rewrite equation \eqref{TE-decomposition2} as
	\begin{align*}
	\sum_{A_j \in \mathbf{A}} \tau_{a_jn.\mathbf{a}_{-j}} \tau_{a_ia_j}= \sum_{A_j \in \mathbf{A}} \sum_{p \in \mathbf{P}^{a_in}_{a_j.\mathbf{a}_{-j}}} \tau^{p},
	\end{align*}
	where $\mathbf{P}^{a_in}_{a_j.\mathbf{a}_{-j}}$ is the set of all directed paths $p$ from $N$ to $A_i$, such that $A_j$ lies on $p$ and $p(N,A_j)$ contains no node from $\mathbf{A}_{-j}$. Clearly, for any two nodes $A_j, A_k \in \mathbf{A}$, $\mathbf{P}^{a_in}_{a_j.\mathbf{a}_{-j}} \cap \mathbf{P}^{a_in}_{a_k.\mathbf{a}_{-k}} = \emptyset$. Since every directed path from $N$ to $A_i$ contains a node in $\mathbf{A}$, it follows that $\bigcup_{A_j \in \mathbf{A}} \mathbf{P}^{a_in}_{a_j.\mathbf{a}_{-j}}$ is a partition of $\mathbf{P}^{a_in}$ and therefore,
	\begin{align*}
		\tau_{a_in} = \sum_{p \in \mathbf{P}^{a_in}} \tau^p = 
		\sum_{A_j \in \mathbf{A}} \sum_{p \in \mathbf{P'}^{a_in}_{a_j}} \tau^{p} = 
		\sum_{A_j \in \mathbf{A}} \tau_{a_jn.\mathbf{a}_{-j}} \tau_{a_ia_j}. 
	\end{align*}
\end{proof}

\section{Proof of Theorem \ref{cor:bignew12}}

\begin{proofof}[Theorem \ref{cor:bignew12}] Let $\mathbf{X}=\{X_1,\dots,X_{k_x}\}$ and $\mathbf{Y}=\{Y_1,\dots,Y_{k_y}\}$ be disjoint node sets in a causal \mpdag{} $\g=(\mathbf{V},\mathbf{E})$ and let $\mathbf{V}$ follow a causal linear model compatible with $\g$. Let $\mathbf{Z_1}$ and $\mathbf{Z_2}$ be two valid adjustment sets relative to $(\mathbf{X,Y})$ in $\g$, such that $\mathbf{Y} \perp_{\g} (\mathbf{Z_1} \setminus \mathbf{Z_2}) | \mathbf{X} \cup \mathbf{Z_2} $ and $\mathbf{X} \perp_{\g} (\mathbf{Z_2} \setminus \mathbf{Z_1}) | \mathbf{Z_1}$.
	
	We first consider the case that $\g$ is a causal DAG. By applying Lemma \ref{close} with $\mathbf{T}=\mathbf{Z_1} \setminus \mathbf{Z_2},\mathbf{S}=\mathbf{Z_2} \setminus \mathbf{Z_1}$ and $\mathbf{W} = \mathbf{Z_1} \cap \mathbf{Z_2}$, it follows that
	\begin{align}
	&\sigma_{x_ix_i.\mathbf{x_{-i}\mathbf{z_1}}} \leq \sigma_{x_ix_i.\mathbf{x_{-i}\mathbf{z_2}}} \ \textrm{and} \nonumber \\
	&\sigma_{y_jy_j.\mathbf{xz_1}} \geq \sigma_{y_jy_j.\mathbf{xz_2}}, \nonumber 
	\end{align}
	for all $X_i \in \mathbf{X}$ and $Y_j \in \mathbf{Y}$. Using the asymptotic variance formula from Proposition \ref{avarDAG} it follows that
	\begin{align}
	a.var(\hat{\tau}^{\mathbf{z_2}}_{\mathbf{yx}})_{j,i} &= \frac{\sigma_{y_jy_j.\mathbf{xz_2}}}{\sigma_{x_ix_i.\mathbf{x}_{-i}\mathbf{z_2}}} \nonumber \\
	& \leq \frac{\sigma_{y_jy_j.\mathbf{xz_1}}}{\sigma_{x_ix_i.\mathbf{x}_{-i}\mathbf{z_1}}} 
	= a.var(\hat{\tau}^{\mathbf{z_1}}_{\mathbf{yx}})_{j,i}. \nonumber 
	\end{align}
	
	The proof then extends to the causal \mpdag{} setting with the fact that by Lemma \ref{lemma:dsepp1}, d-separation in a \mpdag{} implies d-separation in every represented DAG, including the true underlying one.
\end{proofof}

	In the multivariate Gaussian setting, the result of Theorem \ref{cor:bignew12} follows by Lemma \ref{close}  and the well known asymptotic variance formula from Lemma \ref{avar.equ} directly and does not require the new result from Proposition \ref{avarDAG}. In this setting it also holds for a larger class of sets, since Lemma \ref{avar.equ}, as opposed to Proposition 
	\ref{avarDAG}, does not require the conditioning set to be a valid adjustment set.

 Lemma \ref{lemma:dsepp1} shows that our definition of d-separation in \mpdag{}s is sensible, in the sense that it is compatible with d-separation in the DAGs represented by a \mpdag{}. It is analogous to Theorem 4.18 in \cite{RichardsonSpirtes02} for m-separation in maximal ancestral graphs (see also Lemma 20 in \citealp{zhang2008causal}) and Lemma 26 in \citet{zhang2008causal} for m-separation in partial ancestral graphs. 
 
 %Since the statement of Lemma \ref{lemma:dsepp1} for CPDAGs follows directly from the proof of Lemma 26 in \citet{zhang2008causal}, we only prove Lemma \ref{lemma:dsepp1} for maximal PDAGs.
 
 %Lemma 20 and 26 in \cite{zhang200 \ref{} in  follows from Lemma~\ref{lemma-2-pdag-dsep}, which itself relies on Lemma \ref{lemma-1-pdag-dsep}. Note that Lemma~\ref{lemma-2-pdag-dsep} is analogous to Lemma C.6 in \cite{perkovic17} and Lemma 2 in \cite{zhang2006causal}, while Lemma~\ref{lemma-1-pdag-dsep} is analogous to Lemma C.5 in \cite{perkovic17} and Lemma 1 in \cite{zhang2006causal}. 
 
 \begin{Lemma} \label{lemma:dsepp1}
 	Let $\mathbf{X}$, $\mathbf{Y}$ and $\mathbf{Z}$ be pairwise disjoint node sets in a \mpdag{} $\g $.
 	Then $\mathbf{X} \perp_{\g[D]} \mathbf{Y} | \mathbf{Z}$ in every DAG $\g[D] \in [\g]$, if and only if $\mathbf{Z}$ blocks every definite status path between any node in $\mathbf{X}$ and any node in $\mathbf{Y}$ in $\g$. 
 \end{Lemma}
 
 \begin{proofof}[Lemma  \ref{lemma:dsepp1}]
 	%Sets $\mathbf{X}$ and $\mathbf{Y}$ are d-connected given $\mathbf{Z}$ in at least one DAG $\g[D]$ in $[\g]$ if and only if there is a definite status path from $\mathbf{X}$ to $\mathbf{Y}$ that is open given $\mathbf{Z}$ in $\g$. 
 	
   We prove this statement by showing that the contrapositive statement is true. 
 	
 	Consider a definite status path $p$ from $\mathbf{X}$ to $\mathbf{Y}$ that is open given $\mathbf{Z}$ in $\g$ and a DAG $\g[D] \in [\g]$. Since $\g[D]$ is in the equivalence class described by  $\g$ it follows that $\g[D]$ has the same adjacencies as $\g$, and every edge $A \rightarrow B$ in $\g$ is also in $\g[D]$. 
 	Let $p^*$ be the corresponding path to $p$ in $\g[D]$. 
 	Since every node on $p$ is of definite status in $\g$, every node on $p^*$ is of the same definite status in $\g[D]$. Since additionally, $p$ is  open given $\mathbf{Z}$ and since for every $V \in \mathbf{V}$, $\de(V,\g) \subseteq \de(V,\g[D])$, $p^*$ is a definite status path from $\mathbf{X}$ to $\mathbf{Y}$ in $\g[D]$ that is open given $\mathbf{Z}$.
 	
 	Conversely, if there is a path $q$ from  $\mathbf{X}$ to $\mathbf{Y}$ that is open given $\mathbf{Z}$  in every DAG $\g[D] \in [\g]$, then by the proof of Lemma 26 in \citet{zhang2008causal} (see also the proof of Lemma 5.1.7 in \citealp{zhang2006causal}), there is a definite status path $q^*$ from $\mathbf{X}$ to $\mathbf{Y}$ that is open given $\mathbf{Z}$ in the CPDAG $\g[C]$ of any such $\g[D]$. 
 	Since $\g$ describes a subset of the Markov equivalence class of $[\g[C]]$ \citep{meek1995causal}, we can conclude with the same argument as above that the corresponding path $q^{**}$ of $q^*$ in $\g$ is a definite status path from $\mathbf{X}$ to $\mathbf{Y}$ in $\g$ that is open given $\mathbf{Z}$.

 \end{proofof}

%
%\subsection{Theorem \ref{A-cor:bignew12} in the multivariate Gaussian setting}
%
%In the multivariate Gaussian setting we can formulate Theorem \ref{A-cor:bignew12} more generally and drop the graphical setting, as the asymptotic variance formula from  Lemma \ref{avar.equ} always holds, as opposed to the one from Theorem \ref{avarDAG} used to prove Theorem \ref{A-cor:bignew12}.
%
%\begin{Theorem} \label{bignew}
%	Let $(\mathbf{X}^T,\mathbf{Y}^T,\mathbf{S}^T, \mathbf{T}^T,\mathbf{W}^T)^T$ be a multivariate Gaussian vector with mean $\mathbf{0}$, such that $\mathbf{X}=(X_1,\dots,X_{k_x})^T$ and $\mathbf{Y} = (Y_1,\dots,Y_{k_y})^T,$ 
%	%$k,l \ge 1$, 
%	with $\mathbf{T}, \mathbf{S}$ and $\mathbf{W}$ possibly of length zero. 
%	
%	If $\mathbf{T} \ci \mathbf{Y} | (\mathbf{W}^T,\mathbf{S}^T,\mathbf{X}^T)^T$ and $\mathbf{S} \ci \mathbf{X} | (\mathbf{W}^T,\mathbf{T}^T)^T$, 
%	then 
%	\[
%	a.var(\hat{\beta}_{y_jx_i.\mathbf{x}_{-i}\mathbf{ws}}) \leq 		
%	a.var(\hat{\beta}_{y_jx_i.\mathbf{x}_{-i}\mathbf{wt}}),
%	\]
%	for all $i \in \{1,\dots,k_x\}$ and $j \in \{1,\dots, k_y\}$.
%\end{Theorem}	
%
%
%Although this result seems rather straightforward and is commonly assumed to hold we did not find it in the literature.
%Note that \textit{not} all regression coefficients whose asymptotic variance can be compared with Theorem \ref{bignew} have to posses the same limit.
%
%\begin{proofof}[Theorem \ref{bignew}]
%By the asymptotic variance formula from Lemma \ref{avar.equ} the claim follows in the same way as in the proof of Theorem \ref{A-cor:bignew12}.
%\end{proofof}

\subsection{Residual linear variance  inequalities}

%\begin{Remark}
%	Throughout this section we use the notation that for multiple random vectors $\mathbf{S}, \mathbf{T}$ and $\mathbf{W_1}, \mathbf{W_2}, \dots, \mathbf{W_m}$, $m >1$, $\beta_{\mathbf{st}.{\mathbf{w_1w_2}\cdots \mathbf{w_m}}} = \beta_{\mathbf{st}.{\mathbf{w}}}$ and $\Sigma_{\mathbf{s}\mathbf{t}.{\mathbf{w_1w_2}\cdots \mathbf{w_m}}} = \Sigma_{\mathbf{s}\mathbf{t}.\mathbf{w}} $, where 
%	$\mathbf{W}=(\mathbf{W_1}^T, \dots, \mathbf{W_m}^T)^T$.
%\end{Remark}

By Lemmas \ref{lemma: buja} and \ref{buja} the asymptotic limit of a least squares regression is a function of the covariance matrix only, even when the regression is misspecified. Hence, the following statements and proofs are essentially linear algebra formulated in statistical terms. They do not depend on any property of the Gaussian distribution.

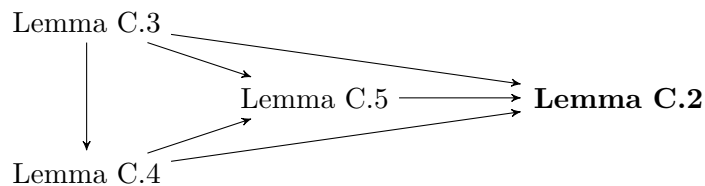
\begin{figure}[ht]
	\centering
	\begin{tikzpicture}[>=stealth',shorten >=1pt,node distance=3cm, main node/.style={minimum size=0.4cm}]
	\node[main node,yshift=0cm,xshift=0cm]         (T32) {\textbf{Lemma \ref{close}}};
		
	\node[main node,xshift=-4cm] (L3) at (T32) {Lemma \ref{Wermuth}};
	\node[main node,yshift=-1cm] (L2)  [left of= L3]   {Lemma \ref{convarsets}};   
	\node[main node,yshift=1cm]  (L1)  [left of= L3]   {Lemma \ref{Cochran}}; 
	
	%\node[main node,yshift=-1cm,xshift=-1cm]            
	%(C1)  at (T32)  {Corollary \ref{non-increasing}};  
	
	%\draw[->] (C1) edge    (T32);
	\draw[->] (L1) edge    (L2);
	\draw[->] (L1) edge    (L3);
	\draw[->] (L1) edge    (T32);
	\draw[->] (L2) edge    (L3);
	%\draw[->] (L2) edge    (C1);
	\draw[->] (L2) edge    (T32);
	\draw[->] (L3) edge    (T32);
	\end{tikzpicture}
	\caption{Proof structure of Lemma \ref{close}.}
	%\label{figproof}
\end{figure}

The following four Lemmas are simple generalizations of already existing results and are primarily given for completeness and conciseness; especially the latter three.
Lemma \ref{close} is a simple extension of Lemma 4 in \citet{kuroki2004selection} from random variables $X$ and $Y$  to random vectors $\mathbf{X}$ and $\mathbf{Y}$, while additionally also allowing $\mathbf{W}$ to be non-empty.
Lemma \ref{Cochran} is a an extension of a well known result from \citet{cochran1938omission} to vectors $\mathbf{X}, \mathbf{Y}$ and $\mathbf{Z}$.
Lemma \ref{convarsets} is a simple generalization of a result by \citet[Lemma 1]{kuroki2003covariate} from random variables $Z$ and $S$ to random vectors $\mathbf{Z}$ and $\mathbf{S}$. Note that it is quite similar to a result presented in Section 2.5 of \citet{anderson1958introduction}.
Lemma \ref{Wermuth} is a generalization of results by \citet[Results 1.2 and 5.2]{wermuth1989moderating} from random variables $X$, $Y$, $S$ and $T$ to random vectors $\mathbf{X},\mathbf{Y},\mathbf{S}$ and $\mathbf{T}$. Further, all of these Lemmas are also generalizations to non-Gaussian random variables using the result from Lemma \ref{lemma: buja}.

\begin{Lemma}
	Let $(\mathbf{X}^T,\mathbf{Y}^T,\mathbf{T}^T,\mathbf{S}^T, \mathbf{W}^T)^T$, with $\mathbf{T}, \mathbf{S}$ and $\mathbf{W}$ possibly of length zero, be
	 a mean $\mathbf{0}$ random vector with finite variance, such that $\mathbf{X}=(X_1,\dots,X_{k_x})^T$ and $\mathbf{Y} = (Y_1,\dots,Y_{k_y})^T$. If $\mathbf{T} \ci \mathbf{Y} | \mathbf{W},\mathbf{S},\mathbf{X}$ and $\mathbf{S} \ci \mathbf{X} | \mathbf{W},\mathbf{T}$, 
	then
	\begin{enumerate}[ label = (\alph*)]
		\item $\sigma_{x_ix_i.\mathbf{x_{-i}wt}} \leq \sigma_{x_ix_i.\mathbf{x_{-i}ws}} \ \textrm{and}$
		\label{x-variance1}
		\item $\sigma_{y_jy_j.\mathbf{xws}} \leq \sigma_{y_jy_j.\mathbf{xwt}},$
		\label{y-variance1}
	\end{enumerate}
for all $i \in \{1,\dots,k_x\}$ and $j \in \{1,\dots, k_y\}$.
	\label{close}
\end{Lemma}

\begin{proof}

We first assume that $\mathbf{T} = \emptyset$. Since $\mathbf{S} \ci \mathbf{X} | \mathbf{W}$ it must also hold that $\mathbf{S} \ci X_i | (\mathbf{W}^T, \mathbf{X_{-i}}^T)^T$ for all $X_i \in \mathbf{X}$, by the weak union property of conditional independence from \cite{dawid1979conditional}. Then,  by Lemma \ref{convarsets}, $\sigma_{y_jy_j.\mathbf{xws}} \leq \sigma_{y_jy_j.\mathbf{xw}}$ and by Lemma \ref{Wermuth}, $\sigma_{x_ix_i.\mathbf{x_{-i}w}} = \sigma_{x_ix_i.\mathbf{x_{-i}ws}}$. 

We now assume that $\mathbf{S} = \emptyset$. By Lemma \ref{convarsets}, $\sigma_{x_ix_i.\mathbf{x_{-i}wt}} \leq \sigma_{x_ix_i.\mathbf{x_{-i}w}}$ and by Lemma \ref{Wermuth}, $\sigma_{y_jy_j.\mathbf{xw}} = \sigma_{y_jy_j.\mathbf{xwt}}$.

We now assume $\mathbf{T} \neq \emptyset$ and $\mathbf{S} \neq \emptyset$.  
%because these imply
%%If $\sigma_{x_ix_i.\mathbf{x_{-i}wt}} \leq \sigma_{x_ix_i.\mathbf{x_{-i}ws}}$ and $\sigma_{y_jy_j.\mathbf{xws}} \leq \sigma_{y_jy_j.\mathbf{xwt}}$ hold for all $X_i \in \mathbf{X}$ and $Y_j \in \mathbf{Y}$,
%\begin{align*}
%a.var(\hat{\beta}_{\mathbf{yx}.\mathbf{ws}}) & = 
%\sum_{Y_j \in \mathbf{Y}} \sum_{X_i \in \mathbf{X}}\frac{\sigma_{y_jy_j.\mathbf{xws}}}{n\sigma_{x_ix_i.\mathbf{x_{-i}ws}}}\\ 
%  & \leq
%\sum_{Y_j \in \mathbf{Y}} \sum_{X_i \in \mathbf{X}} \frac{\sigma_{y_jy_j.\mathbf{xwt}}}{n\sigma_{x_ix_i.\mathbf{x_{-i}wt}}}
%  = a.var(\hat{\beta}_{\mathbf{yx}.\mathbf{wt}}).
%\end{align*}
First, we show inequality \ref{x-variance1}.  Since $\mathbf{S} \ci \mathbf{X} | (\mathbf{W}^T , \mathbf{T}^T)$ it also holds that $\mathbf{S} \ci X_i | (\mathbf{W}^T, \mathbf{T}^T, \mathbf{X}_{-i}^T)^T$ for all $X_i \in \mathbf{X}$ by the weak union property of conditional independence. % The fact that $\mathbf{X}$ and $\mathbf{S}$ are d-separated by $\mathbf{W} \cup \mathbf{T}$ implies that $X_i$ and $\mathbf{S}$ are d-separated by $\mathbf{W} \cup \mathbf{T} \cup \mathbf{X_{-i}}$ for all $X_i \in \mathbf{X}$ by Lemma \ref{d-sep-extension}. 
%NOTE: check if this is standard conditional independence result. 
Thus, $\boldsymbol{\beta}_{x_i\mathbf{s}.\mathbf{tx_{-i}w}} = 0$ by Lemma \ref{lemma: independence} and $\boldsymbol{\beta}_{x_i\mathbf{t}.\mathbf{sx_{-i}w}} = \boldsymbol{\beta}_{x_i\mathbf{t}.\mathbf{x_{-i}w}}$ by Lemma~\ref{Wermuth}. Hence, by Lemma \ref{Cochran}
\[
\begin{aligned}
\boldsymbol{\beta}_{x_i\mathbf{s}.\mathbf{x_{-i}w}} &= \boldsymbol{\beta}_{x_i\mathbf{s}.\mathbf{tx_{-i}w}} + \boldsymbol{\beta}_{x_i\mathbf{t}.\mathbf{sx_{-i}w}}\boldsymbol{\beta}_{\mathbf{ts}.\mathbf{x_{-i}w}} \\ &=\boldsymbol{\beta}_{x_i\mathbf{t}.\mathbf{x_{-i}w}} \boldsymbol{\beta}_{\mathbf{ts}.\mathbf{x_{-i}w}}.
\end{aligned}
\]
Then by Lemma \ref{convarsets},
\[
\begin{aligned}
\sigma_{x_ix_i.\mathbf{x_{-i}ws}} &- \sigma_{x_ix_i.\mathbf{x_{-i}wt}} \\
&=\boldsymbol{\beta}_{x_i\mathbf{t}.\mathbf{x_{-i}w}} \Sigma_{\mathbf{tt}.\mathbf{x_{-i}w}} \boldsymbol{\beta}^T_{x_i\mathbf{t}.\mathbf{x_{-i}w}}  -
\boldsymbol{\beta}_{x_i\mathbf{s}.\mathbf{x_{-i}w}} \Sigma_{\mathbf{ss}.\mathbf{x_{-i}w}} \boldsymbol{\beta}^T_{x_i\mathbf{s}.\mathbf{x_{-i}w}} \\
&= \boldsymbol{\beta}_{x_i\mathbf{t}.\mathbf{x_{-i}w}}( \Sigma_{\mathbf{tt}.\mathbf{x_{-i}w}} 
- \boldsymbol{\beta}_{\mathbf{ts}.\mathbf{x_{-i}w}} \Sigma_{\mathbf{ss}.\mathbf{x_{-i}w}} \boldsymbol{\beta}_{\mathbf{ts}.\mathbf{x_{-i}w}}^T)
\boldsymbol{\beta}^T_{x_i\mathbf{t}.\mathbf{x_{-i}w}} \\
&=  \boldsymbol{\beta}_{x_i\mathbf{t}.\mathbf{x_{-i}w}}
(\Sigma_{\mathbf{tt}.\mathbf{swx_{-i}}})
\boldsymbol{\beta}^T_{x_i\mathbf{t}.\mathbf{x_{-i}w}} \geq 0.
\end{aligned} 
\]

For inequality $\ref{y-variance1}$ we use that by Lemma~\ref{Wermuth}, $\mathbf{Y} \ci \mathbf{T} | (\mathbf{W}^T, \mathbf{S}^T , \mathbf{X}^T)^T$ implies that
$\boldsymbol{\beta}_{y_j\mathbf{t}.\mathbf{xws}} = 0$ by Lemma \ref{lemma: independence} and $\boldsymbol{\beta}_{y_j\mathbf{s}.\mathbf{txw}} = \boldsymbol{\beta}_{y_j\mathbf{s}.\mathbf{xw}}$ by Lemma \ref{Wermuth}. Hence, by Lemma \ref{Cochran}
\[
\boldsymbol{\beta}_{y_j\mathbf{t}.\mathbf{xw}} = \boldsymbol{\beta}_{y_j\mathbf{t}.\mathbf{sxw}} + \boldsymbol{\beta}_{y_j\mathbf{s}.\mathbf{txw}}
\boldsymbol{\beta}_{\mathbf{st}.\mathbf{xw}} =
\boldsymbol{\beta}_{y_j\mathbf{s}.\mathbf{xw}} 
\boldsymbol{\beta}_{\mathbf{st}.\mathbf{xw}}.
\]
Then by Lemma \ref{convarsets},
\begin{align}
\sigma_{y_jy_j.\mathbf{xwt}} - \sigma_{y_jy_j.\mathbf{xws}} &=\boldsymbol{\beta}_{y_j\mathbf{s}.\mathbf{xw}} \Sigma_{\mathbf{ss}.\mathbf{xw}} \boldsymbol{\beta}^T_{y_j\mathbf{s}.\mathbf{xw}} -
\boldsymbol{\beta}_{y_j\mathbf{t}.\mathbf{xw}} \Sigma_{\mathbf{tt}.\mathbf{xw}} \boldsymbol{\beta}^T_{y_j\mathbf{t}.\mathbf{xw}} \nonumber\\
&= \boldsymbol{\beta}_{y_j\mathbf{s}.\mathbf{xw}}( \Sigma_{\mathbf{ss}.\mathbf{xw}} - \boldsymbol{\beta}_{\mathbf{st}.\mathbf{xw}} \Sigma_{\mathbf{tt}.\mathbf{xw}} \boldsymbol{\beta}_{\mathbf{st}.\mathbf{xw}}^T)\boldsymbol{\beta}^T_{y_j\mathbf{s}.\mathbf{xw}} \nonumber\\
&=  \boldsymbol{\beta}_{y_j\mathbf{s}.\mathbf{xw}}(\Sigma_{\mathbf{ss}.\mathbf{xwt}})\boldsymbol{\beta}^T_{y_j\mathbf{s}.\mathbf{xw}} \geq 0. \label{y_variance_calc}
\end{align}

\end{proof}

\begin{Lemma}\label{Cochran}
	Let $\mathbf{V} = (\mathbf{S}^T,\mathbf{T}^T,\mathbf{W}^T,\mathbf{Z}^T)^T$, with $\mathbf{Z}$ possibly of length zero, be a mean $\mathbf{0}$ random vector with finite variance.
Then
%\begin{enumerate}[label=(\roman*)]
%\item
\begin{align}
\label{Cochran1} 
\boldsymbol{\beta}_{\mathbf{tw}.\mathbf{z}} = 
\boldsymbol{\beta}_{\mathbf{tw}.\mathbf{sz}} 
+ \boldsymbol{\beta}_{\mathbf{ts}.\mathbf{wz}} \boldsymbol{\beta}_{\mathbf{sw}.\mathbf{z}}.  
\end{align}

%\item\label{Cochran2} $
%\boldsymbol{\beta}_{\mathbf{tw}} = 
%\boldsymbol{\beta}_{\mathbf{tw}.\mathbf{s}} 
%+ \boldsymbol{\beta}_{\mathbf{ts}.\mathbf{w}} \boldsymbol{\beta}_{\mathbf{sw}}. 
%$
%\end{enumerate}
\end{Lemma}

\begin{proof}%of}[Lemma~\ref{A-Cochran}]
%As the regression coefficient values depend only on the distribution of $\mathbf{V}$ via the covariance matrix $\Sigma_{\mathbf{vv}}$ (see Lemma \ref{lemma: buja}), we can assume without loss of generality that the distribution of $\mathbf{V}$ is multivariate Gaussian.
%	
This proof is based on the uniqueness of the least squares regression. Precisely, by a projection argument, it holds that for any random vector $(\mathbf{Y},\mathbf{X}^T)^T$ with
$\mathbf{Y}=\boldsymbol{\beta}_{\mathbf{yx}} \mathbf{X} + \boldsymbol{\epsilon},$
the least squares regression coefficient $\boldsymbol{\beta}_{\mathbf{yx}}$ is characterized by the property that $\mathbb{E}[\boldsymbol{\epsilon} \mathbf{X}^T] = \boldsymbol{0}$.

Now suppose first that $\mathbf{Z} \neq \emptyset$. Regressing $\mathbf{S}$ on $(\mathbf{Z}^T,\mathbf{W}^T)^T$ yields
%We now apply equation \eqref{multinomial} with $\mathbf{Y}=\mathbf{S}$ and $\mathbf{X}=(\mathbf{Z}^T,\mathbf{W}^T)^T$ so that
\begin{equation}
\mathbf{S} = \boldsymbol{\beta}_{\mathbf{sz}.\mathbf{w}} {\mathbf{Z}} 
+ \boldsymbol{\beta}_{\mathbf{sw}.\mathbf{z}} {\mathbf{W}} 
+ \boldsymbol{\epsilon}_{\mathbf{s}},
\label{jequations}
\end{equation}
with $\mathbb{E}[\boldsymbol{\epsilon}_{\mathbf{s}} (\mathbf{Z}^T,\mathbf{W}^T)] = \boldsymbol{0}$. 
%and let $\boldsymbol{\epsilon}_{\mathbf{s}} = (\epsilon_1, \dots, \epsilon_{k_s})^T$. 
Similarly, regressing $\mathbf{T}$ on $(\mathbf{Z}^T,\mathbf{W}^T,\mathbf{S}^T)^T$ yields
\begin{align}
\mathbf{T} &= \boldsymbol{\beta}_{\mathbf{tz}.\mathbf{s}\mathbf{w}} \mathbf{Z} +
\boldsymbol{\beta}_{\mathbf{tw}.\mathbf{s}\mathbf{z}} \mathbf{W} +
\boldsymbol{\beta}_{\mathbf{ts}.\mathbf{w}\mathbf{z}} \mathbf{S} + \boldsymbol{\epsilon}_{\mathbf{t}},  \label{eq:ti1}
\end{align}
with $\mathbb{E}[\boldsymbol{\epsilon}_{\mathbf{t}} (\mathbf{Z}^T,\mathbf{W}^T,\mathbf{S}^T)] = \boldsymbol{0}$. 
Substituting equation \eqref{jequations} into equation \eqref{eq:ti1} gives
\begin{align}
\mathbf{T} &= \boldsymbol{\beta}_{\mathbf{tz}.\mathbf{s}\mathbf{w}} \mathbf{Z} + 
\boldsymbol{\beta}_{\mathbf{tw}.\mathbf{s}\mathbf{z}} \mathbf{W}
+ \boldsymbol{\beta}_{\mathbf{ts}.\mathbf{w}\mathbf{z}} (\boldsymbol{\beta}_{\mathbf{s}\mathbf{z}.\mathbf{w}} \mathbf{Z} + \boldsymbol{\beta}_{\mathbf{s}\mathbf{w}.\mathbf{z}} \mathbf{W} + \boldsymbol{\epsilon}_{\mathbf{s}}) + \boldsymbol{\epsilon}_{\mathbf{t}} \nonumber\\
\begin{split}
&= (\boldsymbol{\beta}_{\mathbf{tz}.\mathbf{s}\mathbf{w}} + 
\boldsymbol{\beta}_{\mathbf{ts}.\mathbf{w}\mathbf{z}} \boldsymbol{\beta}_{\mathbf{s}\mathbf{z}.\mathbf{w}}) \mathbf{Z} + (\boldsymbol{\beta}_{\mathbf{tw}.\mathbf{s}\mathbf{z}}  + \boldsymbol{\beta}_{\mathbf{ts}.\mathbf{w}\mathbf{z}} \boldsymbol{\beta}_{\mathbf{s}\mathbf{w}.\mathbf{z}}) \mathbf{W}  + \boldsymbol{\beta}_{\mathbf{ts}.\mathbf{wz}} \boldsymbol{\epsilon}_{\mathbf{s}} 
+ \boldsymbol{\epsilon}_{\mathbf{t}}.  
\end{split}
\label{OLShard}
\end{align}
Letting $\boldsymbol{\tilde{\epsilon}}_{\mathbf{t}}=\boldsymbol{\beta}_{\mathbf{ts}.\mathbf{wz}} \boldsymbol{\epsilon}_{\mathbf{s}} + \boldsymbol{\epsilon}_{\mathbf{t}}$ it follows that 
$\mathbb{E}[\boldsymbol{\tilde{\epsilon}}_{\mathbf{t}}(\mathbf{Z}^T,\mathbf{W}^T)]=\boldsymbol{0}.$

On the other hand, regressing $\mathbf{T}$ on $(\mathbf{Z}^T,\mathbf{W}^T)^T$ directly yields
\begin{align}
\mathbf{T} &= \boldsymbol{\beta}_{\mathbf{tz}.\mathbf{w}} \mathbf{Z} +
\boldsymbol{\beta}_{\mathbf{tw}.\mathbf{z}} \mathbf{W} + \boldsymbol{{\epsilon}}_{\mathbf{t}}',  
\label{OLSsimple}
\end{align}
with $\mathbb{E}[\boldsymbol{\epsilon}_{\mathbf{t}}' (\mathbf{Z}^T,\mathbf{W}^T)] = \boldsymbol{0}$.

Comparing equations \eqref{OLShard} and \eqref{OLSsimple}, combined with the uniqueness of the least squares regression coefficient, implies that
\begin{equation}
\boldsymbol{\beta}_{\mathbf{tw}.\mathbf{z}} = \boldsymbol{\beta}_{\mathbf{tw}.\mathbf{s}\mathbf{z}} 
+ \boldsymbol{\beta}_{\mathbf{ts}.\mathbf{w}\mathbf{z}} \boldsymbol{\beta}_{\mathbf{s}\mathbf{w}.\mathbf{z}}.
\nonumber 
\end{equation}

If $\mathbf{Z}=\emptyset$, one can simply drop all terms involving $\mathbf{Z}$. 
%The proof in the case $\mathbf{Z}$ of length zero is the same, except for dropping all $\mathbf{Z}$ terms. 

\end{proof}%of}

\begin{Lemma}
Let $\mathbf{V}=(\mathbf{S}^T,\mathbf{W}^T,\mathbf{Z}^T)^T$, with $\mathbf{S}$ possibly of length zero, be a mean $\mathbf{0}$ random vector with finite variance. Then 
\begin{align*}
%\item\label{convar1} $\Sigma_{\mathbf{zz}.\mathbf{r}} = \Sigma_{\mathbf{zz}} - \boldsymbol{\beta}_{\mathbf{zr}} \Sigma_{\mathbf{rr}} \boldsymbol{\beta}^T_{\mathbf{zr}}$, and
\Sigma_{\mathbf{zz}.\mathbf{sw}} = \Sigma_{\mathbf{zz}.\mathbf{s}} - \boldsymbol{\beta}_{\mathbf{zw}.\mathbf{s}} \Sigma_{\mathbf{ww}.\mathbf{s}} \boldsymbol{\beta}^T_{\mathbf{zw}.\mathbf{s}}.
\end{align*}
\label{convarsets}
\end{Lemma}

\begin{proof}%of}[Lemma~\ref{A-convarsets}] 
Let $\mathbf{R} = (\mathbf{S}^T,\mathbf{W}^T)^T$. By Lemma \ref{lemma: buja} it holds that
 $\Sigma_{\mathbf{zr}}=\boldsymbol{\beta}_{\mathbf{zr}}\Sigma_{\mathbf{rr}}$. Combining this with the fact that
$
    \Sigma_{\mathbf{zz}.\mathbf{r}} = \Sigma_{\mathbf{zz}} - \Sigma_{\mathbf{zr}} \Sigma^{-1}_{\mathbf{rr}} \Sigma^T_{\mathbf{zr}}$, it follows that 
 \begin{align}  
  \Sigma_{\mathbf{zz}} = \Sigma_{\mathbf{zz}.\mathbf{r}} - \boldsymbol{\beta}_{\mathbf{zr}} \Sigma_{\mathbf{rr}} \boldsymbol{\beta}^T_{\mathbf{zr}}.
    \label{ConVarStandard}
\end{align}

If $\mathbf{S}=\emptyset$ then $\mathbf{R}=\mathbf{W}$ and our claim follows. So suppose that $\mathbf{S}\neq \emptyset$. Since $\mathbf{R} = (\mathbf{S}^T,\mathbf{W}^T)^T$, $\Sigma_{\mathbf{zz}.\mathbf{sw}} = \Sigma_{\mathbf{zz}.\mathbf{r}}$. Note that 
$\Sigma_{\mathbf{rr}} = \begin{pmatrix} \Sigma_{\mathbf{ss}}& \Sigma_{\mathbf{sw}}\\\Sigma_{\mathbf{ws}}&\Sigma_{\mathbf{ww}}\end{pmatrix}  $ 
and
$ \boldsymbol{\beta}_{\mathbf{z}\mathbf{r}} = \begin{pmatrix} \boldsymbol{\beta}_{\mathbf{zs}.\mathbf{w}} & \boldsymbol{\beta}_{\mathbf{zw}.\mathbf{s}} \end{pmatrix}.$
Plugging this into equation \eqref{ConVarStandard} yields
\begin{align}
\begin{split}
\Sigma_{\mathbf{zz}.\mathbf{sw}} = 
\Sigma_{\mathbf{zz}.\mathbf{r}} &=
\Sigma_{\mathbf{zz}} - \boldsymbol{\beta}_{\mathbf{zw}.\mathbf{s}} \Sigma_{\mathbf{ww}} \boldsymbol{\beta}^{T}_{\mathbf{zw}.\mathbf{s}}
- \boldsymbol{\beta}_{\mathbf{zw}.\mathbf{s}} \Sigma_{\mathbf{ws}} \boldsymbol{\beta}^{T}_{\mathbf{zs}.\mathbf{w}}\\
&\quad \quad   - \boldsymbol{\beta}_{\mathbf{zs}.\mathbf{w}} \Sigma_{\mathbf{sw}} \boldsymbol{\beta}^{T}_{\mathbf{zw}.\mathbf{s}} 
-\boldsymbol{\beta}_{\mathbf{zs}.\mathbf{w}} \Sigma_{\mathbf{ss}} \boldsymbol{\beta}^{T}_{\mathbf{zs}.\mathbf{w}}. 
\end{split}
\label{eq:convar1}
\end{align}

Using $\Sigma_{\mathbf{ws}} = \boldsymbol{\beta}_{\mathbf{ws}} \Sigma_{\mathbf{ss}}$ and $\Sigma_{\mathbf{sw}} = \Sigma_{\mathbf{ss}} \boldsymbol{\beta}^T_{\mathbf{ws}}$, we can rewrite equation \eqref{eq:convar1} as
\begin{align}
\begin{split}
\Sigma_{\mathbf{zz}.\mathbf{sw}} &= 
\Sigma_{\mathbf{zz}} - \boldsymbol{\beta}_{\mathbf{zw}.\mathbf{s}} \Sigma_{\mathbf{ww}} \boldsymbol{\beta}^{T}_{\mathbf{zw}.\mathbf{s}} 
+ \boldsymbol{\beta}_{\mathbf{zw}.\mathbf{s}} \boldsymbol{\beta}_{\mathbf{ws}} \Sigma_{\mathbf{ss}} \boldsymbol{\beta}^T_{\mathbf{ws}} \boldsymbol{\beta}^T_{\mathbf{zw}.\mathbf{s}} \\
&\quad\quad - (\boldsymbol{\beta}_{\mathbf{zs}.\mathbf{w}} + \boldsymbol{\beta}_{\mathbf{zw}.\mathbf{s}} \boldsymbol{\beta}_{\mathbf{ws}}) \Sigma_{\mathbf{ss}}
(\boldsymbol{\beta}_{\mathbf{zs}.\mathbf{w}} + \boldsymbol{\beta}_{\mathbf{zw}.\mathbf{s}} \boldsymbol{\beta}_{\mathbf{ws}})^T. \label{eq:convar2}
\end{split}
\end{align}

By Lemma \ref{Cochran} it holds that $\boldsymbol{\beta}_{\mathbf{zs}} = \boldsymbol{\beta}_{\mathbf{zs}.\mathbf{w}} + \boldsymbol{\beta}_{\mathbf{zw}.\mathbf{s}} \boldsymbol{\beta}_{\mathbf{ws}}$. Plugging this into equation \eqref{eq:convar2} and then using equation \eqref{ConVarStandard} twice, we arrive at
\begin{align*}
\Sigma_{\mathbf{zz}.\mathbf{sw}} &= 
\Sigma_{\mathbf{zz}} - \boldsymbol{\beta}_{\mathbf{zs}} \Sigma_{\mathbf{ss}} \boldsymbol{\beta}^{T}_{\mathbf{zs}} 
- \boldsymbol{\beta}_{\mathbf{zw}.\mathbf{s}}(\Sigma_{\mathbf{ww}} - \boldsymbol{\beta}_{\mathbf{ws}} \Sigma_{\mathbf{ss}} \boldsymbol{\beta}^T_{\mathbf{ws}} )\boldsymbol{\beta}^{T}_{\mathbf{zw}.\mathbf{s}} \\
&= \Sigma_{\mathbf{zz}.\mathbf{s}} - \boldsymbol{\beta}_{\mathbf{zw}.\mathbf{s}} \Sigma_{\mathbf{ww}.\mathbf{s}} \boldsymbol{\beta}^T_{\mathbf{zw}.\mathbf{s}}.
\end{align*}
\end{proof}%of}

%\begin{Corollary}
%Let $(X,\mathbf{W},\mathbf{S})$ be normally distributed and pairwise disjoint. Then $\sigma_{xx.\mathbf{ws}} \leq \sigma_{xx.\mathbf{s}}$.
%\label{non-increasing}
%\end{Corollary}

\begin{Lemma}
Let $(\mathbf{X}^T,\mathbf{Y}^T,\mathbf{S}^T,\mathbf{T}^T)^T$, with $\mathbf{S}$ possibly of length zero, be a mean $\mathbf{0}$ random vector with finite variance. If $\mathbf{T} \ci \mathbf{X} | \mathbf{S}$ or $\mathbf{T} \ci \mathbf{Y} | \mathbf{X}, \mathbf{S}$, then
$
\boldsymbol{\beta}_{\mathbf{yx}.\mathbf{s}} = \boldsymbol{\beta}_{\mathbf{yx}.\mathbf{st}}.
$
Furthermore, if $\mathbf{T} \ci \mathbf{Y}| \mathbf{X} , \mathbf{S}$, then
$
\Sigma_{\mathbf{yy}.\mathbf{xst}} = \Sigma_{\mathbf{yy}.\mathbf{xs}}.$

\label{Wermuth}
\end{Lemma}

\begin{proof} If $\mathbf{T} \ci \mathbf{X} | \mathbf{S}$ or $\mathbf{T} \ci \mathbf{Y} | \mathbf{X}, \mathbf{S}$, then $\boldsymbol{\beta}_{\mathbf{tx}.\mathbf{s}} =0$ or $\boldsymbol{\beta}_{\mathbf{yt}.\mathbf{xs}} =0$ respectively by Lemma \ref{lemma: independence}. Then, using equation  \eqref{Cochran1} in Lemma \ref{Cochran}, we have
\[
\boldsymbol{\beta}_{\mathbf{yx}.\mathbf{s}} = \boldsymbol{\beta}_{\mathbf{yx}.\mathbf{st}} + \boldsymbol{\beta}_{\mathbf{yt}.\mathbf{xs}} \boldsymbol{\beta}_{\mathbf{tx}.\mathbf{s}} = \boldsymbol{\beta}_{\mathbf{yx}.\mathbf{st}}.
\]

Now, assume that $\mathbf{T} \ci \mathbf{Y}| \mathbf{X} , \mathbf{S}$. If $\mathbf{S}\neq \emptyset$, let $\mathbf{S}^{\prime}= (\mathbf{X}^T,\mathbf{S}^T)^T$ and otherwise let $\mathbf{S}' =\mathbf{X}$. By Lemma \ref{lemma: independence}, $\boldsymbol{\beta}_{\mathbf{yt}.\mathbf{s}^{\prime}}=0$ and by the already shown statements we have
$
\boldsymbol{\beta}_{\mathbf{ys}^{\prime}} = 
\boldsymbol{\beta}_{\mathbf{ys}^{\prime}.\mathbf{t}}.$
With $\boldsymbol{\beta}_{\mathbf{yt}.\mathbf{s}^{\prime}}=0$, $
\boldsymbol{\beta}_{\mathbf{ys}^{\prime}} = 
\boldsymbol{\beta}_{\mathbf{ys}^{\prime}.\mathbf{t}}$ and Lemma \ref{convarsets} it follows that
\begin{align*}
\Sigma_{\mathbf{yy}.\mathbf{xst}} &= \Sigma_{\mathbf{yy}.\mathbf{s}^{\prime}\mathbf{t}}
 = \Sigma_{\mathbf{yy}} - 
\begin{pmatrix}\boldsymbol{\beta}_{\mathbf{y}\mathbf{s}^{\prime}.\mathbf{t}} & \boldsymbol{\beta}_{\mathbf{y}\mathbf{t}.\mathbf{s}^{\prime}} \end{pmatrix}
\left(
\begin{array}{cc}
\Sigma_{\mathbf{s}^{\prime}\mathbf{s}^{\prime}}  & \Sigma_{\mathbf{s}^{\prime}\mathbf{t}} \\
\Sigma_{\mathbf{t}\mathbf{s}^{\prime}} & \Sigma_{\mathbf{t}\mathbf{t}} 
\end{array}
\right)
\begin{pmatrix}\boldsymbol{\beta}_{\mathbf{y}\mathbf{s}^{\prime}.\mathbf{t}}^T \\ \boldsymbol{\beta}_{\mathbf{y}\mathbf{t}.\mathbf{s}^{\prime}}^T \end{pmatrix}
\nonumber\\ 
&= \Sigma_{\mathbf{yy}} - \boldsymbol{\beta}_{\mathbf{y}\mathbf{s}^{\prime}.\mathbf{t}} \Sigma_{\mathbf{s}^{\prime}\mathbf{s}^{\prime}} \boldsymbol{\beta}^T_{\mathbf{y}\mathbf{s}^{\prime}.\mathbf{t}}  - \boldsymbol{\beta}_{\mathbf{y}\mathbf{t}.\mathbf{s}^{\prime}} \Sigma_{\mathbf{t}\mathbf{s}^{\prime}} \boldsymbol{\beta}^{T}_{\mathbf{y}\mathbf{s}^{\prime}.t} \nonumber \\
& \quad - \boldsymbol{\beta}_{\mathbf{y}\mathbf{s}^{\prime}.\mathbf{t}}  \Sigma_{\mathbf{s}^{\prime}\mathbf{t}} \boldsymbol{\beta}^T_{\mathbf{y}\mathbf{t}.\mathbf{s}^{\prime}}
- \boldsymbol{\beta}_{\mathbf{y}\mathbf{t}.\mathbf{s}^{\prime}}  \Sigma_{\mathbf{tt}} \boldsymbol{\beta}^T_{\mathbf{y}\mathbf{t}.\mathbf{s}^{\prime}}
\nonumber\\
&= \Sigma_{\mathbf{yy}} - \boldsymbol{\beta}_{\mathbf{y}\mathbf{s}^{\prime}} \Sigma_{\mathbf{s}^{\prime}\mathbf{s}^{\prime}} \boldsymbol{\beta}^T_{\mathbf{y}\mathbf{s}^{\prime}} 
= \Sigma_{\mathbf{yy}.\mathbf{s}^{\prime}} =\Sigma_{\mathbf{yy}.\mathbf{xs}}.
\end{align*}

\end{proof}

\subsection{Proof of Corollaries \ref{cor:badparents}, \ref{cor:goodparents} and \ref{cor:pretreatment}}

\begin{proofof}[Corollary \ref{cor:badparents}]
Let $\mathbf{X}$ and $\mathbf{Y}$ be disjoint node sets in a \mpdag{} $\g = (\mathbf{V},\mathbf{E})$ and let $\mathbf{V}$ follow a causal linear model compatible with $\g$. Let $\mathbf{Z}$ and $\mathbf{Z'} = \mathbf{Z} \setminus \{P\}$, with $P \in (\pa(\mathbf{X},\g) \cap \mathbf{Z})$, be each a valid adjustment set relative to $(\mathbf{X},\mathbf{Y})$ in $\g$.

Consider a DAG $\g[D]$ compatible with $\g$. By the completeness of the adjustment criterion, $\mathbf{Z}$ and $\mathbf{Z'}$ are also valid adjustment sets relative to $(\mathbf{X},\mathbf{Y})$ in $\g[D]$. Since $\pa(\mathbf{X},\g) \subseteq \pa(\mathbf{X},\g[D])$ it also follows that $P \in (\pa(\mathbf{X},\g[D]) \cap \mathbf{Z})$. We can therefore without loss of generality consider $\g[D]$ rather than $\g$.
	
We apply Theorem \ref{cor:bignew12} with $\mathbf{Z_1} = \mathbf{Z}$ and $\mathbf{Z_2} = \mathbf{Z'}$. Since $\mathbf{Z'} \subset \mathbf{Z}$ it only remains to show that $\{P\} \perp_{\g[D]} \mathbf{Y} | \mathbf{X} \cup \mathbf{Z'}$. We prove this by showing that the existence of a path from $P$ to $\mathbf{Y}$ that is open given $\mathbf{X} \cup \mathbf{Z'}$ in $\g[D]$ contradicts the assumption that $\mathbf{Z'}$  is a valid adjustment set. 
Let $p$ be such a path and $p'$ be the path $X \leftarrow P, X\in\mathbf{X}$ which exists by construction. Let $I$ be the node closest to X on $p'$ which also lies on $p$ and consider $q = p'(X,I) \oplus p(I,Y)$. If $I=X$, then $q$ is a subpath of $p$ and is hence open given $\mathbf{X}\cup \mathbf{Z'}$. As we assume $\mathbf{Z'}$ to be a valid adjustment set, $q$ must also be open given $\mathbf{Z'}$ by Lemma \ref{lemma:collider} and therefore has to be causal. This, however, implies that $X$ is a non-collider on $p$, which contradicts our starting assumption that $p$ is open given $\mathbf{X}\cup \mathbf{Z'}$. If $I=P$, then $q$ is a non-causal path from $\mathbf{X}$ to $\mathbf{Y}$. Since $P \notin \mathbf{X} \cup \mathbf{Z'}$ is a non-collider on $q$, $q$ is open given $\mathbf{X} \cup \mathbf{Z'}$ and hence by Lemma \ref{lemma:collider} and our assumptions, it must also be open given $\mathbf{Z'}$. Further, $q$ is either proper or any node from $\mathbf{X}$ it contains, is a collider. Let $X'$ be the node closest to $Y$ on $q$. Then $q(X',Y)$ is a proper, non-causal path from $\mathbf{X}$ to $\mathbf{Y}$ that is open given $\mathbf{Z'}$.

%This holds as we assume that $\mathbf{Z'}$ is a valid adjustment set.
\end{proofof}

\begin{proofof}[Corollary \ref{cor:goodparents}]
	Let $\mathbf{X}$ and $\mathbf{Y}$ be disjoint node sets in a \mpdag{} $\g = (\mathbf{V},\mathbf{E})$ and let $\mathbf{V}$ follow a causal linear model compatible with $\g$. Let $\mathbf{Z}$  be a valid adjustment set relative to $(\mathbf{X},\mathbf{Y})$ in $\g$ and let $\mathbf{Z'} = \mathbf{Z} \cup \{R\}$, with $R \in \pa(\mathbf{Y},\g)$. 
	
	Consider a DAG $\g[D]$ compatible with $\g$. By the completeness of the adjustment criterion, $\mathbf{Z}$ and $\mathbf{Z'}$ are also valid adjustment sets relative to $(\mathbf{X},\mathbf{Y})$ in $\g[D]$. Since $\pa(\mathbf{Y},\g) \subseteq \pa(\mathbf{Y},\g[D])$, it follows that $R \in \pa(\mathbf{Y},\g[D]) \setminus \mathbf{Z}$. We can therefore without loss of generality consider $\g[D]$ rather than $\g$.

	We now apply Theorem \ref{cor:bignew12} with $\mathbf{Z_1} = \mathbf{Z}$ and $\mathbf{Z_2} = \mathbf{Z'}$. Since $\mathbf{Z} \subset \mathbf{Z'}$ it only remains to show that $\{R\} \perp_{\g} \mathbf{X} | \mathbf{Z}$. We will now show that the existence of a non-causal path $p$ from $\mathbf{X}$ to $R$ that is open given $\mathbf{Z}$ contradicts the assumption that $\mathbf{Z}$ is a valid adjustment set. Since the existence of a causal path would contradict that $R \notin \fb{\g[D]}$, implicit to the assumption that $\mathbf{Z'}$ is a valid adjustment set, this suffices to show our claim.
	
%	We can therefore without loss of generality consider $\g[D]$ rather than $\g$.
%	We will first show that $\mathbf{Z'}$ is a valid adjustment set relative to $(\mathbf{X},\mathbf{Y})$ in $\g[D]$. Clearly, $\mathbf{Z'} \cap \fb{\g} = \emptyset$ so it remains to show that $\mathbf{Z'}$ blocks all proper non-causal paths from $\mathbf{X}$ to $\mathbf{Y}$. Assume that there exists a path $p$ contradicting this. 
%	We now show that the existence of such a path $p$ contradicts the assumption that $\mathbf{Z}$ is a valid adjustment set. 
	
	Consider a non-causal path $p$ from $\mathbf{X}$ to $R$ and suppose that it is open given the valid adjustment set $\mathbf{Z}$. Let $X \in \mathbf{X}$ be the last such node on $p$. Suppose that $p(X,R)$ is directed from $X$ to $R$ and therefore $R \in \de(\mathbf{X},\g)$. Since $R \in \pa(\mathbf{Y},\g)$ this contradicts that $R \notin \fb{\g}$. Hence, $p'=p(X, \dots, R)$ must be a proper, non-causal path from $\mathbf{X}$ to $R$ that is open given $\mathbf{Z}$. 
	Let $p''$ be the path $R \rightarrow Y$ that exists by assumption and let $I$ be the node closest to $Y$ on $p''$ which also lies on $p'$ and consider $q = p'(X,I) \oplus p''(I,Y)$.  
	
	If $I=Y$, then $q$ is a subpath of $p'$ and hence open given $\mathbf{Z}$. Therefore, it has to be causal as otherwise its existence would contradict the assumption that $\mathbf{Z}$ is a valid adjustment set. Since $R \in \pa(Y,\g)$, $Y$ or one of its descendants is then a collider on $p'$. By the causality of $q$, $\de(Y,\g) \subseteq \fb{\g}$ and thus, $p' $ may not be open given $\mathbf{Z}$ yielding a contradiction.
	% by the assumption that $\mathbf{Z} \cap \fb{\g} = \emptyset$, yielding a contradiction. 
	
	If $I=R$, then $R \notin \mathbf{Z}$ is a non-collider on $q$ and hence $q$ is a proper path from $\mathbf{X}$ to $\mathbf{Y}$ that is open given $\mathbf{Z}$. Since $p'$ is non-causal and a subpath of $q$, it follows that $q$ is non-causal. Hence, $q$ is a proper non-causal path from $\mathbf{X}$ to $\mathbf{Y}$ that is open given $\mathbf{Z}$ again yielding a contradiction.
	
%	 Thus $\mathbf{Z'}$ is a valid adjustment set relative to $(\mathbf{X},\mathbf{Y})$ in $\g[D]$. Since this holds for any DAG compatible with $\g$ we can conclude by the completeness of the adjustment criterion for \mpdag{}s that $\mathbf{Z'}$ is a valid adjustment relative to $(\mathbf{X},\mathbf{Y})$ in $\g$ (see Lemma C.3 in \citet{perkovic17} and Lemma \ref{lemma:forbIsSame}).

%We prove that this is the case by showing that the existence of a path $p$ from $\mathbf{X}$ to $R$ that is open given $\mathbf{Z}$ contradicts the assumption that $\mathbf{Z}$ is a valid adjustment set. Let $p'$ be the path $R \rightarrow Y$. Let $I$ be the node closest to $Y$ on $p'$ which also lies on $p$ and consider $q = p(X,I) \oplus p'(I,Y)$.  
%
%If $I=Y$, then $q$ is a subpath of $p$ and hence open given $\mathbf{Z}$. Therefore, it has to be causal as otherwise its existence would contradict the assumption that $\mathbf{Z}$ is a valid adjustment set. Since $R \in \pa(Y,\g)$, $Y$ is then a collider on $p$. By the causality of $q$, $Y \in \fb{\g}$ and thus, $p$ may not be open given $\mathbf{Z}$.
%% by the assumption that $\mathbf{Z} \cap \fb{\g} = \emptyset$, yielding a contradiction. 
%
%If $I=R$, then $R \notin \mathbf{Z}$ is a non-collider on $q$ and hence $q$ is a proper path from $\mathbf{X}$ to $\mathbf{Y}$ that is open given $\mathbf{Z}$. Again it follows by the assumption that $R \notin \fb{\g}$ that $q$ may not be causal. Hence, $q$ is a proper non-causal path from $\mathbf{X}$ to $\mathbf{Y}$ that is open given $\mathbf{Z}$.

\end{proofof}

\begin{proofof}[Corollary \ref{cor:pretreatment}]
	Let $\mathbf{X}$ and $\mathbf{Y}$ be disjoint node sets in a DAG $\g = (\mathbf{V},\mathbf{E})$, such that $\pa(\mathbf{X},\g)=\emptyset$ and let $\mathbf{V}$ follow a causal linear model compatible with $\g$. Let $\mathbf{Z}$ and $\mathbf{Z'}$ be node sets in $\g$, such that $\mathbf{Z} \cap (\de(\mathbf{X},\g)\cup\mathbf{Y}) = \emptyset$ and $\mathbf{Z'} \cap (\de(\mathbf{X},\g)\cup\mathbf{Y}) = \emptyset$. 
	
	We first show that $\mathbf{Z}$ and $\mathbf{Z'}$ are valid adjustment sets relative to $(\mathbf{X},\mathbf{Y})$ in $\g$. Clearly, $\fb{\g} \subseteq \de(\mathbf{X},\g)$ and therefore $\mathbf{Z} \cap \fb{\g} = \emptyset$. Let $p$ be a proper non-causal path from some node $X\in\mathbf{X}$ to some node $Y\in\mathbf{Y}$ in $\g$. By the assumption that $\pa(\mathbf{X},\g)=\emptyset$ it follows that $p$ must be of the form $X \rightarrow \dots Y$. As $p$ is non-causal this implies the existence of a collider $C$ on $\g$, such that $C \in \de(\mathbf{X},\g)$. As $\mathbf{Z} \cap \de(\mathbf{X},\g) = \emptyset$ it follows that $\mathbf{Z} \cap \de(C,\g)=\emptyset$ and thus $p$ is blocked by $\mathbf{Z}$. The same reasoning holds for $\mathbf{Z'}$.
	
%	Then $\mathbf{Z}'$ and $\mathbf{Z}'$ are valid adjustment sets relative to $(\mathbf{X},\mathbf{Y})$ in $\g$. 
Suppose now that $\mathbf{Z} \subseteq \mathbf{Z'}$. We will show that $\mathbf{Z'}\setminus \mathbf{Z} \perp_{\g} \mathbf{X} | \mathbf{Z}$. 
%	by showing that $\mathbf{Z'}\setminus \mathbf{Z} \not\perp_{\g} \mathbf{X} | \mathbf{Z},\pa(\mathbf{X},\g)=\emptyset,\pa(\mathbf{X},\g)=\emptyset,\mathbf{Z} \cap \de(\mathbf{X},\g) = \emptyset$ and $\mathbf{Z'}\cap \de(\mathbf{X},\g) = \emptyset$ leads to a contradiction. 
Consider a proper path $p$ from some node $Z' \in \mathbf{Z'}\setminus \mathbf{Z}$ to some node $X\in\mathbf{X}$. As $\pa(\mathbf{X},\g)=\emptyset$, the path $p$ must be of the form $Z' \dots \leftarrow X$. Since $Z'\notin \de(\mathbf{X},\g)$, $p$ must contain a collider $C$, such that $C \in \de(\mathbf{X},\g)$. But since $\mathbf{Z}\cap\de(\mathbf{X},\g)=\emptyset$ it follows that $\mathbf{Z} \cap \de(C,\g)=\emptyset$ and therefore $p$ is blocked given $\mathbf{Z}$. Our claim then follows by Theorem \ref{cor:bignew12}.

\end{proofof}

%\section{Results related to Algorithm \ref{A-Algorithm}}

\section{Proof of Theorem \ref{prop:order-indep}}

\label{section-D}

\begin{proofof}[Theorem \ref{prop:order-indep}]
	Let $\mathbf{X}$ and $\mathbf{Y}$ be disjoint node sets in a \mpdag{} $\g = (\mathbf{V},\mathbf{E})$ and let $\mathbf{Z}$ be a valid adjustment set relative to $(\mathbf{X,Y})$ in $\g$.

	%	 $\mathbf{Y} \perp_{\g} \mathbf{Z''_1} | \mathbf{Z'_1} \cup \mathbf{Z_2} \cup \mathbf{X}$ for any two disjoint subsets $\mathbf{Z'_1}, \mathbf{Z''_1} \subseteq \mathbf{Z_1}$. 

	As a simple corollary of Lemma \ref{lemma max size unique} there is a unique partition $\mathbf{Z}=\mathbf{Z_1} \cup \mathbf{Z_2}$ such that $\mathbf{Y} \perp_{\g} \mathbf{Z_1} | \mathbf{Z_2} \cup \mathbf{X}$ and $\mathbf{Z_1}$ is of maximal size or equivalently $\mathbf{Z_2}$ is of minimal size. 
	By Lemma \ref{lemma: weak union} $\mathbf{Y} \perp_{\g} \mathbf{Z''_1} | \mathbf{Z'_1} \cup \mathbf{Z_2} \cup \mathbf{X}$ for any two disjoint subsets $\mathbf{Z'_1}, \mathbf{Z''_1} \subseteq \mathbf{Z_1}$. 
%		and Lemma \ref{lemma: d-sep VAS} any two partitions $\mathbf{Z} = \mathbf{Z_1} \cup \mathbf{Z_2}=\mathbf{Z'_1} \cup \mathbf{Z'_2}$,
	Jointly, these two results imply that given any subset $\mathbf{Z'} \subseteq \mathbf{Z}$, fulfilling $\mathbf{Z_2} \subseteq \mathbf{Z'}$, it holds that for any node $Z \in \mathbf{Z'}$
	\begin{enumerate}
		\item if $Z \in \mathbf{Z'} \cap \mathbf{Z_1}$, then $\mathbf{Y} \perp_{\g} Z | \mathbf{Z'}_{-z} \cup \mathbf{X}$ and 
		\item if $Z \in \mathbf{Z'} \cap \mathbf{Z_2}$, then $\mathbf{Y} \not\perp_{\g} Z | \mathbf{Z'}_{-z} \cup \mathbf{X}$.
	\end{enumerate} 
 	From this it follows that the output of Algorithm \ref{Algorithm} is $\mathbf{Z_2} \subseteq \mathbf{Z}$, independently of the order at which the nodes in $\mathbf{Z}$ are considered. By the same d-separation result and Lemma \ref{lemma: d-sep VAS}, both $\mathbf{Z_2}$ and any possible intermediate $\mathbf{Z'_2} \supseteq \mathbf{Z_2}$ are valid adjustment sets relative to $(\mathbf{X},\mathbf{Y})$ in $\g$.

	If $\mathbf{V}$ follows a causal linear model compatible with $\g$ we can apply Theorem \ref{cor:bignew12} to conclude that 
	$a.var(\tau^{\mathbf{z_2}}_{\mathbf{xy}}) \leq a.var(\tau^{\mathbf{z}}_{\mathbf{xy}})$. Further, by the minimality of $\mathbf{Z_2}$, no other subset of $\mathbf{Z}$ is guaranteed to have a better asymptotic variance than $\mathbf{Z_2}$ by Theorem \ref{cor:bignew12}.
\end{proofof}

The following result is very similar to Proposition 2 from \citet{vanderweele2011new}, albeit in the graphical rather than the potential outcomes framework.

\begin{Lemma}
	\label{lemma: d-sep VAS}
	Consider pairwise disjoint node sets $\mathbf{X},\mathbf{Y},\mathbf{Z_1}$ and $\mathbf{Z_2}$ in a causal \mpdag{} $\g$, such that $\mathbf{Z} = \mathbf{Z_1} \cup \mathbf{Z_2}$ is a valid adjustment set relative to $(\mathbf{X},\mathbf{Y})$ in $\g$. If $\mathbf{Y} \perp_{\g} \mathbf{Z_1} | \mathbf{Z_2} \cup \mathbf{X}$, then $\mathbf{Z_2}$ is a valid adjustment set relative to $(\mathbf{X},\mathbf{Y})$.
\end{Lemma}

\begin{proof}

	We prove the statement by contraposition, so assume that $\mathbf{Z}$ is a valid adjustment set relative to $(\mathbf{X}, \mathbf{Y})$ in $\g$, whereas $\mathbf{Z_2}$ is not. Since $\mathbf{Z_2} \subseteq \mathbf{Z}$, $\mathbf{Z_2} \cap \fb{\g} = \emptyset$ trivially holds. Thus, there must be a proper, non-causal, definite status path $p$ from $\mathbf{X}$ to $\mathbf{Y}$ that is open given $\mathbf{Z_2}$, while being blocked by $\mathbf{Z}$. Therefore, $p$ contains at least one non-collider $N \in (\mathbf{Z} \setminus \mathbf{Z_2}) = \mathbf{Z_1}$, where we choose the $N$ closest to the endpoint node $Y \in \mathbf{Y}$. Then the subpath path $p(N,Y)$ of $p$ is a definite status path from $\mathbf{Z_1}$ to $\mathbf{Y}$ that is open given $\mathbf{Z_2} $. As $p$ is proper, $p(N,Y)$ does not contain any nodes in $\mathbf{X}$ and is therefore also open given $\mathbf{Z_2} \cup \mathbf{X}$. Thus, $\mathbf{Y} \not\perp_{\g} \mathbf{Z_1} | \mathbf{Z_2} \cup \mathbf{X}$. 
\end{proof}

%merge with Lemma from section C and say it is just weak union property and put some Pearl paper citation in \citet{verma1990causal}.

%\begin{proof}
%	We prove the statement by contraposition, so assume that there exists a proper path $p$ from $\mathbf{Z'_1}$ to $\mathbf{Y}$ that is open given $\mathbf{Z''_1} \cup \mathbf{Z_2} \cup \mathbf{X}$. By assumption the path $p$ is closed given $\mathbf{Z_2} \cup \mathbf{X}$ and hence $p$ must contain a collider $C$, such that $\de(C,\g) \cap \mathbf{Z''_1} = \emptyset$, but this implies the existence of a path from $\mathbf{Z''_1} \subset \mathbf{Z_1}$ to $\mathbf{Y}$ that is open given $\mathbf{Z_2} \cup \mathbf{X}$, so $\mathbf{Y} \perp_{\g} \mathbf{Z_1} | \mathbf{Z_2} \cup \mathbf{X}$ does not hold.
%\end{proof}

\begin{Lemma}
	\label{lemma max size unique}
	Consider pairwise disjoint node sets $\mathbf{X},\mathbf{Y}$ and $\mathbf{Z}$ in a causal \mpdag{} $\g$, such that 
%	$\mathbf{Z} = \mathbf{Z_1} \cup \mathbf{Z_2}$ is a valid adjustment set relative to $(\mathbf{X},\mathbf{Y})$ in $\g$ and 
$\mathbf{Y} \perp_{\g} \mathbf{Z_1} | \mathbf{Z_2} \cup \mathbf{X}$. 
	%and $\mathbf{Z_1} \subseteq \mathbf{Z}$ is of maximal size. 
	Then given any other partition $\mathbf{Z} = \mathbf{Z'_1} \cup \mathbf{Z'_2}$, such that $\mathbf{Y} \perp_{\g} \mathbf{Z'_1} | \mathbf{Z'_2} \cup \mathbf{X}$, it also holds that $\mathbf{Y} \perp_{\g} \mathbf{Z_1} \cup \mathbf{Z'_1} | (\mathbf{Z_2} \cap \mathbf{Z'_2}) \cup \mathbf{X}$.
\end{Lemma}

\begin{proof}
% 	We prove the statement by contradiction, so assume there exists a set $\mathbf{\hat{Z}} \subseteq \mathbf{Z}$, such that $\mathbf{\tilde{Z}} = \mathbf{\hat{Z}} \setminus \mathbf{Z_1} = \mathbf{\hat{Z}} \cap \mathbf{Z_2} \neq \emptyset$. 
%	%By Lemma \ref{lemma: order independence} with $\mathbf{Z'_1}=Z$ and $\mathbf{Z''_1}=\emptyset$, $\mathbf{Y} \perp_{\g} Z | \mathbf{Z_2} \cup \mathbf{X}$.  
%	We now show that 
%	\[
%	\mathbf{Y} \perp_{\g} \mathbf{Z_1} \cup \mathbf{\tilde{Z}} | (\mathbf{Z_2} \setminus \mathbf{\tilde{Z}}) \cup \mathbf{X}.
%	\]
	 We first show that $\mathbf{Y} \perp_{\g} \mathbf{Z_2} \cap \mathbf{Z'_1}| (\mathbf{Z_2} \cap \mathbf{Z'_2}) \cup \mathbf{X}$ by contradiction. Assume that $\mathbf{Y} \perp_{\g} \mathbf{Z_1} | \mathbf{Z_2}  \cup \mathbf{X}$ and $\mathbf{Y} \perp_{\g} \mathbf{Z'_1} | \mathbf{Z'_2}  \cup \mathbf{X}$. Let $p$ be a proper, definite status path from $\mathbf{Z_2} \cap \mathbf{Z'_1} \subseteq \mathbf{Z'_1}$ to some node $Y \in \mathbf{Y}$ and assume that $p$ is open given $(\mathbf{Z_2} \cap \mathbf{Z'_2}) \cup \mathbf{X}$. By assumption $p$ is blocked by $\mathbf{Z'_2}  \cup \mathbf{X}$. As $((\mathbf{Z_2} \cap \mathbf{Z'_2}) \cup \mathbf{X}) \subseteq (\mathbf{Z'_2}\cup \mathbf{X})$, this can only be the case if there exists a non-collider $N \in (\mathbf{Z'_2} \setminus \mathbf{Z_2}) \subseteq \mathbf{Z_1}$ on $p$. As a subpath of $p$, $p(N,Y)$ is trivially open given $(\mathbf{Z_2} \cap \mathbf{Z'_2}) \cup \mathbf{X}$ and since $p$ is proper, $p(N,Y)$ contains no node in $\mathbf{Z_2} \cap \mathbf{Z'_1}$. Therefore, $p(N,Y)$ is also open given $(\mathbf{Z_2} \cap \mathbf{Z'_1}) \cup (\mathbf{Z_2} \cap \mathbf{Z'_2}) \cup \mathbf{X}= \mathbf{Z_2} \cup \mathbf{X}$. But this contradicts that $\mathbf{Y} \perp_{\g} \mathbf{Z_1} | \mathbf{Z_2} \cup \mathbf{X}$. Thus, any such $p$ must be blocked by $(\mathbf{Z_2} \cap \mathbf{Z'_2}) \cup \mathbf{X}$.
	  
	  By Lemma \ref{lemma: contraction}, the two d-separation statements 
	  \begin{enumerate}
	  	\item $\mathbf{Y} \perp_{\g} \mathbf{Z_2} \cap \mathbf{Z'_1}| (\mathbf{Z_2} \cap \mathbf{Z'_2}) \cup \mathbf{X}$ and
	  	\item $\mathbf{Y} \perp_{\g} \mathbf{Z_1} | \mathbf{Z_2} \cup \mathbf{X}$
	  \end{enumerate}  jointly imply that 
	  $\mathbf{Y} \perp_{\g} \mathbf{Z_1} \cup \mathbf{Z'_1} | (\mathbf{Z_2} \cap \mathbf{Z'_2}) \cup \mathbf{X}$.
	 
	 % and therefore $\mathbf{Y} \perp_{\g} \mathbf{Z_1} \cup \mathbf{\hat{Z}} | \mathbf{Z_2} \setminus \mathbf{\hat{Z}} \cup \mathbf{X}$.

\end{proof}

The two following Lemmas are general properties of the d-separation criterion from \citet{verma1990causal}. They extend to the causal \mpdag{} setting with the fact that by Lemma \ref{lemma:dsepp1}, d-separation in a \mpdag{} implies d-separation in every represented DAG and vice versa. 

%and that similarly, any valid adjustment set in a \mpdag{} is a valid adjustment set in every represented DAG, including the true underlying one.

\begin{Lemma} (Weak union and decomposition property)
	\label{lemma: weak union}
	Let $\mathbf{X},\mathbf{Y}$ and $\mathbf{Z}$ be pairwise disjoint node sets in a DAG $\g$. If $\mathbf{Y} \perp_{\g} \mathbf{Z} | \mathbf{X}$, then for any $\mathbf{Z'} \subseteq \mathbf{Z}$ and $\mathbf{Z''} \subset \mathbf{Z}$ such that $\mathbf{Z'} \cap \mathbf{Z''} = \emptyset$, $\mathbf{Y} \perp_{G} \mathbf{Z''}|\mathbf{Z'} \cup \mathbf{X}$. 
\end{Lemma}

\begin{Lemma} (Contraction property)
	\label{lemma: contraction}
	Let $\mathbf{X},\mathbf{Y},\mathbf{Z}$ and $\mathbf{W}$ be pairwise disjoint node sets in a DAG $\g$. If $\mathbf{X} \perp_{\g} \mathbf{Y} | \mathbf{Z}$ and $\mathbf{X} \perp_{\g} \mathbf{W} | \mathbf{Z} \cup \mathbf{Y}$ then $\mathbf{X} \perp_{\g} \mathbf{Y} \cup \mathbf{W} | \mathbf{Z}$.
\end{Lemma}

\section{Proof of Theorem~\ref{thm:optimalSetCPDAG}}
\label{section-E}

\subsection{Theorem~\ref{thm:optimalSetCPDAG} for DAGs}

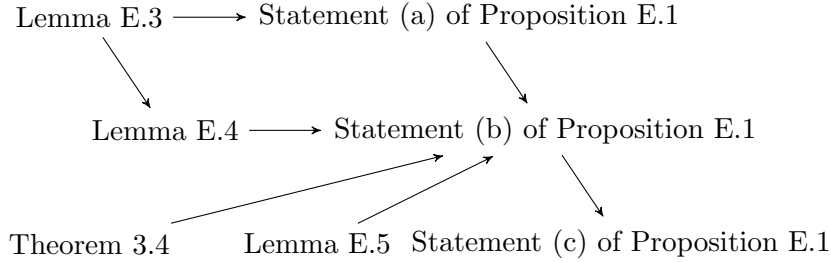
\begin{figure}[ht]
	\centering
	\begin{tikzpicture}[>=stealth',shorten >=1pt,node distance=3cm, main node/.style={minimum size=0.4cm}]
	
	\node[main node,yshift=0cm,xshift=0cm]         (C2) {Statement \ref{optimal} of Proposition~\ref{almost}};	
	\node[main node,yshift=1.5cm,xshift=-1cm]  (C1)  at (C2)  {Statement \ref{VAS} of Proposition~\ref{almost}};	
	\node[main node,yshift=-1.5cm,xshift=1cm] (C3) at (C2) {Statement \ref{smallest} of Proposition~\ref{almost}};	
		
%	\node[main node,yshift=0cm,xshift=1cm]         (R1) [right of= C2] {\textbf{Theorem~\ref{almost}}};
	
    \node[main node,xshift=-2cm]            (L2) [left of= C2]   {Lemma \ref{y-separation}}; 
	\node[main node,yshift=1.5cm,xshift=-1cm]            (L1) at (L2)   {Lemma \ref{lemma:eadd}};    
	\node[main node,yshift=-1.5cm,xshift=2cm]            (L3)  at (L2)   {Lemma \ref{x-separation}};   
	\node[main node,xshift=-3cm]            (L4) at (L3)   {Theorem \ref{cor:bignew12}};

	\draw[->] (L1) edge    (C1);
    
    \draw[->] (L1) edge    (L2);
    
%    \draw[->] (L2) edge    (C2)
%    \draw[->] (L3) edge    (C2);
%    \draw[->] (C1) edge    (C2);

	\draw[->] (L2) edge    (C2);
	\draw[->] (L3) edge    (C2);
	\draw[->] (L4) edge    (C2);
	\draw[->] (C1) edge    (C2);

	\draw[->] (C2) edge    (C3);

%	\draw[->] (C1) edge    (R1);
%	\draw[->] (C2) edge    (R1);
%	\draw[->] (C3) edge    (R1);
	
	\end{tikzpicture}
	\caption{Proof structure of Proposition \ref{almost}.}
	\label{figproofDAG}
\end{figure}

Figure~\ref{figproofDAG} shows the structure of the proof of Proposition \ref{almost}, which is Theorem~\ref{thm:optimalSetCPDAG} for DAGs only.

\begin{Proposition}
	Let $\mathbf{X}$ and $\mathbf{Y}$ be disjoint node sets in a causal DAG $\g=(\mathbf{V},\mathbf{E})$, such that $\mathbf{Y} \subseteq \de(\mathbf{X},\g)$. Let the density $f$ of $\mathbf{V}$ be compatible with $\g$ and let $\mathbf{O}=\optb{\g}$. Then the following three statements hold:
	\begin{enumerate}[label = (\alph*)]
		\item The set $\mathbf{O}$ is a valid adjustment set relative to $(\mathbf{X,Y})$ in $\g$ if and only if there exists a valid adjustment set relative to $(\mathbf{X,Y})$ in $\g$. 
		\label{VAS}
		\item Let $\mathbf{Z}$ be a valid adjustment set relative to $(\mathbf{X},\mathbf{Y})$ in $\g$. If $\mathbf{V}$ follows a causal linear model compatible with $\g$ then
		\[
		a.var(\hat{\tau}^{\mathbf{o}}_{\mathbf{yx}}) \leq a.var(\hat{\tau}^{\mathbf{z}}_{\mathbf{yx}}).
		\]
		\label{optimal} 
		\item Let $\mathbf{Z}$ be a valid adjustment set relative to $(\mathbf{X},\mathbf{Y})$ in $\g$, such that 
		\[
		a.var(\hat{\tau}^{\mathbf{o}}_{\mathbf{yx}}) = a.var(\hat{\tau}^{\mathbf{z}}_{\mathbf{yx}}).
		\]
		If $\mathbf{V}$ follows a causal linear model compatible with $\g$ and $f$ is faithful to $\g$ then $\mathbf{O} \subseteq \mathbf{Z}$. 
		\label{smallest}
	\end{enumerate}
	
	\label{almost}
\end{Proposition}

%We first introduce preparatory Lemma~\ref{loops} and then prove the three statements of Theorem \ref{A-almost}.
%
%\begin{Lemma}
%	Let $\mathbf{X}$ and $\mathbf{Y}$ be disjoint node sets in a DAG $\g$, such that $\mathbf{Y} \subseteq \de(\mathbf{X},\g)$. There exists a valid adjustment set relative to $(\mathbf{X},\mathbf{Y})$ in $\g$ if and only if $\mathbf{X} \cap \de(\cnb{G},\g) = \emptyset$.
%	\label{loops}
%\end{Lemma}\
%
%\begin{proof}
%	We first prove that $\mathbf{X} \cap \de(\cnb{\g},\g) \neq \emptyset$ implies that there cannot a exist a valid adjustment set. Let $X \in \mathbf{X} \cap \de(\cnb{\g},\g)$. Then there exists a path $p$, of the form
%	\[
%	X \leftarrow \dots \leftarrow C \rightarrow \dots \rightarrow Y,
%	\]
%	where $C \in \cnb{\g}$ and $Y \in \mathbf{Y} \cap \cnb{\g}$. Since $p$ is a non-collider, non-causal path from $\mathbf{X}$ to $\mathbf{Y}$ that consists in its entirety of forbidden nodes, no set of covariates may fulfill the blocking condition of Definition \ref{A-adjustment}.
%	
%	We now prove converse statement.
%	As we assume $\mathbf{Y} \subseteq \de(\mathbf{X},\g)$, $\mathbf{Y} \subseteq \cnb{G}$. The condition $\mathbf{X} \cap \de(\cnb{\g},\g) = \emptyset$ implies that there are no proper causal paths from any $C \in \cnb{\g}$ to any $X_i \in \mathbf{X}$. Thus paths of the form 
%	\[
%	X_i \leftarrow \dots \leftarrow C \rightarrow \dots \rightarrow Y_j,
%	\]
%	with $C \in \cnb{\g}$ may not exist. By Theorem 5.9 from \citet{perkovic16} it follows that a valid adjustment set exists.
%\end{proof}
%

We prove each of the three Statements of Proposition \ref{almost} separately, due to the complexity of the proofs.

\begin{proofof}[Statement \ref{VAS} of Proposition~\ref{almost}]
	As we are considering a DAG, the amenability condition is trivially fulfilled.
	By construction $\mathbf{O}$ fulfills the forbidden set condition relative to $(\mathbf{X,Y})$ in $\g$ in Definition \ref{adjustment}. By the assumption that there exists a valid adjustment set and Corollary \ref{cor:noforbx} it follows that $\mathbf{X} \cap \cnb{\g}  = \emptyset$. Hence, we can invoke Lemma~\ref{lemma:eadd} with $\mathbf{T} = \mathbf{X}$, to conclude that $\mathbf{O}$ fulfills the blocking condition relative to $(\mathbf{X,Y})$ in $\g$ from Definition \ref{adjustment}. 
	
\end{proofof}

\begin{proofof}[Statement \ref{optimal} of Proposition~\ref{almost}]
	Let $\mathbf{W} = \mathbf{O} \cap \mathbf{Z}$, $\mathbf{T} = \mathbf{Z} \setminus \mathbf{O}$ and $\mathbf{S} = \mathbf{O} \setminus \mathbf{Z}$. 
	The sets $\mathbf{T}$ and $\mathbf{W} \cup \mathbf{S} = \mathbf{O}$ satisfy Lemma \ref{y-separation} and hence $\mathbf{Y} \perp_{\g} \mathbf{T}|\mathbf{O} \cup \mathbf{X}$. 
	Additionally, $\mathbf{S} \subseteq \mathbf{O}$ and $\mathbf{W} \cup \mathbf{T} = \mathbf{Z}$ satisfy the conditions of Lemma \ref{x-separation} and hence $\mathbf{X} \perp_{\g} \mathbf{S}|\mathbf{Z}$. By Theorem \ref{cor:bignew12} it thus follows that
	\[
	\begin{aligned} 
	a.var(\hat{\tau}^\mathbf{o}_{\mathbf{yx}})
	& 
	= a.var(\hat{\beta}_{\mathbf{yx}.\mathbf{o}}) 
	= a.var(\hat{\beta}_{\mathbf{yx}.\mathbf{ws}}) \\ 
	&
	\leq a.var(\hat{\beta}_{\mathbf{yx}.\mathbf{wt}})  = a.var(\hat{\beta}_{\mathbf{yx}.\mathbf{z}}) =
	a.var(\hat{\tau}^\mathbf{z}_{\mathbf{yx}}).
	\end{aligned}
	\]
	%Is this formulation OK?
\end{proofof}

\begin{proofof}[of Statement \ref{smallest} of Proposition \ref{almost}]
	We prove the claim by contradiction. Let $\mathbf{Z}$ be an asymptotically optimal valid adjustment set that is not a superset of $\mathbf{O}$. Let $\mathbf{W}=\mathbf{O} \cap \mathbf{Z}$ and $\mathbf{S} =\mathbf{O} \setminus \mathbf{Z} \neq \emptyset$. 
	
	We now show that $\mathbf{S}$ and $\mathbf{Y}$ are not d-separated by $\mathbf{W} \cup \mathbf{X}$. Let $S \in\mathbf{S} \subseteq \mathbf{O}$. Then there exists a directed  path $p$ from $S$ to $\mathbf{Y}$ that consists of causal nodes, except for $S$. Hence, $p$ cannot be blocked by $\mathbf{X} \cup \mathbf{W}$, as $(\mathbf{X} \cup \mathbf{W}) \cap (\cnb{G} \cup \mathbf{S}) = \emptyset$. As the underlying distribution $f$ is assumed to be faithful to $\g$ it thus follows that $\boldsymbol{\beta}_{\mathbf{ys}.\mathbf{xw}} \neq 0$.
	
	Within the proof of Statement \ref{optimal} of Proposition \ref{almost} we have already shown that $\sigma_{x_ix_i.\mathbf{x_{-i}z}} \leq \sigma_{x_ix_i.\mathbf{x_{-i}o}}$ for all $X_i \in \mathbf{X}$, for any valid adjustment set. Therefore, for $\mathbf{Z}$ to be asymptotically optimal, $\sigma_{y_jy_j.\mathbf{xz}} \leq \sigma_{y_jy_j.\mathbf{xo}}$ has to hold for all $Y_j \in \mathbf{Y}$. 
	By equation \eqref{y_variance_calc}
	\[
	\sigma_{y_jy_j.\mathbf{xz}} - \sigma_{y_jy_j.\mathbf{xo}} = \boldsymbol{\beta}_{y_j\mathbf{s}.\mathbf{xw}} \Sigma_{\mathbf{ss}.\mathbf{xwt}} \boldsymbol{\beta}^T_{y_j\mathbf{s}.\mathbf{xw}},
	\]
 for all $Y_j \in \mathbf{Y}$, 	where $\mathbf{T} = \mathbf{Z} \setminus \mathbf{O} \neq \emptyset$. The case $\mathbf{T}= \emptyset$ is the equivalent statement, with the convention that empty conditioning sets are omitted and  follows directly from Lemma \ref{convarsets} jointly with the fact that in this case $\mathbf{O}=\mathbf{Z} \cup \mathbf{S}$. As the conditional distributions are assumed to be non-degenerate, 
	$\Sigma_{\mathbf{ss}.\mathbf{xwt}}$ is positive definite. We have already shown that $\boldsymbol{\beta}_{\mathbf{ys}.\mathbf{xw}} \neq 0$, so it follows that $\sigma_{y_jy_j.\mathbf{xz}} > \sigma_{y_jy_j.\mathbf{xo}}$ for some $Y_j \in \mathbf{Y}$, which yields a contradiction. 
\end{proofof}

While Lemma \ref{lemma:eadd} is a rather technical result, Lemma \ref{y-separation} and Lemma \ref{x-separation} have an intuitive interpretation. Lemma \ref{y-separation} states that given a valid adjustment set $\mathbf{Z}$, this set may not contain more information on $\mathbf{Y}$ conditionally on $\mathbf{X}$ than $\mathbf{O}$. Lemma \ref{x-separation} states that $\mathbf{Z}$ cannot contain less information on $\mathbf{X}$ than $\mathbf{O}$. The surprising fact that $\mathbf{O}$ has both of these properties simultaneously is the key to its asymptotic optimality.

We now introduce an object that is used in the proof of Lemma \ref{lemma:eadd}.

\begin{Definition}
	Consider two disjoint node sets $\mathbf{X}$ and $\mathbf{Y}$ in a causal DAG $\g$, such that $\mathbf{Y} \subseteq \de(\mathbf{X},\g)$ and let $p$  be a path to some $Y\in\mathbf{Y}$. Then the \emph{maximal causal segment} of $p$ with respect to $(\mathbf{X},\mathbf{Y})$, is the longest subpath $p(C,Y)$ of $p$, such that all nodes on $p(C,Y)$ are in $\cnb{\g}$ and $p(C,Y)$ is directed towards $Y$. 
\end{Definition} 

For an example consider Figure \ref{AllinOne} and consider the path $p=(X,A_1,A_2,Y,F)$. The longest maximal causal segment with respect to $(X,F)$ of $p$ is $Y \rightarrow F$. It is the longest directed subpath of $p$ that consists of causal nodes and ends in $F$. 

\begin{Lemma}
	Let $\mathbf{X}$ and $\mathbf{Y}$ be disjoint node sets in a causal DAG $\g$, such that $\mathbf{Y} \subseteq \de(\mathbf{X},\g)$ and let $\mathbf{O}=\optb{\g}$. Let $\mathbf{T}$ be a node set such that $\mathbf{T} \cap \de(\cnb{\g},\g) = \emptyset$ and $\mathbf{T} \cap \mathbf{O} = \emptyset$. If an adjustment set relative to $(\mathbf{X,Y})$ in $\g$ exists, then all proper non-causal paths from $\mathbf{T}$ to $\mathbf{Y}$ in $\g$ that contain no nodes from $\mathbf{X} \setminus \mathbf{T}$ are blocked by $\mathbf{O} \cup (\mathbf{X} \setminus \mathbf{T})$.  

	\label{lemma:eadd}
\end{Lemma}

\begin{proof}
	Let $p$ be a proper non-causal path from $T \in \mathbf{T}$ to $Y \in \mathbf{Y}$ that is open given $\mathbf{O} \cup (\mathbf{X} \setminus \mathbf{T})$ and contains no nodes in $\mathbf{X} \setminus \mathbf{T}$. 
%	Without loss of generality we can assume that $p$ contains only one node in $\mathbf{Y}$ (otherwise consider an appropriate subpath).
%	It is enough to show that $\mathbf{O}$ blocks all non-causal paths from $\mathbf{T}$ to $\mathbf{Y}$ that contain no nodes from $\mathbf{X} \setminus \mathbf{T}$ and only one node in $\mathbf{Y}$.
%	Let $p$ be one such path from $T \in \mathbf{T}$ to $Y_j \in \mathbf{Y}$ in $\g$. 
	Let $c_{Y} = p(C_1,Y)$ be the maximal causal segment of $p$ with respect to $(\mathbf{X,Y})$, where we use that $\mathbf{Y} \subseteq \de(\mathbf{X},\g)$ implies that $\mathbf{Y} \subseteq \cnb{\g}$.
%	, with possibly $C_1 = Y$. 
Then $p$ is of the form
	\begin{enumerate}[label=(\alph*)]
	\item
	$T \cdots V \rightarrow C_1 \rightarrow \dots \rightarrow Y$ \ \textit{or} \label{rightR1}
	\item
	$T \cdots V \leftarrow C_1 \rightarrow \dots \rightarrow Y$. \label{leftR1}
	\end{enumerate}
	
	If $p$ is of the form \ref{rightR1} and $V \in \mathbf{O}$, then $p$ is blocked by $\mathbf{O}\cup (\mathbf{X} \setminus \mathbf{T})$. We now show that $V \in \mathbf{O}$ does in fact hold by contradiction, so suppose that $V \notin \mathbf{O}$. Note that $V \neq T$, as otherwise $p$ would be a causal path from $T$ to $Y$. By assumption $p$ is proper and contains no nodes from $\mathbf{X} \setminus \mathbf{T}$ and hence $V \notin \mathbf{X}$. As $V$ is a parent of the causal node $C_1$, $V \notin \mathbf{O}$ can only hold if $\mathbf{V} \in \fb{\g}$. As $ V \notin \mathbf{X}$ it follows that $V \in \de(\cnb{\g},\g)$. Then there is a proper directed path from $\mathbf{X}$ to $V$. Additionally, $p(V,Y)$ is a directed path towards $Y$ that does not contain a node in $\mathbf{X}$, so $V \in \cnb{\g}$, which contradicts that $c_{Y}$ is the maximal causal segment of $p$ with respect to $(\mathbf{X,Y})$.
	
	If $p$ is of the form \ref{leftR1}, $V \in \de(\cnb{\g},\g)$ which implies $V \neq T$. Further, $p$ has to contain at least one collider, as otherwise $T \in \de(V,\g) \subseteq \de(\cnb{\g},\g)$. Let $V'$ be the collider on $p$ that is closest to $V$.  Then $V' \in \de(V,\g)$, so $V' \in \fb{\g}$. Therefore, $\de(V',\g) \cap \mathbf{O} = \emptyset$. 
	By the assumption that there is an adjustment set relative to $(\mathbf{X,Y})$ in $\g$ and Corollary~\ref{cor:noforbx} it follows that $\mathbf{X} \cap \de(\cn{\g},\g) = \emptyset$. Therefore, it also holds that $\de(V',\g) \cap (\mathbf{X} \setminus \mathbf{T}) = \emptyset$. 
%	Additionally, since $\de(V'_1,\g) \subseteq \de(\cn{\g},\g)$ it also holds that $\de(V'_1,\g) \cap (\mathbf{X} \setminus \mathbf{T}) = \emptyset$. Otherwise $\mathbf{X} \cap \de(\cn{\g},\g) \neq \emptyset$ which by Corollary~\ref{cor:noforbx} contradicts that there is an adjustment set relative to $(\mathbf{X,Y})$ in $\g$. 
	Since $V'$ is a collider on $p$ and $\de(V',\g) \cap ( \mathbf{O} \cup (\mathbf{X} \setminus \mathbf{T})) = \emptyset$, $p$ is blocked by $\mathbf{O} \cup (\mathbf{X} \setminus \mathbf{T})$.
\end{proof}

\begin{Lemma}
	Let $\mathbf{X}$ and $\mathbf{Y}$ be disjoint node sets in a causal DAG $\g$, such that $\mathbf{Y} \subseteq \de(\mathbf{X},\g)$ and let $\mathbf{O}=\optb{\g}$. Let $\mathbf{T}$ be a set such that 
	$\mathbf{T} \cap (\fb{\g} \cup \mathbf{O}) = \emptyset$. If there exists a valid adjustment set relative to $(\mathbf{X},\mathbf{Y})$ in $\g$, then $\mathbf{Y} \perp_{\g}\mathbf{T}|\mathbf{O} \cup \mathbf{X}$.  
	\label{y-separation}
\end{Lemma}

\begin{proof}%[Lemma~\ref{y-separation}]
	It is enough to show that all paths from $\mathbf{T}$ to $\mathbf{Y}$ that are proper, are blocked by $\mathbf{O} \cup \mathbf{X}$. Let $p$ be such a path from from $T \in \mathbf{T}$ to $Y \in \mathbf{Y}$. % We will consider the cases, that $p$ contains nodes from $\mathbf{X}$ and that it does not separately.
	
	First, suppose that no node from $\mathbf{X}$ lies on $p$. If $p$ is a non-causal path from $T$ to $Y$, then by Lemma~\ref{lemma:eadd}, $p$ is blocked by  $\mathbf{O} \cup \mathbf{X}$. If $p$ is causal from $T$ to $Y$,  then by the fact that $\an(\mathbf{Y},\g) \cap \fb{\g} = \cnb{\g}$, the non-forbidden node $O$ closest to $Y$ on $p$ is in $\mathbf{O}$. Since $T \notin (\fb{\g} \cup \mathbf{O})$ such a node $O$ exists and it holds that $O \neq T$. Clearly, $O \neq Y$ and therefore, $O$ is a non-collider on $p$. Hence $p$ is blocked given $\mathbf{O} \cup \mathbf{X}$.
	
	Now, assume that $p$ contains at least one node from $\mathbf{X}$. If a node in $\mathbf{X}$ is a non-collider on $p$, $p$ is blocked by $\mathbf{O} \cup \mathbf{X}$. So we can assume that all nodes on $p$ that are in $\mathbf{X}$ are colliders. Let $X \in \mathbf{X}$ be the collider on $p$ that is closest to $Y$. Then $p(X,Y)$ is a proper non-causal
	path from $\mathbf{X}$ to $\mathbf{Y}$ and since, by the already proven Statement \ref{VAS} in Theorem~\ref{almost}, $\mathbf{O}$ is a valid adjustment set relative to 
	$(\mathbf{X,Y})$ in $\g$, $\mathbf{O}$ blocks $p(X,Y)$. Now assume that $p(X,Y)$ is open given $\mathbf{O} \cup \mathbf{X}$ while being blocked by $\mathbf{O}$. By Lemma \ref{lemma:collider} this contradicts that $\mathbf{O}$ is a valid adjustment. Hence, we can conclude that $p$ is blocked by $\mathbf{O} \cup \mathbf{X}$.
	
\end{proof}

\begin{Lemma}
	Let $\mathbf{X},\mathbf{Y},\mathbf{S}$ and $\mathbf{W}$ be pairwise disjoint node sets in a causal DAG $\g$, such that $\mathbf{Y} \subseteq \de(\mathbf{X},\g)$, let $\mathbf{O}=\optb{\g}$ and assume that $\mathbf{S}  \subseteq \mathbf{O}$. If $\mathbf{W}$ is a valid adjustment set relative to $(\mathbf{X},\mathbf{Y})$ in $\g$, then $\mathbf{X} \perp_{\g} \mathbf{S}| \mathbf{W}$. 
	\label{x-separation}
\end{Lemma}

\begin{proof}%[Lemma~\ref{A-x-separation}]
	For contraposition, suppose that $\mathbf{X} \not\perp_{\g}\mathbf{S}|\mathbf{W}$ and that there exists a valid adjustment set relative to $(\mathbf{X},\mathbf{Y})$ in $\g$. We will show that this implies the existence of a proper non-causal path from $\mathbf{X}$ to $\mathbf{Y}$ that is open given $\mathbf{W}$ and hence $\mathbf{W}$ is not a valid adjustment set. %This shows that $\mathbf{W}$ is not a valid adjustment set. 
	
%	If there is a path from $\mathbf{X}$ to $\mathbf{S}$ that is open given $\mathbf{W}$, then there is also at least one such path that is proper.
	Let $p$ be a proper path from $X \in \mathbf{X}$ to $S \in \mathbf{S}$ that is open given $\mathbf{W}$. Since $S \in \mathbf{S} \subseteq \mathbf{O}$, there exists a directed path $p'$ from $S$ to some $Y \in \mathbf{Y}$ that consists of $S$ and nodes in $\cnb{\g}$.
	As $\mathbf{W} \cap \cnb{\g} = \emptyset$, $p'$ must be open given $\mathbf{W}$. 
	
	Let $I$ be the the node closest to  $X$ on $p$ that is also on $p'$ and consider the path $q=p(X, I) \oplus p'(I, Y)$. Since $I$ is on $p'$ and $(\mathbf{S} \cup \cnb{\g}) \cap \mathbf{X}=\emptyset$, $I \neq X$. 
	Suppose now that $I=Y$. Then $q$ is a subpath of $p$ and since $p$ is open given $\mathbf{W}$, $q$ is a path from $\mathbf{X}$ to $\mathbf{Y}$ that is open given $\mathbf{W}$. Suppose, now that $I \neq Y$. As both $p(X,I)$ and $p'(I,Y)$ are open given $\mathbf{W}$ and since $p'(X,I)$ is directed towards $Y$, $I$ is a  non-collider on $q$. With the fact that $I$ is on $p'$ and $p'$ contains no node in $\mathbf{W}$ we can thus conclude that $q$ is a path from $\mathbf{X}$ to $\mathbf{Y}$ that is open given $\mathbf{W}$.
	
	We now show that $q$ is proper.
	By the assumption that there is an adjustment set relative to $(\mathbf{X,Y})$ in $\g$ and Corollary~\ref{cor:noforbx} it follows that $\mathbf{X} \cap \de(\cn{\g},\g) = \emptyset$. Hence, $p'$ does not contain a node in $\mathbf{X}$ and as $p$ is proper, it follows that $q$ is a proper path from $\mathbf{X}$ to $\mathbf{Y}$.
	
	It is left to show that $q$ is a non-causal path. For contradiction, suppose that $q$ is a causal path. Then $p$ must be of the form 
	$
	X \rightarrow \dots \rightarrow I \cdots S.
	$
	Since $S \in \mathbf{O} \subseteq \an(\mathbf{Y},\g) \setminus \fb{\g}$ it follows that $S \notin \de(\mathbf{X},\g)$ and hence, $S \neq I$. Thus, a descendant of $I$ must be a collider on $p$. Since $S \neq I$ and all nodes expect for $S$ on $q$ are in $\cnb{\g}$, it follows that $I \in \cnb{\g}$. Thus, there exists a collider on $p$ that is in $\de(\cnb{\g},\g)$. Since $\mathbf{W} \cap \de(\cnb{\g},\g) = \emptyset$, this contradicts our assumption that $p$ open given $\mathbf{W}$. 
\end{proof}

\subsection{Extension of Theorem~\ref{thm:optimalSetCPDAG} to \mpdag{}s} 

%Figure~\ref{figproof1} shows the structure of the proof of Theorem~\ref{A-thm:optimalSetCPDAG} for \mpdag{}s.

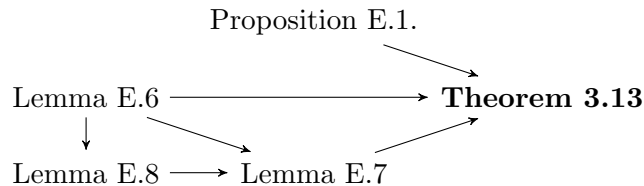
\begin{figure}[ht]
	\centering
	\begin{tikzpicture}[>=stealth',shorten >=1pt,node distance=3cm, main node/.style={minimum size=0.4cm}]
	
	\node[main node,yshift=0cm,xshift=0cm]         (T65) {\textbf{Theorem~\ref{thm:optimalSetCPDAG}}};
	\node[main node,yshift=1cm,xshift=-3cm]  (T313)  at (T65)  {Proposition \ref{almost}.};
	
	\node[main node,yshift=-1cm,xshift=-3 cm] (L66) at (T65) {Lemma \ref{lemma:optimalSetEquiv}};
	\node[main node,xshift=-3cm]            (L68) [left of= T65]   {Lemma \ref{lemma:additional-1}};    
	\node[main node,yshift=-1cm,xshift=0cm]            (L67)  at (L68)   {Lemma \ref{lemma:forbIsSame}};

	%   \draw[->] (T37) edge    (T65);
	\draw[->] (T313) edge    (T65);
	\draw[->] (L68) edge    (T65);
	\draw[->] (L66) edge    (T65);
	\draw[->] (L67) edge    (L66);
	\draw[->] (L68) edge    (L67);
	\draw[->] (L68) edge    (L66);
	
	\end{tikzpicture}
	\caption{Proof structure of Theorem \ref{thm:optimalSetCPDAG} for CPDAGs and maximal PDAGs.}
	\label{figproof1}
\end{figure}

\begin{proofof}[Theorem~\ref{thm:optimalSetCPDAG} for maximal PDAGs]
	We first prove Statement \ref{res1-cpdag}. Let $\g[D]$ be a DAG in $[\g]$. By Lemma \ref{lemma:additional-1}, $\mathbf{Y} \subseteq \possde(\mathbf{X},\g)$ also implies that $\mathbf{Y} \subseteq \de(\mathbf{X},\g)$. Hence, for any DAG $\g[D]$ in $[\g]$, $\mathbf{Y} \subseteq \de(\mathbf{X},\g[D])$. 
	By Lemma~\ref{lemma:optimalSetEquiv}, $\optb{\g[D]} = \optb{\g}$ and thus $\optb{\g}$ satisfies the generalized adjustment criterion relative to $(\mathbf{X,Y})$ in $\g[D]$ by Statement \ref{VAS} in Theorem~\ref{almost}. Since this holds for every DAG in $[\g]$, 
	%$\optb{\g}$ satisfies the generalized adjustment criterion relative to $(\mathbf{X,Y})$ in $\g$. Hence, by Theorem~\ref{A-thm:adjustment-set}, 
	$\optb{\g}$  is a valid adjustment set relative to $(\mathbf{X,Y})$ in $\g$. 
	
	Statement \ref{res2-cpdag} follows from the proof of Statement \ref{res1-cpdag}, and Statement \ref{optimal} in Theorem~\ref{almost}, while
	Statement \ref{res3-cpdag} follows from the proof of Statement \ref{res1-cpdag} and Statement \ref{smallest} in Theorem~\ref{almost}.
\end{proofof}

\begin{Lemma}
	Let $\mathbf{X}$ and $\mathbf{Y}$ be disjoint node sets in a \mpdag{} $\g$, such that $\g$ is amenable relative to $(\mathbf{X,Y})$. Then for every node $W \in \possde(\posscnb{\g},\g)$
%	\fb{\g} \setminus \mathbf{X}$ 
	there is a path $p_W = ( X = V_1, V_2, \dots ,V_k=W ), k > 1$ from $\mathbf{X}$ to $W$, such that 
	\begin{enumerate}[label=(\arabic*)]
		
		\item\label{l:sm-l1} $p_W$ is proper, and 
		\item\label{l:sm-l2} $p_W = X \rightarrow V_1 \rightarrow \dots \rightarrow W$. 
	\end{enumerate}
	Note that if $\mathbf{Y} \subseteq \possde(\mathbf{X},\g)$, this implies that $\mathbf{Y} \subseteq \de(\mathbf{X},\g)$ and $\mathbf{Y} \subseteq \cnb{\g}$.
	\label{lemma:additional-1}
\end{Lemma}

\begin{proofof}[Lemma~\ref{lemma:additional-1}]
	We prove the lemma separately for: 
	
	\begin{enumerate}[label=(\alph*)]
	\item\label{small-lemma1} $W$ $\in$ $\posscnb{\g}$, and  
	
	\item\label{small-lemma2} $W \in$  $\possde(\posscnb{\g},\g)$, $W$ $\notin$ $\posscnb{\g}$.
		\end{enumerate}
	
		\ref{small-lemma1} Since $W \in \posscnb{\g}$, $W$ lies on a proper possibly causal path $p$ from $X \in \mathbf{X}$ to $Y \in \mathbf{Y}$. Let $\overline{p(X,W)}$ be a subsequence of $p(X,W)$ as defined in Lemma~\ref{lemma:unshielded-analog}. Then $\overline{p(X,W)}$ satisfies \ref{l:sm-l1}. 
		Since $\overline{p(X,W)} \oplus p(W,Y)$ is a possibly causal path from $X$ to $Y$ that contains no non-endpoint node in $\mathbf{X}$ in $\g$, and since $\g$ is amenable relative to $(\mathbf{X,Y})$, $\overline{p(X,W)}$ starts with a directed edge out of $X$. Then since $\overline{p(X,W)}$ is an unshielded possibly causal path from $X$ to $W$ that starts with a directed edge out of $X$, $\overline{p(X,W)}$ must be a causal path from $X$ to $W$.
		
		\ref{small-lemma2}  Suppose $W \in \possde(\posscnb{\g},\g)$ and $W \notin \posscnb{\g}$. Let $V$ be a node in $\posscnb{\g}$ such that $W \in \possDe(V,\g)$.
		Since $V \in \posscnb{\g}$, by \ref{small-lemma1} there is a causal path $p$ from some $X \in \mathbf{X}$ to $V$ in $\g$ that is proper. Then $p$ is of the form $X \rightarrow D_1 \rightarrow \dots \rightarrow V $ with possibly $V=D_1$. Since $W \in \possDe(V,\g)$, we can choose an unshielded possibly causal path $p'$ from $V$ to $W$ in $ \g$ (Lemma~\ref{lemma:unshielded-analog}). 
		
		We now concatenate $p$ and $p'$. Hence, let $I$ be the node closest to $X$ on $p$ such that $I$ is also on $p'$.  
		Then possibly $I = D_1$. Additionally, $I \neq X$, { since otherwise $X \rightarrow D_1$ and $p(D_1,V) \oplus p'(V,X)$ is a possibly causal path from $D_1$ to $X$ in $\g$, which contradicts the definition of possibly causal paths in $\g$. } %  (Lemma~\ref{lemmamarlcycle}). 
		Hence, $q = p(X,I) \oplus p'(I,W)$ is a possibly causal path from $X$ to $W$ that starts with $X \rightarrow D_1$ in $\g$.
		
		Let $\overline{q} = ( X ,V_1, \dots V_{k}=W), k \ge 1$ be the subsequence of $q$ that forms an unshielded path (Lemma~\ref{lemma:unshielded-analog}). We next show that $\overline{q}$ satisfies \ref{l:sm-l2}.
		The edge between $X$ and $V_1$ is $X - V_1$, or $X \rightarrow V_1$. Since $V_1$ is also on $p$, $p(X,V_1)$ is of the form $X \rightarrow D_1 \cdots V_1$ (possibly $V_1=D_1$). Hence, $X - V_1$ cannot be in $\g$, as otherwise the existence of the edge $X \rightarrow D_1$ and $q(D_1,V_1) \oplus ( V_1,X)$ being a possibly causal path from $V_1$ to $X$ in $\g$ contradict the definition of possibly causal paths. 
		Hence, $\overline{q}$ is a possible causal path from $X$ to $W$, beginning with a directed edge out of $X$ and is hence causal.	%would contradict Lemma~\ref{lemmamarlcycle}. 
		
		Finally, $\overline{q}$ is either a proper causal path from $\mathbf{X}$  to $W$, in which case it satisfies \ref{l:sm-l1}, or there is a node $X' \in \mathbf{X}$ on $\overline{q}$ such that $\overline{q}(X',W)$ is.
%		of the form $X' \rightarrow \dots \rightarrow W$  and contains no non-endpoint node in $\mathbf{X}$. 

\end{proofof}

\begin{Lemma}
	Let $\mathbf{X}$ and $\mathbf{Y}$ be disjoint node sets in a \mpdag{} $\g$, such that $\g$ is amenable relative to $(\mathbf{X,Y})$, $\mathbf{Y} \subseteq \possde(\mathbf{X},\g)$ and let $\g[D]$ be a DAG in $[\g]$.
	Then  
	$\optb{\g[D]} = \optb{\g}$.
	\label{lemma:optimalSetEquiv}
\end{Lemma}

\begin{proofof}[Lemma~\ref{lemma:optimalSetEquiv}] %\color{blue}
	By definition, $\cnb{\g} \subseteq \cnb{\g[D]}$. Furthermore, it holds that $\pa(\cnb{\g},\g)$ $\subseteq$ $ \pa(\cnb{\g[D]},\g[D])$ and by Lemma~\ref{lemma:forbIsSame}, it also holds that $\fb{\g}$ $=$ $\fb{\g[D]}$. Since $\optb{\g}$ $=$ $\pa(\cnb{\g},\g)\setminus \fb{\g}$ and $\optb{\g[D]}$ $=$ $\pa(\cnb{\g[D]},\g[D])$ $\setminus$ $\fb{\g[D]}$, it thus follows that $\optb{\g}$ $\subseteq$ $\optb{\g[D]}$.
	
	We therefore only need to show $\optb{\g[D]} \subseteq \optb{\g}$. Hence, consider a node $A \in \optb{\g[D]}$. Since $A \in \optb{\g[D]}$, $A \notin \fb{\g[D]}$ and by Lemma~\ref{lemma:forbIsSame}, $A \notin \fb{\g}$. Thus, it is enough to show that $A \in \pa(\cnb{\g},\g)$.
	
	Since $A \in \pa(\cnb{\g[D]},\g[D])$, let $\mathbf{B}$ be the set of all nodes $B $, such that $B \in \posscnb{\g}$ and either $A \rightarrow B$ or $A-B$ is in $\g$. If $A-B$ is in $\g$ for any $B \in \mathbf{B}$, then since $A \in \possde(B,\g)$, $A \in \fb{\g}$, which is a contradiction. Hence, there exists a path $A \rightarrow B$ for every $B \in \mathbf{B}$ and $A \in \pa(\posscnb{\g},\g)$.
	Since $\mathbf{B} \subseteq \posscn{\g}$, for every node $B \in \mathbf{B}$, there is a possibly causal path from $\mathbf{B}$ to $\mathbf{Y}$ that does not contain a node in $\mathbf{X}$. Choose  the $B' \in \mathbf{B}$ that has a shortest such path among all nodes in $\mathbf{B}$. 
	We will show that $A \in \pa(\cnb{\g},\g)$ by showing that $B' \in \cnb{\g}$.
	
	If $B' \in \mathbf{Y}$,  it also holds that $B' \in \cnb{\g}$, by Lemma~\ref{lemma:additional-1} and hence we can assume $B' \notin \mathbf{Y}$. Since $B' \in \posscnb{\g}$ there exists a possibly causal path from $B'$ to $\mathbf{Y}$. Let $q$ be a shortest possibly causal path from $B'$ to $\mathbf{Y}$ that does not contain a node in $\mathbf{X}$ in $\g$. Then $q = ( B' =V_1, V_2,\dots, V_k=Y)$, $k >1$, $Y\in \mathbf{Y}$, with possibly $V_2 =Y$ and $q$ is an unshielded path. Further, $V_i \in \fb{\g}$, for every $i \in \{1,\dots,k\}$. 
	We will first show that $q$ is a causal path from $B'$ to $Y$.
	If $q$ starts with the edge $B' \rightarrow V_1$, then $q$ is causal from $B'$ to $Y$ in $\g$ (Lemma~\ref{lemma:unshielded-edges-analog}).  Otherwise, suppose $q$ starts with $B' - V_1$. Then $A \rightarrow B' -V_1$ is in $\g$, so $A -V_1$ or $A \rightarrow V_1$ is also $\g$ (by definition of \mpdag{}s in \citet{meek1995causal}, see Section~\ref{prelim}).
	However, $A - V_1$ implies $A \in \fb{\g}$ and $A \rightarrow V_1$, contradicts the choice of $B'$. 
	Since $B' \in \posscnb{\g}$, we can pick a proper causal path $p_{B'}$ from $\mathbf{X}$ to $B'$ in $\g$ (Lemma~\ref{lemma:additional-1}). 
%	Since $\g$ is acyclic, no node other than $B'$ is on both $p_{B'}$ and $q$. 
Then, $ p_{B'} \oplus q$ is a proper causal path from $\mathbf{X}$ to $\mathbf{Y}$ in $\g$ and hence, $B' \in \cnb{\g}$.
		
\end{proofof}

\begin{Lemma}
	Let $\mathbf{X}$ and $\mathbf{Y}$ be disjoint node sets in a \mpdag{} $\g$, such that $\g$ is amenable relative to $(\mathbf{X,Y})$ and let $\g[D]$ be a DAG in $[\g]$. Then
	$\fb{\g} = \fb{\g[D]}$.
	\label{lemma:forbIsSame}
\end{Lemma}

\begin{proofof}[Lemma~\ref{lemma:forbIsSame}]
	%\ema[Not done yet. This is the proof that requires the most work, but I did a version of it for the PAGs already, so it should be doable.]
	As by definition $\fb{\g}=$ $\possde(\posscnb{\g},\g)$ $\cup$ $\mathbf{X}$ and $\fb{\g[D]}= \de(\cnb{\g[D]},\g[D]) \cup \mathbf{X}$, it suffices to show that \[\possde(\posscnb{\g},\g)\subseteq \de(\cnb{\g[D]},\g[D]).\]
	%Since $\de(\cnb{\g},\g) \subseteq \possde(\posscnb{\g},\g)$ and $\de(\cnb{\g},\g) \subseteq \de(\cnb{\g[D]},\g[D])$, 
	It is enough to show that if $A \in  \possde(\posscnb{\g},\g)\setminus \de(\cnb{\g},\g)$, then $A$ $\in$ $\de(\cnb{\g[D]},\g[D])$.
	We divide the proof in two parts:
	\begin{enumerate}[label = (\alph*)]
		\item\label{forbissame1} $\posscnb{\g} \subseteq \de(\cnb{\g[D]},\g[D])$, and
		\item\label{forbissame2} $\possde(\posscnb{\g},\g) \subseteq \de(\cnb{\g[D]},\g[D])$.
	\end{enumerate}
	
	\ref{forbissame1} Let $A \in \posscnb{\g} \setminus \cnb{\g}$. We will show that $A \in  \de(\cnb{\g[D]},\g[D])$.
	Since $A$ is a possibly causal node relative to $(\mathbf{X,Y})$ in $\g$, there exists a proper causal path $p$ from $\mathbf{X}$ to $A$  by Lemma~\ref{lemma:additional-1}. 
%	Then $p$ is of the form $X \rightarrow \dots D \rightarrow A$, $X \in \mathbf{X}$ and possibly $X =D$. 
	Additionally, let $q$ be a shortest possibly causal path from $A$ to $\mathbf{Y}$ that does not contain a node in $\mathbf{X}$ in $\g$. Then $q$ is an unshielded possibly causal path from $A$ to $Y \in \mathbf{Y}$ and since $A \notin  \cnb{\g}$, $q$ is of the form $A - V_1 \dots V_k$, $V_k= Y$. Then $Y \in \possde(\posscnb{\g})$ and hence, by Lemma E.6, $Y \in \cnb{\g}$.
	
	Since $q$ is an unshielded possibly causal path from $A$ to $Y$ in $\g$, the corresponding path to $q$ in $\g[D]$ is of the form $A \rightarrow \dots \rightarrow Y$, or $A \leftarrow \dots \leftarrow Y$, or $A \leftarrow \dots \leftarrow V_j \rightarrow \dots \rightarrow Y$, for some $ 1 < j < k$. 
	
	If $q$ corresponds to  $A \rightarrow \dots \rightarrow Y$ in $\g[D]$, we can concatenate $r$ and $q$ in $\g[D]$, so that $A \in \de(\cnb{\g[D]},\g[D])$. If $q$ corresponds to $A \leftarrow \dots \leftarrow Y$ in $\g[D]$, then $A \in \de(Y, \g[D])$. By Lemma \ref{lemma:additional-1} and the fact that $Y \in \cnb{\g}$ it follows that $A \in \de(\cnb{\g[D]},\g[D])$ in this case. 
	
		Lastly, suppose $q$ corresponds to $A \leftarrow \dots \leftarrow V_j \rightarrow \dots \rightarrow Y$ in $\g[D]$.  We will show that $V_j \in \de(\cnb{\g[D]},\g[D])$, since then from $A \in \de(V_j, \g[D])$, it follows that $A \in  \de(\cnb{\g[D]},\g[D])$.  Since $V_j \in \possde(A,\g)$, $V_j \in \fb{\g}$ and since $V_j$ is on $q$, $V_j \neq X$. Hence, $V_j \in \fb{\g} \setminus \mathbf{X}$. Let $p_{V_{j}}$ be a proper causal path from $\mathbf{X}$ to $V_j$ in $\g$ (Lemma~\ref{lemma:additional-1}) and let ${p}_{V_{j}}^*$ and $q^*$ be the paths corresponding to $p_{V_{j}}$ and $q$ in $\g[D]$. Then we can concatenate ${p}_{V_{j}}^*$ and $q^{*}(V_j,Y)$ so that $V_j \in \de(\cnb{\g[D]},\g[D])$. %  Let $V_i$, $1 < i \le k$ be the node closest to $A$ on $q$ such that $q(V_i,Y)$ is a causal path from $V_i$ to $Y$, possibly $V_i = Y$. Then $j \le i$.  Since $X \rightarrow \dots \rightarrow D \rightarrow V_j$ is in $\g[D]$, by Lemma~\ref{lemmamarlcycle}, no node on $q(V_j,Y)$ is on $p(X,D)$. Hence, $p(X,D) \oplus \langle D, V_j \rangle \oplus q(V_j, Y)$ corresponds to a causal path from $X$ to $Y$ in $\g[D]$.}
		
		\ref{forbissame2} From \ref{forbissame1}, $\posscnb{\g} \subseteq \de(\cnb{\g[D]},\g[D])$. 
		Thus, by property of descendant sets, $\de(\posscnb{\g}, \g[D])$ $\subseteq$ $\de(\cnb{\g[D]},\g[D])$. Additionally, since 
		\[\de(\posscnb{\g}, \g) \subseteq \de(\posscnb{\g},\g[D]),\]
		 it follows that $\de(\posscnb{\g}, \g) \subseteq \de(\cnb{\g[D]},\g[D])$. 
		
		To finish the proof we only need to show 
		\[\possde(\posscnb{\g},\g) \subseteq \de(\posscnb{\g}, \g). \]
		 Hence, let $A \in \possde(\posscnb{\g},\g)$.
		Let $B$ be a node in $\posscnb{\g}$ such that $A \in \possde(B,\g)$. If $A \in \de(B,\g)$, then $A \in \de(\posscnb{\g}, \g)$. 
		Otherwise, let $r = ( B=V_1,\dots, V_k=A)$, $k >1$ be a shortest possibly causal path from $B$ to $A$. 
		Then $r$ is an unshielded path of the form $B - \dots -D_1 \rightarrow \dots \rightarrow A$, with possibly $D_1 = A$.
		It is only left to show that $D_1 \in \de(\posscnb{\g},\g)$, since then $A \in \de(\posscnb{\g},\g)$.
		
		Since $D_1 \in \possde(\posscnb{\g},\g)$ and by definition 
		\[\possde(\posscnb{\g},\g)= \fb{\g} \setminus \mathbf{X},\] 
		let $p_{D_1}$ be a proper causal path from $\mathbf{X}$ to $D_1$. In order to prove that $D_1$ is in the set $\de(\posscnb{\g},\g)$, we only need to show that there is a possibly causal path from $D_1$ to $Y$ that does not contain a node in $\mathbf{X}$.
		
		Hence, let $-r = ( A, \dots, D_1, \dots,B)$. Since $r$ is a possibly causal unshielded path from $B$ to $A$, $r$ and $-r$ are paths of definite status. Then since $-r(D_1,B)$ is of the form $D_1 - \dots - B$, $-r(D_1,B)$ is a possibly causal path from $D_1$ to $B$.
		Since $B \in \posscnb{\g}$, let $s$ be a possibly causal path from $B$ to $Y \in \mathbf{Y}$ that does not contain a node in $\mathbf{X}$. % We now concatenate $-r(D_1,B)$ and $s$. 
		Let $D$ be the node on $-r(D_1,B)$ closest to $D_1$ that is also on $s$. Then $q= -r(D_1,D) \oplus s(D,Y)$ is a possibly causal path from $D_1$ to $Y$ in $\g$. 
		
		Lastly, since $s$ does not contain a node in $\mathbf{X}$ any node in $\mathbf{X}$ on $q$ would need to be on $D_1 - \dots - D$. However, this would contradict the amenability of $\g$ relative to $(\mathbf{X,Y})$. Hence, $q$ does not contain a node in $\mathbf{X}$.
	
	%Since $V_i \in \possde(\posscnb{\g},\g)$ and $\possde(\posscnb{\g},\g)= \fb{\g} \setminus \mathbf{X}$, let $p$ be a proper causal path from $\mathbf{X}$ to $V_i$ (Lemma~\ref{lemma:additional-1}) such that $p =(X, \dots, V_i ), X \in \mathbf{X}$. We now concatenate $p$ and $s$. Hence, let $W$ be the closest node to $X$ on $p$ that is also on $s$. Then $p(X,W) \oplus s(W,Y)$ is a proper possibly causal path from $\mathbf{X}$ to $\mathbf{Y}$ in $\g$ and since $V_i \in \de(W,\g)$, $V_i \in \de(\posscnb{\g},\g)$.}
\end{proofof}

\begin{Lemma}(\citealp[Cf.\ Lemma~B.1 in][]{zhang2008completeness} and \citealp[Lemma~3.6 in][]{perkovic17})
	Let $X$ and $Y$ be two nodes in a \mpdag{} $\g$. If $p$ is a directed causal path from $X$ to $Y$ in $\g$, then a subsequence $\overline{p}$ of $p$ forms an unshielded directed causal path from $X$ to $Y$ in~$\g$.
	\label{lemma:unshielded-analog}
\end{Lemma}

\begin{Lemma}(\citealp[Lemma~7.2][]{maathuis2013generalized} and \citealp[Lemma~B.1][]{perkovic17})
	Let $p =(V_1, \dots , V_k)$ be a possibly causal definite status path in a \mpdag{} $\g$. If there is a node $V_i, i \in \{1, \dots, k-1\}$ such that $V_i \rightarrow V_{i+1}$, then $p(V_i, V_k)$ is a causal path in $\g$.
	\label{lemma:unshielded-edges-analog}
\end{Lemma}

\section{Supplement: Simulation study}

\subsection{Setup}
\label{setupMSE}

\noindent{\bf Technical details:} For our simulation we use {\em R (3.5.2)} and the R-package {\em pcalg (2.6-11)} by \citet{Kalisch2012}.

\noindent {\bf Graph:} We uniformly draw the number of nodes in the DAG from $\{10,20,50,100\}$, the expected neighborhood size from $\{2,3,4,5\}$ and the graph type to be either Erd\"os-R\'enyi or power law. We use the function {\tt randDAG} from the R-package {\em pcalg} to generate 10'000 such graphs.

\noindent {\bf Causal linear model:} To each graph $\g$ we associate a causal linear model in the following manner. We draw the model's error distribution be either normal, logistic, uniform or a t-distribution with 5 degrees of freedom, with probabilities $\frac{1}{2},\frac{1}{6},\frac{1}{6}$ and $\frac{1}{6}$. For each node, we then draw error parameters ensuring that the mean is 0 and the variance between $0.5$ and $1.5$. Specifically, for the normal distribution, the variance parameter is drawn uniformly from $[0.5,1.5]$ and the mean set to 0. For the logistic distribution, the scale parameter is drawn uniformly from $[0.4,0.7]$ and the location parameter is set to 0. For the uniform distribution, the sampling interval is of the form $[-a,a]$, with $a$ drawn uniformly from $[1.2,2.1]$. Finally, the t-distribution is first normalized and then multiplied with the square root of a parameter $v$ uniformly drawn from $[0.5,1.5]$. Finally, we draw coefficients for each edge in $\g$ from a uniform distribution on $[-2,-0.1] \cup [0.1,2]$.

\noindent {\bf The pair} $(\mathbf{X},Y)${\bf:} For each DAG $\g[D]$, we draw one pair $(\mathbf{X},Y)$ in the following manner. We first draw $|\mathbf{X}|$ from $\{1,2,3\}$, with probabilities $\frac{1}{2},\frac{1}{4}$ and $\frac{1}{4}$.  We then uniformly draw node sets $\mathbf{X}$ of the specified size, until $\bigcap_{X_i \in \mathbf{X}} \de(X_i,\g[D]) \neq \emptyset$ and then draw $Y$ uniformly from $\bigcap_{X_i \in \mathbf{X}} \de(X_i,\g[D])$. We then verify whether $\mathbf{X} \cap \de(\text{cn}(\mathbf{X},Y,\g[D]),\g[D]) = \emptyset$ and whether the true CPDAG is amenable with respect to $(\mathbf{X},Y)$.
%the being done using the ``dag2cpdag" function from ``pcalg". 
This is done to ensure that valid adjustment sets exist in both the true DAG and CPDAG. If either is not the case, we discard $(\mathbf{X},Y)$ and repeat the procedure. If no $(\mathbf{X},Y)$ is found, we draw a new DAG with the same parameters as the original. 

\noindent {\bf Total effect estimation:}  For each DAG $\g[D]$ we uniformly draw a sample size $n \in \{125,500,2000,10000\}$. We then sample 100 data sets of size $n$ from the causal linear model corresponding to $\g[D]$. For each data set we compute a graph estimate $\widehat{\g}$. 
When the error distribution is normal, a CPDAG is estimated with the Greedy Equivalence Search (GES) algorithm by \citet{Chickering02}. This is done with the {\tt ges} function from the R-package {\em pcalg} without tuning any of the parameters. When the error distribution is not-normal, we estimate a DAG with the Linear Non-Gaussian Acyclic Models (LiNGAM) algorithm by \citet{shimizu2006linear}. This is done with the {\tt lingam}  function from the R-package {\em pcalg}. 
%a DAG was estimated using ``lingam".  
%
%
%using the function ``ges" 
%
%When the error distribution was not-normal a DAG was estimated using ``lingam".  
%
%
%We then computed 
%
%the estimated causal graph $\widehat{\g}$ was then used to compute the sets $\text{O}(\mathbf{X},Y,\g[\widehat{G}]), \pa(\mathbf{X},\g[\widehat{G}]) \setminus \text{forb}(\mathbf{X},Y,\g[\widehat{G}])$ and $\text{adjust}(\mathbf{X},Y,\g[\widehat{G}])$.
%
%
%
%Only once all required conditions were fulfilled, did we proceed to compute the sets $\text{O}(\mathbf{X},Y,\g), \pa(\mathbf{X},\g) \setminus \text{forb}(\mathbf{X},Y,\g)$ and $\text{adjust}(\mathbf{X},Y,\g)$. 
%
%
%\noindent{\bf Total effect estimation:} 

We then proceed to estimate the total effect of $\mathbf{X}$ on $Y$ via covariate adjustment. We considere the graphical adjustment sets  $\text{Adjust}(\mathbf{X},Y,\g),\pa(\mathbf{X},\g) \setminus \textrm{forb}(\mathbf{X},Y,\g)$ and $\textbf{O}(\mathbf{X},Y,\g)$, computing them both from the true DAG $\g[D]$ and the graph estimate $\widehat{\g}$. Further, we considere the non-graphical empty set. The sets $\text{Adjust}(\mathbf{X},Y,\g)$ and $\textbf{O}(\mathbf{X},Y,\g)$ are guaranteed to be valid adjustment sets and the set $\pa(\mathbf{X},\g) \setminus \textrm{forb}(\mathbf{X},Y,\g)$ is guaranteed to be a valid adjustment set for the case $|\mathbf{X}=1|$, but not if $|\mathbf{X}|>1$. The empty set is generally not a valid adjustment set. 
We then estimate the total effect in the following manner: 
%Depending on whether the considered covariate set $\mathbf{Z}$ depended on $\g$ or $\g[\widehat{\g}]$, the total effect estimation differed slightly:
\begin{enumerate}
	\item When the considered adjustment set $\mathbf{Z}$ is computed from $\g[D]$,
	\[ 
	\hat{\tau}^{\mathbf{z}}_{y\mathbf{x}}=\boldsymbol{\hat{\beta}}_{y\mathbf{x}.\mathbf{z}}.
	\]
	\label{eq:DAG}
	\item When the considered adjustment set $\mathbf{Z}$ is computed from $\g[\widehat{\g}]$,
	\begin{numcases}{\hat{\tau}^{\mathbf{z}}_{y\mathbf{x}} =}
		\boldsymbol{\hat{\beta}}_{y\mathbf{x}.\mathbf{z}}, & \parbox{8cm}{if $\g[\widehat{\g}]$ amenable, $\mathbf{X} \cap \de(\text{cn}(\mathbf{X},Y,\g[\widehat{\g}]),\g[\widehat{\g}]) = \emptyset$ and 
		$Y \in \possde(\mathbf{X},\g[\widehat{\g}]),$} \label{fine} \\
		\boldsymbol{0}, & if $Y \notin \possde(\mathbf{X},\g[\widehat{\g}]),$ \label{noDesc} \\  
		\text{NA}, & else. \label{noVAS}
\end{numcases}
	\label{eq:CPDAG}
\end{enumerate}
For comparability, we also treat the non-graphical empty set as if it were graphical, coming from both $\g[D]$ and $\widehat{\g}$. Accordingly, we estimate total effects with both procedure (\ref{eq:DAG}) and (\ref{eq:CPDAG}). In total, we thus obtain 8 total effect estimates.

%To maintain comparability, we treated the two non-graphical sets $\emptyset$ and $\mathbf{V} \setminus (\mathbf{X} \cup \{Y\})$ as if they had been graphical, coming from both $\g$ and $\hat{\g}$. Accordingly, we estimated the total effects with procedures \ref{eq:DAG} and \ref{eq:CPDAG}. In total, we thus obtained 10 total effect estimates.
The difference between procedure (\ref{eq:DAG}) and (\ref{eq:CPDAG}) arises from the fact that $(\mathbf{X},Y)$ is sampled in a manner ensuring that the two cases \eqref{noDesc} and \eqref{noVAS} do not occur for the true DAG and the corresponding true CPDAG.  The decision to return $\boldsymbol{0}$ in \eqref{noDesc} is based on the fact that the total effect on a non-descendant is $\mathbf{0}$. Since it affects all estimates with respect to $\widehat{\g }$ equally, it has the effect of making their average output more similar. We chose to return ``NA" in \eqref{noVAS}, effectively discarding it, as in this case no valid adjustment set exists. When this occurs, we recommend the use of alternative total effect estimators such as the {\em IDA} algorithm by \citet{maathuis2009estimating} and the {\em jointIDA} algorithm by \citet{nandy2017estimating}.

%The pair $(\mathbf{X},\g)$ was drawn in a manner ensuring that 

%Whenever no valid adjustment set relative to $(\mathbf{X},Y)$ existed in $\hat{\g}$ we returned ``NA", essentially ignoring these cases. In setting such as these, we recommend alternative causal effect estimators such as instrumental variable estimators, or the jointIDA algorithm by \citet{nandy2017estimating}. Hence the decision to disregard them for our mean squared error computation. This is also the difference between the non-graphical sets $\emptyset$ and $\emptyset(\g[\widehat{\g}])$, as well as $\mathbf{V} \setminus (\mathbf{X} \cup \mathbf{Y})$ and $\mathbf{V} \setminus (\mathbf{X} \cup \mathbf{Y})(\g[\widehat{\g}])$. We also set $\hat{\tau}^{\mathbf{z}}_{y\mathbf{x}}= \boldsymbol{0}$, whenever When $Y \notin \possde(\mathbf{X},\g[\widehat{\g}])$, . We did so to reflect that the total effect (joint or not) of a random variable on a non-descendant is always $0$ and hence when a researcher encounters an estimated causal graph $\widehat{\g}$, such that $Y \notin \possde(\mathbf{X},\g[\widehat{\g}])$ setting $\hat{\tau}^{\mathbf{z}}_{y\mathbf{x}}= \boldsymbol{0}$ is the natural choice (see Remark \ref{A-remark: joint vs non}). 

\noindent{\bf Mean squared error computation:} 
For each graph $\g$, associated causal linear model and each of the corresponding 100 data sets, we compute 8 different estimates of $\tau_{y\mathbf{x}}$ as explained above. 
%We look at each $X$-component of our estimates individually. 
%$\hat{tau}_{y\mathbf{x}}
%Aggregating the results over all data sets yields a total of up to $100$ estimates for each method. 
We compute the empirical mean squared error with respect to the true total effect of each estimator over these 100 estimates, where we look at each $(\hat{\tau}_{y\mathbf{x}})_i$ for $X_i \in \mathbf{X}$ separately. 
Here, the true total effect of $\mathbf{X}$ on $Y$ is computed from the path coefficients of the causal linear model. 

%We use the term method rather than covariate set, as for the procedures depending on the estimated causal graph, the specific set will vary for different data sets.

\begin{figure}[p!]
	\centering
	\includegraphics[ height=17.4cm, width=13.5cm]{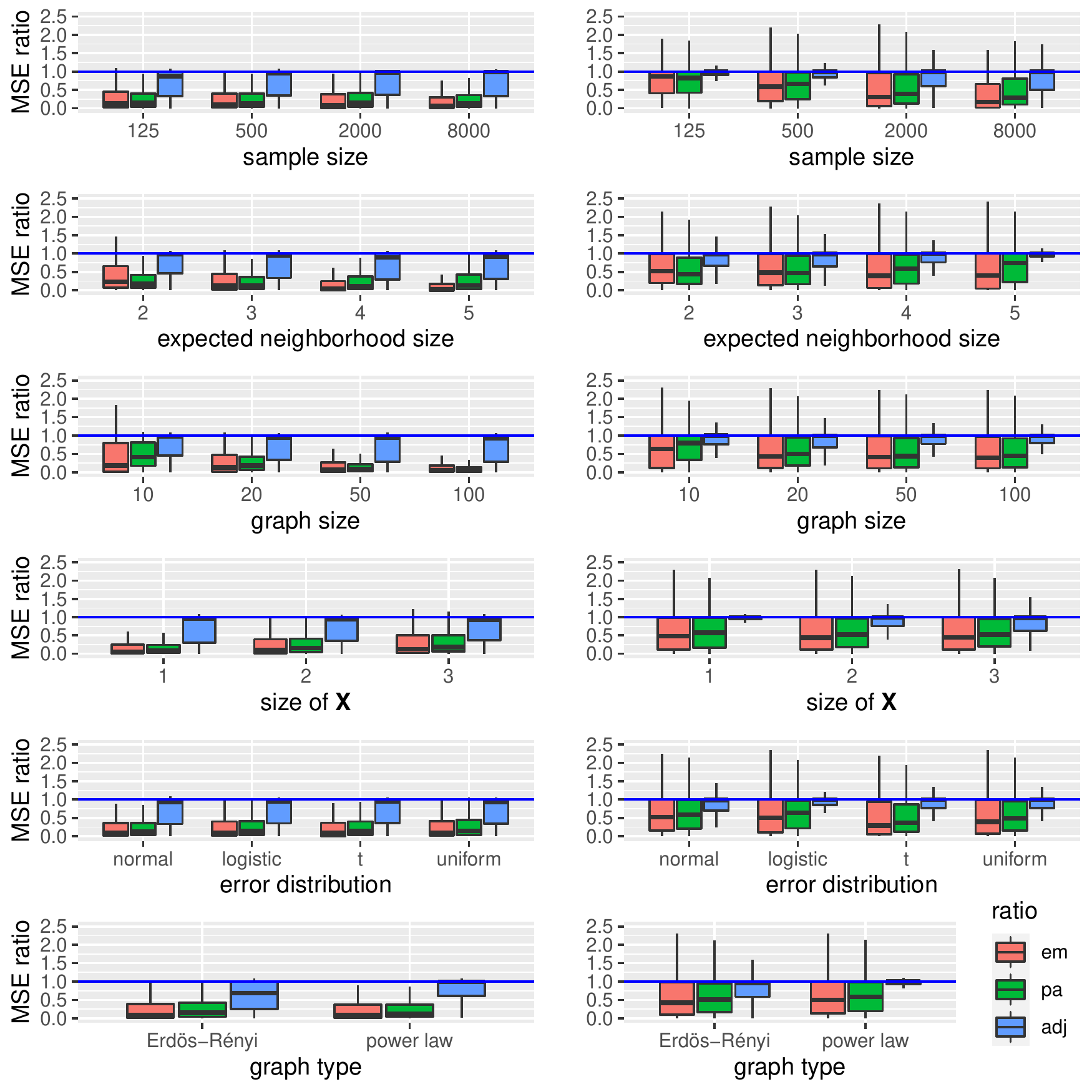}
	\caption{
		Boxplots of the ratios of the mean squared errors provided by $\text{\textbf{O}}(\mathbf{X},Y,\g)$ and the three alternative adjustment sets $\emptyset, \pa(\mathbf{X},\g) \setminus \textrm{forb}(\mathbf{X},Y,\g)$ and $\textrm{Adjust}(\mathbf{X},Y,\g)$, denoted respectively by ``em", ``pa" and ``adj", as a function of sample size, expected neighborhood size, graph size, size of $\mathbf{X}$, error distribution and graph type.
		The cases where the true causal DAG is used are given on the left and the ones where the causal graph is estimated on the right.}
	\label{moreBox}
\end{figure}

\begin{figure}[p!]
	\centering
	\captionsetup[subfigure]{labelformat=empty} 
	\subfloat[]{\includegraphics[height=7.7cm, width=7cm]{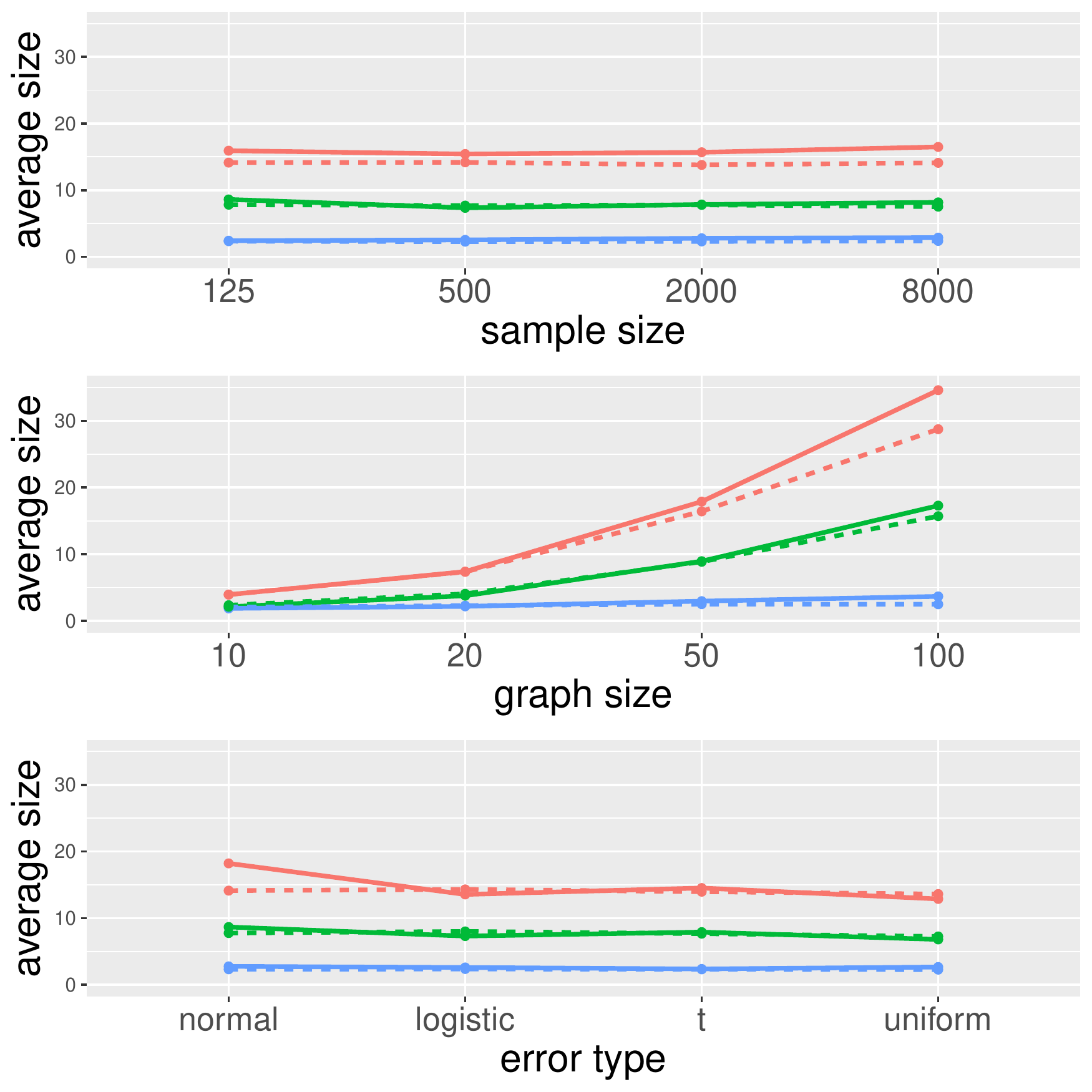}}
	\subfloat[]{\includegraphics[height=7.7cm, width=7cm]{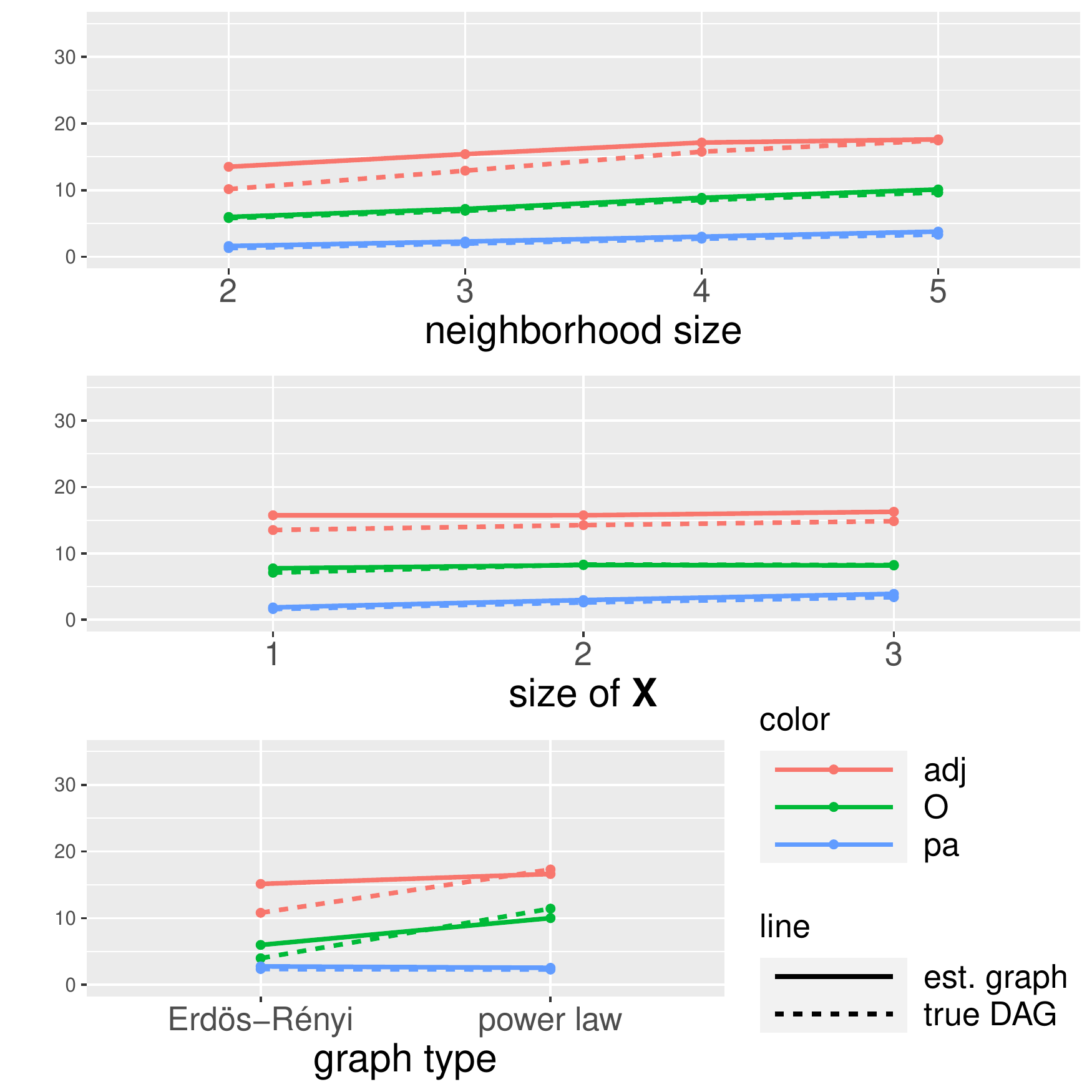}}
	\vspace{-18pt}
	\caption{The average size of $\pa(\mathbf{X},\g) \setminus \textrm{forb}(\mathbf{X},Y,\g), \text{\textbf{O}}(\mathbf{X},Y,\g)$ and $\textrm{Adjust}(\mathbf{X},Y,\g)$, denoted ``pa",``O" and ``adj" respectively,  in the true causal DAG and the estimated graphs as a function of sample size, expected neighborhood size, graph size, size of $\mathbf{X}$, error distribution and graph type.}
	\label{sizediff}
\end{figure}

\begin{figure}[p!]
	\centering
	\captionsetup[subfigure]{labelformat=empty}
	\subfloat[]{\includegraphics[height=7.7cm, width=7cm]{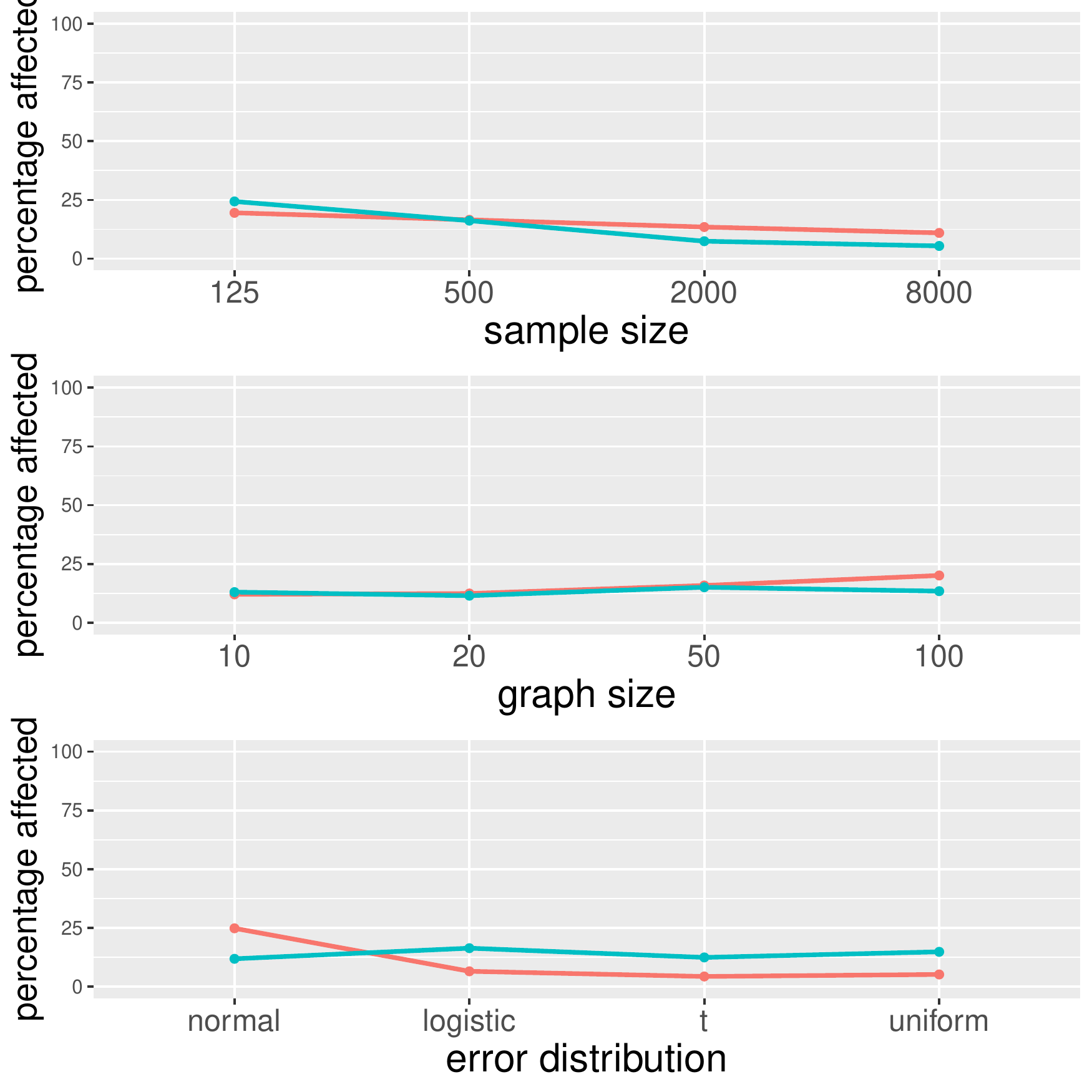}} 
	\subfloat[]{\includegraphics[height=7.7cm, width=7cm]{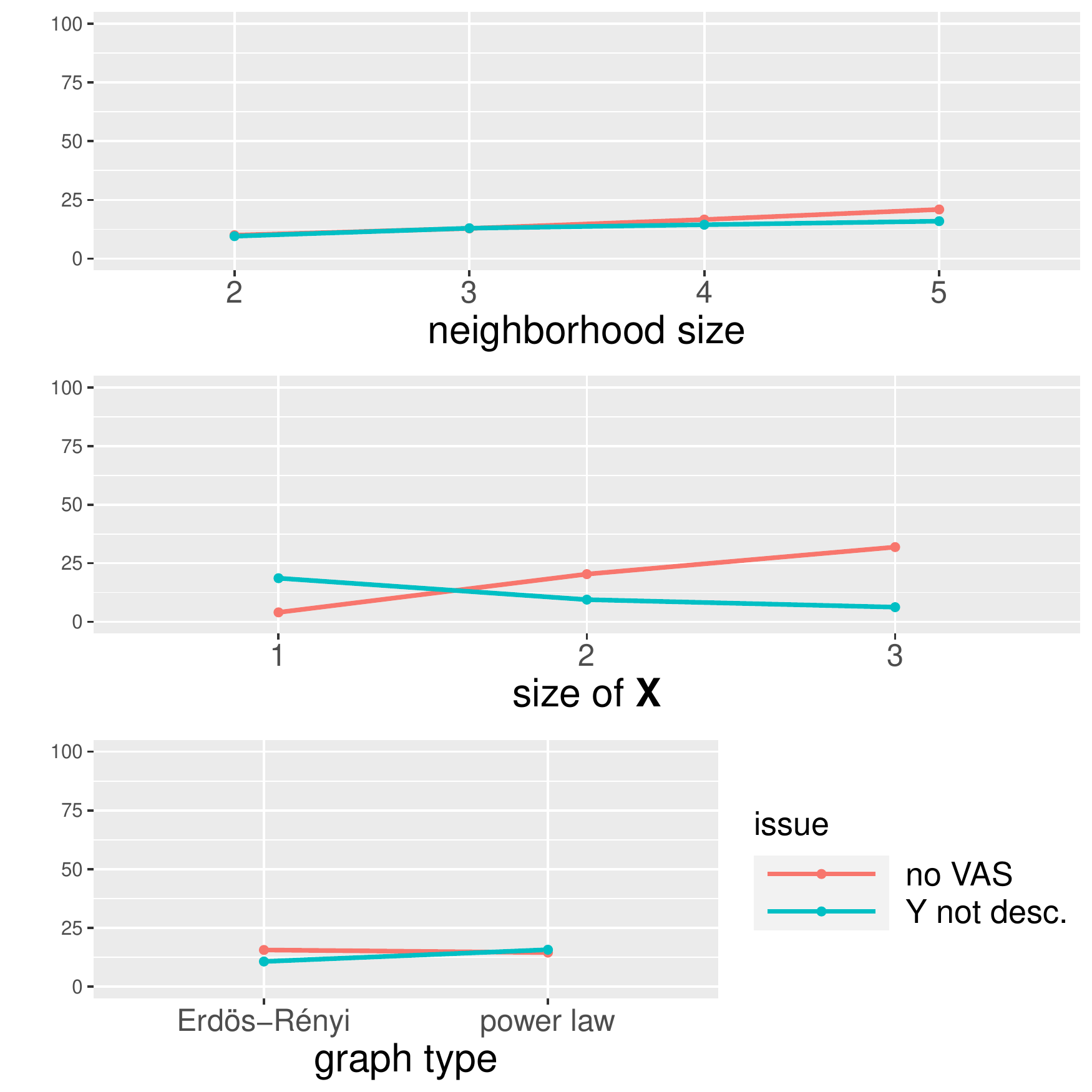}}
	\vspace{-18pt}
	\caption{The average percentage of estimated causal graphs $\widehat{\g}$ that do not have a valid adjustment set relative to $(\mathbf{X},Y)$, denoted ``no VAS" or where $Y \notin \possde(\mathbf{X},\widehat{\g})$, denoted ``Y not desc.", as a function of sample size, expected neighborhood size, graph size, size of $\mathbf{X}$, error distribution and graph type.}
	\label{issues}
\end{figure}

\subsection{Additional results}
\label{resultsMSE}

To understand how different settings impact the performance of each considered adjustment set,  Figure \ref{moreBox} shows boxplots of the MSE ratios as a function of sample size, expected neighborhood size, graph size, size of $\mathbf{X}$, error distribution and graph type. This plot reveals some interesting patterns.

When the true DAG is used, the ratios are generally stable across varying settings, with one major exception. The more complex the graph, i.e., the larger the graph and expected neighborhood size, the smaller the MSE ratios. This is probably due to the fact that all three alternative adjustment sets are less likely to be similar to $\text{\textbf{O}}(\mathbf{X},Y,\g[D])$ for larger or denser DAGs. Since $\text{\textbf{O}}(\mathbf{X},Y,\g[D])$ is guaranteed to be unbiased and to provide the smallest asymptotic variance, this leads to smaller ratios. 

When the graph is estimated some of these effects disappear, since graph estimation is more challenging for larger and denser graphs, especially when the sample size is small. This indicates that $\text{\textbf{O}}(\mathbf{X},Y,\widehat{\g})$ is more affected by graph estimation errors than the alternative adjustment sets. Even with these difficulties there are few especially large ratios, while, except for the comparison with $\textrm{Adjust}(\mathbf{X},Y,\widehat{\g})$, there is a respectable number of ratios smaller than 0.5 in all settings. 

%In the setting where the causal graph was estimated the picture is more complex. Larger sample sizes lead to smaller ratios. This is not surprising as a larger sample size leads to a better graph estimate. Error prone graph estimates add essentially random bias terms to all estimators pushing the ratios closer to 1.
%Still, this does also indicate that the optimal set is possibly more sensitive to graph estimation errors than the alternative adjustment sets. Considering this, it is unsurprising that larger average neighborhood sizes do not lead to smaller ratios in this setting. Typically causal discovery algorithms only perform well for sparse graphs. Larger graph sizes, on the other hand, do lead to smaller ratios. More complex graphs and less similarity between the optimal set and the alternative sets seem to only lead to smaller ratios, if that complexity does not adversely impact the graph estimation.  

The only set competitive with the optimal set is $\textrm{Adjust}(\mathbf{X},Y,\g)$. However, the set  $\textrm{Adjust}(\mathbf{X},Y,\g)$ is also by construction a superset of the optimal set.  This is reflected in the average sizes plotted in Figure \ref{sizediff}. The average size of $\textrm{Adjust}(\mathbf{X},Y,\g)$ is nearly twice the size of $\text{\textbf{O}}(\mathbf{X},Y,\g)$. Similarly, $\text{\textbf{O}}(\mathbf{X},Y,\g)$ is on average larger than both $\pa(\mathbf{X},\g) \setminus \fb{\g}$ and of course also the empty set. However, while there is a gain in moving from the empty set to $\pa(\mathbf{X},\g) \setminus \fb{\g}$ and from there  to $\text{\textbf{O}}(\mathbf{X},Y,\g)$, there is no corresponding gain in any setting from the increased size of $\textrm{Adjust}(\mathbf{X},Y,\g)$ compared to $\text{\textbf{O}}(\mathbf{X},Y,\g)$.

%Overall, these results are encouraging, as the optimal set performs better than all 4 alternative sets considered here. However, bad graph estimates do adversely affect the optimal set's performance. A sufficient amount of data has to be collected, in order to reliably estimate the underlying causal graph, independently of which adjustment set is used.

\subsection{Issues related to graph estimation}

\label{graphIssues}

In the course of the mean squared error simulations we had to estimate graphs. As explained in the setup, there are two issues that may arise. 
Firstly, there may not be a valid adjustment set relative to $(\mathbf{X},Y)$ in the estimated causal graph $\widehat{\g}$, in which case we return ``NA" for all adjustment sets. Secondly, it can happen that $Y \notin \possde(\mathbf{X},\widehat{\g})$, in which case we return $\boldsymbol{0}$ for all adjustment sets. 
%If this occurs the natural``estimate" of the total effect $\boldsymbol{0}$ may be very biased. 
%In either case, covariate adjustment will not yield a good total effect estimate irrespective of the valid adjustment set used. 

Figure \ref{issues} shows the average percentage of estimated graphs affected by either of the two issues. We see that they occur mostly for i) low sample sizes, ii) normal errors, as we can only estimate a CPDAG in this case and iii) in the case $|\mathbf{X}|>1$,  as only in this case $\mathbf{X} \cap \fb{\g} \neq \emptyset$ may occur (see Corollary \ref{cor:noforbx}). 

%These issues represent a significant obstacle to total effect estimation via adjustment, illustrating the importance of causal discovery for total effect estimation. 

%This also indicates that in practice complementary methods for total effect estimation are often required. We point out that the pair $(\mathbf{X},Y)$ was sampled in a manner ensuring that neither of these two problems arise for either the true causal DAG or the causal CPDAG.

%Note that these issues impact all adjustment sets equally.

%\include{parts/Rcode}

%\bibliography{parts/mybib}
%\bibliographystyle{rss}
%
%\end{document}

\end{document}